\documentclass[12pt]{amsart}
\usepackage[top=72pt,bottom=96pt,left=72pt,right=72pt]{geometry}
\usepackage{amsmath}
\usepackage{mathtools}
\usepackage{amssymb}
\usepackage{bbm}
\usepackage{extarrows}
\def\A{\mathrm{Ass}}
\def\C{\mathbb{C}}
\def\R{\mathbb{R}}
\def\E{\mathbb{E}}
\def\S{\mathbb{S}}
\def\GTS{\mathbf{qGTS}}
\def\F{\mathfrak{F}}
\def\VB{\mathbf{qVB}}

\def\PB{\mathbf{qPB}}
\def\Hor{\mathrm{Hor}}
\def\id{\mathrm{id}}
\def\Ker{\mathrm{Ker}}
\def\Im{\mathrm{Im}}

\def\class{\mathrm{class}}
\def\ad{\mathrm{ad}}
\def\H{\mathrm{H}}
\def\Mor{\textsc{Mor}}
\def\N{\mathbb{N}}
\def\Ad{\mathrm{Ad}}
\def\ad{\mathrm{ad}}
\def\NT{\textsc{nt}}
\def\Obj{\textsc{Obj}}

\def\Z{\mathbb{Z}}
\def\Rep{\mathbf{qRep}}
\def\triv{\mathrm{triv}}

\def\G{\mathcal{G}}

\def\Vert{\mathrm{Ver}}
\def\V{\mathrm{V}}
\def\u{\mathbbm{u}}
\def\U{\mathcal{U}}
\def\SU{\mathcal{SU}}
\def\torus{\mathbb{T}}
\def\G{\mathcal{G}}
\def\T{\mathcal{T}}

\theoremstyle{definition}
\newtheorem{Definition}{Definition}[section]

\theoremstyle{definition}
\newtheorem{Remark}[Definition]{Remark}

\theoremstyle{plain}
\newtheorem{Theorem}[Definition]{Theorem}

\theoremstyle{plain}
\newtheorem{Proposition}[Definition]{Proposition}

\numberwithin{equation}{section}

\begin{document}

\title{Functoriality of Quantum Principal Bundles and Quantum Connections}
\author{Gustavo Amilcar Salda\~na Moncada}
\address{Gustavo Amilcar Salda\~na Moncada\\
Mathematics Research Center, CIMAT}
\email{gamilcar@ciencias.unam.mx}
\begin{abstract}
Within the framework of Category Theory, we study the association between finite--dimensional representations of a compact quantum group and quantum vector bundles endowed with quantum linear connections, for a fixed quantum principal bundle equipped with a regular quantum principal connection. In particular, we establish a categorical equivalence between such \emph{quantum association functors} and quantum principal bundles with a particular regular and multiplicative quantum principal connection.

 \begin{center}
  \parbox{300pt}{\textit{MSC 2010:}\ 46L87, 58B99.}
  \\[5pt]
  \parbox{300pt}{\textit{Keywords:}\ Quantum Connections, Hermitian Structures, Quantum Gauge Group.}
 \end{center}
\end{abstract}
\maketitle

\section{Introduction}
In Differential Geometry, the study of principal bundles and principal connections is a fundamental topic. A key result in this theory is that, given a smooth principal $G$--bundle $\pi:GM\longrightarrow M$ over a manifold $M$, one can associate a fiber bundle over $M$ to every smooth manifold equipped with a smooth $G$--action. This construction defines a covariant functor (usually called the \emph{association functor}) between the category of manifolds with $G$--actions and the category of fiber bundles over $M$. In Physics, this framework provides the basis for the development of Yang--Mills models and field theory~\cite{diff}.

In \cite{2}, M. Nori shows a characterization of these functors; however, this paper was written in the framework of Algebraic Geometry. In \cite{saldgreg}, the authors present a kind of generalization of \cite{2} in the framework of Differential Geometry by considering principal connections and induced connections as well. In particular, it is shown that every covariant functor between these categories that satisfy certain properties is naturally isomorphic to an association functor for a unique (up to isomorphisms) principal bundle $GM$ over $M$ with a principal connection $\omega$. This implies a categorical equivalence between the category of principal bundles over $M$ with principal connections and the category of gauge theory sectors over $M$, category whose objects are these association functors \cite{saldgreg}.

The theory of quantum principal bundles (qpb's) was developed in order to generalize principal bundles into the framework of Non--Commutative Geometry. Many authors have contributed to this theory with their own developing, for example \cite{4,5,6,micho1,micho2,micho3}. All these formulations are intrinsically related by the theory of {\it Hopf--Galois extensions} \cite{8,9}. A lot of work has been developed in this framework, for example \cite{10,11,12}. In particular, M. Durdevich in \cite{micho1,micho2,micho3} used the notion of quantum group  $\G$ published by S. L. Woronowicz in \cite{woro1}  to play the role of the structure group of the qpb, but with another differential calculus that the one showed in \cite{woro2}. With all this, there is a natural definition of quantum principal connections (qpc's). This theory was extended later in order to {\it embrace} other \emph{classical} notions of principal bundles, for example characteristic classes \cite{15,16}.

 On the other hand, the Serre--Swan theorem gives us a natural way to generalize  the concept of finite--dimensional vector bundle to Non--Commutative Geometry: a finitely generated projective module; however, it is possible to use left modules, right modules, or even bimodules. A. Connes, Dubois--Violette, and others have studied in a profound way the concept of quantum vector bundles (qvb's) and quantum linear connections (qlc's) \cite{con,18}. 

This paper arises from the following question: \emph{is it possible, in some way, to recover the categorical equivalence of reference \cite{saldgreg} in the non--commutative geometrical setting?} More concretely, the aim of this paper is to establish a non--commutative geometrical counterpart, for compact spaces, of the categorical equivalence between principal bundles with principal connections and gauge theory sectors described in reference \cite{saldgreg}. To achieve this, we rely on Durdevich's theory, Woronowicz's representation theory, and Dubois--Violette's framework. 

This paper is not the only one to study quantum principal bundles from a categorical point of view. For example, the reader may consult \cite{11,michokrein} and Section 5 of reference \cite{libro}. The difference between our work and these references lies, of course, in the categorical equivalence presented in this paper for quantum bundles and quantum connections.

It is worth noticing that, even though the underlying context of this paper is that of $C^\ast$--algebras (see Theorem  \ref{theoremrep} and Remark \ref{rema}), the philosophy of this work is to develop the theory in purely algebraic--geometric terms, for quantum principal bundles, their differential calculus, and their quantum principal connections. For this reason, we have chosen to use Durdevich's formulation of quantum principal bundles. More concretely, Durdevich's theory allows one to obtain differential calculus and quantum connections in purely \emph{functorial} terms, as the reader can see in Sections 3.5, 4.2.

The approach presented in this paper is important because it provides a stronger foundation for the general theory and unifies three lines of research, naturally leading to the study of Gauge Theory within the framework of quantum principal bundles, as illustrated in references \cite{sald2,saldym,saldcon,saldhopf,saldpoints}. Furthermore, working with categories and functors inherently involves natural constructions that promote a common language. For example, the fact that we can recover the \emph{classical} categorical equivalence suggests that we are dealing with appropriate definitions of qvb, qlc, qpb, and qpc, among other geometric concepts.

This paper consists of $5$ sections. Following this introduction, Section~2 is divided into five subsections. 

The first subsection introduces the category $\Rep_H$ of finite--dimensional corepresentations (or quantum representations) of a $\ast$--Hopf algebra $H$. The second subsection is devoted to the category $\VB_{B}$ of quantum vector bundles over a fixed quantum space $B$, and the third subsection to the category $\PB_{B}$ of quantum principal bundles over $B$. In Section 2.4, we present the definition of associated quantum vector bundles in Durdevich’s formulation, and in Section 2.5 we introduce the quantum association functor $$\A_\zeta: \Rep_H\longrightarrow \VB_{B},$$ for a given quantum principal bundle $\zeta$. It is worth emphasizing that, in this section, we do not consider any differential structure or quantum connections on the spaces.

In Section~3, we introduce differential structures on the spaces; this section is also divided into five subsections. In the first subsection, we introduce the notion of the universal differential envelope $\ast$--calculus of a $\ast$--Hopf algebra $H$ (see~\cite{micho1,stheve}), which plays the role of differential forms on $H$. In Subsection~3.2, we introduce the notion of differential calculus on a quantum vector bundle, as well as the notion of quantum linear connections $\nabla$, thereby defining the category $\VB^\nabla_{\Omega^\bullet(B)}$. In the next subsection, we study the notion of differential calculus on a quantum principal bundle, together with the concepts of quantum principal connections, their covariant derivatives, and their curvatures. In particular, we introduce the notions of regular and multiplicative quantum principal connections.

In order to obtain the desired categorical equivalence, in Subsection~3.4 we impose a condition on the quantum base space $B$ of every quantum principal bundle. More concretely, we require that $B$ be a $\ast$--algebra stable under holomorphic calculus. In addition, in this subsection we study the implications of this hypothesis in terms of strong connections in the formulation of quantum principal bundles presented in reference \cite{libro}.

On the other hand and also in order to obtain the desired categorical equivalence, in Subsection~3.5 we present a method to construct differential calculus on quantum principal bundles from certain \emph{functorially reproducible} data. With this, the category $\PB^{\omega^c}_{\Omega^\bullet(B)}$ of quantum principal bundles with the regular and multiplicative quantum principal connection $\omega^c$ is defined. This subsection is based on Section~6.5 of~\cite{micho2}.

Building on Sections~2 and~3, in Section~4 we present the quantum association functor $$\A^{\omega^c}_{\zeta}:\Rep_H\longrightarrow \VB^{\nabla}_{\Omega^\bullet(B)}$$ for a quantum principal bundle $\zeta$ over $B$ equipped with $\omega^c$. In Subsection~4.1, we discuss some properties of this functor, while in Subsection~4.2 we establish the categorical equivalence, which is the main objective of this paper.

In particular, we emphasize Theorem~\ref{teo1}, which provides a method to reconstruct the quantum principal bundle $\zeta$ and the quantum principal connection $\omega^c$ from a contravariant bar functor between $\Rep_H$ and $\VB^{\nabla}_{\Omega^\bullet(B)}$. Moreover, Theorem~\ref{teo2} establishes the categorical equivalence between $\PB^{\omega^c}_{\Omega^\bullet(B)}$ and the category $\GTS^{\nabla}_{\Omega^\bullet(B)}$ of quantum gauge theory sectors over $B$. This latter category consists of all contravariant bar functors between $\Rep_H$ and $\VB^{\nabla}_{\Omega^\bullet(B)}$.

To conclude the paper, in Section~5 we present two classes of examples illustrating our theory. In Section~5.1, we show that the theory developed in this paper applies to the \emph{usual} quantum principal $U(1)$--bundle over the non--commutative $n$--torus, while in Section~5.2 we show that it also applies to homogeneous quantum principal bundles. Moreover, in Appendix~A we provide a brief summary of reference \cite{saldgreg}; so the reader can compare the results of this papers with the ones in Differential Geometry. In addition, in Appendix~B we recall some definitions from category theory needed for our purposes. Finally, in Appendix~C we present a $C^\ast$--algebraic generalization of the theory developed in this paper (at degree~0 and without quantum connections, as in Section~2.5), thereby establishing a link with previous works.

Before continue, let us talk about the notation in this paper. First, we will formally represent quantum spaces as an associative unital $\ast$--algebras $$(X,m,\mathbbm{1},\ast)$$ over $\C$,  where $m:X\otimes X\longrightarrow X$ is the product, $\mathbbm{1}$ is the unit of the algebra and  $\ast:X\longrightarrow X$ is the antilinear involution.  In general, we will omit the words {\it associative} and {\it unital}. Besides, all our $\ast$--algebra morphisms will be unital. Furthermore, when we work with {\it quantum structures} we will point out how we are going to denote them. 

In the whole text we will use Sweedler's notation and given an arbitrary category $\mathbf{C}$, the class of objects of $\mathbf{C}$ will be denoted by $\Obj(\mathbf{C})$  and given $c_1,\, c_2 \,\in \,  \Obj(\mathbf{C})$, we are going to denote by $\Mor(c_1,c_2)$ the class of all morphisms in $\mathbf{C}$ between $c_1$ and $c_2$. 

It is worth noting that, in this paper, we do not use the \emph{traditional} notation for the categories of finite--dimensional corepresentations and quantum vector bundles (finitely generated projective bimodules) commonly found in the literature (see, for example,~\cite{libro}). Instead, we work with the notation $\Rep_H$ and $\VB_B$. This choice is motivated by the following considerations:
\begin{enumerate}
\item In Differential Geometry, the category of finite--dimensional linear representations of a Lie group $G$ is usually denoted by $\mathbf{Rep}_G$. Likewise, the category of vector bundles over a fixed manifold $M$ is typically denoted by $\mathbf{VB}_M$. Hence, the notation $\Rep_H$ and $\VB_B$ is adopted here to \emph{emphasize} that the results of this paper can be viewed as non--commutative geometrical counterparts, \emph{at purely algebraic level}, of those in reference \cite{saldgreg} for compact manifolds, which reflects the underlying philosophy of the work. The same reasoning applies to the categories $\mathbf{PB}^\omega_M$ and $\mathbf{GTS}^\nabla_M$ (see Appendix~A), as well as to $\PB_B$, $\PB^\omega_{\Omega^\bullet(B)}$, and $\GTS^\nabla_{\Omega^\bullet(B)}$ (see Section~4.2).
\item The notation $\VB_B$ for the category of quantum vector bundles can be naturally extended to $\VB^{\nabla}_{\Omega^\bullet(B)}$, the category of quantum vector bundles endowed with quantum linear connections over a graded differential $\ast$--algebra $\Omega^\bullet(B)$ with $\Omega^0(B)=B$. A similar extension applies to $\PB_B$ and $\PB^{\omega^c}_{\Omega^\bullet(B)}$.
\end{enumerate}
Thus, we expect that these notational choices do not pose a significant difficulty for the reader.

\section{The Quantum Association Functor without Quantum Connections}

As mentioned in the introductory section, the desired quantum association functor studied in this paper has already been considered in several previous publications; see, for example, \cite{11,michokrein,libro}. The purpose of this section is to provide a reasonably self--contained summary of the quantum association functor without quantum connections (based on these  references); while introducing the basic concepts used throughout the paper and the definitions specific to Durdevich's formulation of quantum principal bundles.

\subsection{$\ast$--Hopf Algebras and Corepresentations}

As in Differential Geometry, the concept of a Lie group is fundamental for the study of principal bundles. The same holds in the \emph{quantum} setting. In this way, this subsection is primarily based on \cite{librogroups}, which provides a \emph{modern} approach to the study of compact quantum groups. However, more precisely, our interest lies not in compact quantum groups themselves, but in their canonically associated $\ast$--Hopf algebras.

\begin{Definition}
    \label{hopf}
    We say that a $\ast$--algebra $(H,m,\mathbbm{1},\ast)$ is a $\ast$--Hopf algebra if there exist $\ast$--algebra morphisms $$\Delta:H\longrightarrow H\otimes H,\qquad \epsilon:H \longrightarrow \C $$ called {\it the coproduct} and {\it the counit}, respectively, and there exists a linear map
    $$S:H\longrightarrow H$$ called {\it the antipode}
    such that 
\begin{equation*}
(\id_{H}\otimes \Delta)\circ  \Delta =(\Delta \otimes \id_H)\circ \Delta,\quad  (\epsilon \otimes \id_H)\circ \Delta = \id_H= (\id_H\otimes \epsilon)\circ \Delta, 
\end{equation*}
\begin{equation*}
m\circ (S\otimes \id_H)\circ \Delta=\eta \circ \epsilon \;\;\;\;\; \mbox{ and } \;\;\;\;\;m\circ (\id_H\otimes S )\circ \Delta=\eta \circ \epsilon,
\end{equation*}
where $\eta: \C\longrightarrow H$ is the linear map defined by $\lambda\longmapsto \lambda \mathbbm{1}$. A $\ast$--Hopf algebra will be represented by $H^\infty=(H,m,\mathbbm{1},\Delta,\epsilon,S,\ast)$.
\end{Definition}

The next step is to define the notion of \emph{quantum $H$--representations} \cite{librogroups}.

\begin{Definition}
\label{corep}
    Let $H^\infty$ be a $\ast$--Hopf algebra and let $V$ be a $\C$--vector space. A right $H$--corepresentation (or a right $H$--coaction or a quantum $H$--representation) on $V$ is a linear map
    \begin{equation}
        \label{ec.2.-1}
        \delta^V: V\longrightarrow V\otimes H
    \end{equation}
    such that $$ (\delta^V\otimes \id_H)\circ \delta^V=(\id_V\otimes \Delta)\circ \delta^V, \qquad (\id_V\otimes \epsilon)\circ \delta^V=\id_V.$$ We say that $\delta^V$ is finite--dimensional if $n_V:=\mathrm{dim}(V)\in \N_0:=\N\cup \{ 0\}$. 
    
    Equivalently, a right $H$--corepresentation on $V$ is an invertible element 
    \begin{equation}
        \label{ec.2.0}
        \delta^V\, \in\, B(V)\otimes H,
    \end{equation}
    with $B(V)$ the space of all linear endomorphisms of $V$, such that $$(\mathbbm{1}\otimes \Delta)\delta^V=\delta^V_{12}\,\delta^V_{13} \;\in\; B(V)\otimes H\otimes H,$$ where we have used leg--numbering notation \cite{librogroups}. 
    
    The space of all  finite--dimensional $H$--corepresentations will be denoted by  $$\Obj(\Rep_H).$$  
\end{Definition}

Also, we have

\begin{Definition}
\label{corepmor}
    Given two $H$--corepresentations $\delta^V$, $\delta^W$, a corepresentation morphism between $\delta^V$ and  $\delta^W$ is a linear map $$f:V\longrightarrow W$$ such that $$\delta^W \circ f=(f\otimes \id_H)\circ \delta^V.$$ The set of all  $H$--corepresentation morphisms between $\delta^V$ and  $\delta^W$ will be denoted by  $$\Mor(\delta^V,\delta^W).$$ 
\end{Definition}

\begin{Remark}
    \label{remacorep}
     Finite--dimensional $H$--corepresentation form a  category, denoted by $$\Rep_H.$$ The notions of monomorphism, epimorphism, and isomorphism of $H$--corepresentations should be clear.
\end{Remark}

 Once $V$ is equipped with an inner product $\langle - | - \rangle$, we say that a finite--dimensional $H$--corepresentation $\delta^V$ is \emph{unitary} if $\delta^V$, regarded as an element of $B(V)\otimes H$, is unitary.

According to~\cite{librogroups}, for every finite--dimensional $H$--corepresentation $\delta^V$ there exists an inner product on $V$ (not necessarily unique) with respect to which $\delta^V$ becomes unitary. Therefore, from now on, we will assume that every element of $\Obj(\Rep_H)$ is unitary.

 Additionally, we say that $\delta^V$ is {\it reducible} if there exists a non--trivial subspace $L$ $(L\not=\{0\}, V)$ such that $\delta^V(L)\subseteq L\otimes H$. Of course, $\delta^V$ is irreducible if it is not reducible.

Let $\Obj(\Rep_H)^{\mathrm{Irr}}$ be the set of all equivalence classes\footnote{With respect to the equivalence relation: two irreducible  $H$--corepresentations are related if and only they are isomorphic.} of irreducible (necessarily finite--dimensional) $H$--corepresentations. For each $[\delta^W] \in \Obj(\Rep_H)^{\mathrm{Irr}}$, choose an element $\delta^V \in [\delta^W]$ and let
\begin{equation}
    \label{ec.irre}
    \T
\end{equation}
be the set of all such chosen $H$--corepresentations, one for each equivalence class $[\delta^W]$. In the rest of the text, we will refer to $\T$ as a complete set of mutually non--equivalent, irreducible $H$--corepresentations and we will always assume that $\delta^\C_\triv$ $\in$ $\T$, where 
$$\delta^\C_\triv:\C\longrightarrow \C\otimes H,\qquad z\longmapsto z\otimes \mathbbm{1}.$$ This $H$--corepresentation is called the \emph{trivial} $H$--corepresentation on $\C$.

Let $\delta^V$ $\in$ $\T$ and fix an orthonormal linear basis $\{e_i \}^{n_V}_{i=1}$ (with respect to the inner product for which $\delta^V$ is unitary). Then
\begin{equation}
    \label{ec.2.3}
    \delta^V(e_j)=\sum^{n_V}_{i=1} e_i\otimes g^V_{ij}\; \in\; V\otimes H.
\end{equation}
The elements $\{g^V_{ij} \}^{n_V}_{i,j=1}$ are called \emph{matrix coefficients} of $\delta^V$ (with respect to $\{e_i \}^{n_V}_{i=1}$) and satisfy
\begin{equation}
\label{2.5properties}
\begin{aligned}
    \Delta(g^{V}_{ij})= \sum^{n_{V}}_{k=1} &g^{V}_{ik}\otimes g^{V}_{kj},\qquad S(g^{V}_{ij})=g^{V\,\ast}_{ji},\qquad \epsilon(g^{V}_{ij})=\delta_{ij}
    \\ &\sum^{n_{V}}_{k=1}S(g^{V}_{ik})\,g^{V}_{kj}=\sum^{n_{V}}_{k=1}g^{V}_{ik}\,S(g^{V}_{kj})=\delta_{ij}\mathbbm{1}
\end{aligned}
\end{equation}
with $\delta_{ij}$ being the Kronecker delta, among other properties \cite{woro1}.

\begin{Remark}
    \label{hopf1}
    In the rest of this work, we will assume that $\{\{g^V_{ij} \}^{n_V}_{i,j=1}\mid \delta^V \in \T\}$ is a linear basis of $H$.
\end{Remark}
Intuitively, the hypothesis of Remark \ref{hopf1} is introduced in order to ensure that there are \empty{sufficient} irreducible $H$--corepresentations to reconstruct $H$ from $\mathcal{T}$. This reflects the philosophy of the Tannaka--Krein duality, and the reader may consult Section 2 of reference \cite{librogroups} for a modern formulation of this result. It is worth mentioning that this is also the philosophy of this paper as well, since our aim is to reconstruct a quantum principal $H$--bundle and a quantum principal connection from $\mathcal{T}$ (together with additional information).

The proof of the following theorem can be found in Theorem 1.6.7 of reference \cite{librogroups}.

\begin{Theorem}
    \label{theoremrep}
    Under the assumption of Remark \ref{hopf}, there exists a $C^\ast$--completion $\G$ of $H$ and a extension of $\Delta$ to $$\Delta:\G\longrightarrow \G\otimes_{\mathrm{min}} \G $$ (here, $\otimes_{\mathrm{min}}$ denotes the minimal tensor product of $C^\ast$--algebras) such that $(\G,\Delta)$ is a compact quantum group.  
\end{Theorem}

It is worth mentioning that an \emph{opposite} construction is also possible  \cite{librogroups}. That is, one may start with a compact quantum group $\G$, the concept of quantum $\G$--representation and then get a dense $\ast$--Hopf algebra $H^\infty=(H,m,\mathbbm{1},\Delta,\epsilon,S,\ast)$ by taking the linear span of the matrix coefficients of all finite--dimensional quantum $\G$--representations (and hence Remark \ref{hopf} is satisfied).

Before continue, let us recall that $\Rep_H$ admits direct sums $\oplus$, tensor products $\otimes$ and complex conjugates. 

In fact, let $\delta^V$, $\delta^W$ be two finite--dimensional $H$--corepresentations, where
$$\delta^V(e^V_j)=\sum_{k} e^V_k\otimes g^V_{kj}, \qquad \delta^W(e^W_j)=\sum_{s}e^W_s\otimes g^W_{sj},$$ with $\{ e^V_j\}$, $\{ e^W_j\}$ orthonormal linear basis of $V$ and $W$, respectively. We define the direct sum of $\delta^V$ with $\delta^W$ as the (unitary) $H$--corepresentation 
\begin{equation}
    \label{ec.2.4}
    \delta^V\oplus \delta^W: V\oplus W\longrightarrow (V\oplus W)\otimes H
\end{equation}
such that  $$ (\delta^V\oplus\delta^W)(e^V_j,e^W_k)=\sum_{l} (e^V_l,0)\otimes g^V_{lj}+\sum_{s}(0,e^W_s)\otimes g^W_{sk}.$$ Furthermore, we define the direct sum of $H$--corepresentation morphisms $f_1$  $\in$ $\Mor(\delta^V,\delta^W)$, $f_2$  $\in$ $\Mor(\delta^U,\delta^Z)$ as the element of $\Mor(\delta^V\oplus \delta^U,\delta^W\oplus\delta^Z)$ given by
\begin{equation}
    \label{ec.2.5}
    \oplus(f_1,f_2):=f_1\oplus f_2: V\oplus U\longrightarrow W\oplus Z,\qquad (v,u)\longmapsto (f_1(v),f_2(u)).
\end{equation}

On the other hand, we define the tensor product of $\delta^V$ with $\delta^W$ as the (unitary) $H$--corepresentation
\begin{equation}
    \label{ec.2.6}
    \delta^V\otimes \delta^W: V\otimes W\longrightarrow (V\otimes W)\otimes H
\end{equation}
such that  $$ (\delta^V\otimes\delta^W)(e^V_j\otimes e^W_k)=\sum_{l,s} e^V_l\otimes e^W_s\otimes g^V_{kj}\,g^W_{sj}.$$
Moreover, we define the tensor product of  $H$--corepresentation morphisms $f_1$  $\in$ $\Mor(\delta^V,\delta^W)$, $f_2$  $\in$ $\Mor(\delta^U,\delta^Z)$ as the element of $\Mor(\delta^V\otimes \delta^U,\delta^W\otimes\delta^Z)$ given by
\begin{equation}
    \label{ec.2.7}
    \otimes(f_1,f_2):=f_1\otimes f_2: V\otimes U\longrightarrow W\otimes Z,\qquad v\otimes u\longmapsto f_1(v)\otimes f_2(u).
\end{equation}

Finally, given a $\C$--vector space $V$ consider its conjugate vector space $\overline{V}$. This space has the same underlying elements and additive structure as $V$, but the scalar multiplication is given by $\lambda\cdot v:=\lambda^\ast v,$ where $\lambda^\ast$ is the complex conjugate of $\lambda$ $\in$ $\C$. It is usual to denote the elements of $\overline{V}$ as $\overline{v}$. In this way,  we  define the complex conjugate $H$--corepresentation of $\delta^V$ as the (unitary) $H$--corepresentation 
\begin{equation}
    \label{ec.2.8}
    \delta^{\overline{V}}:\overline{V}\longrightarrow \overline{V}\otimes H
\end{equation}
given by $$\delta^{\overline{V}}(\overline{e}^V_j)=\sum_{k} \overline{e}^V_k\otimes g^{V\,\ast}_{kj}.$$ In addition, we define the complex conjugate of a $H$--corepresentation $f$ $\in$ $\Mor(\delta^V,\delta^W)$ as the element of  $\Mor(\delta^{\overline{V}},\delta^{\overline{W}})$ given by
$$\overline{f}:\overline{V}\longrightarrow \overline{W},\qquad \overline{v}\longrightarrow \overline{f}(\overline{v}):=\overline{f(v)}.$$ For more details, see \cite{woro1}. 

Notice that 
\begin{equation}
    \label{ec.2.10}
    \Mor(\delta^{\overline{V}},\delta^W)=\{ T:V\longrightarrow W\mid T \mbox{ is antilinear and } (T\otimes \ast)\circ \delta^V=\delta^W\circ T \}.
\end{equation}

\begin{Proposition}
    \label{prop1}
    The category $\Rep_H$ is a bar category (see Definition \ref{b.3} in Appendix B).
\end{Proposition}
\begin{proof}
    This is statement has been proven in  Section 2 of reference \cite{barcategories}, where the functor $${\bf{-}}:\Rep_H\longrightarrow \Rep_H $$ is given by $$\delta^V\longmapsto \delta^{\overline{V}},\qquad (f:V\longmapsto W)\longmapsto (\overline{f}:\overline{V}\longrightarrow \overline{W});$$ and $$\mathrm{bb}_{\delta^V}: V\longrightarrow \overline{\overline{V}}, \qquad v\longrightarrow \overline{\overline{v}},$$
    
    $$\Xi^\oplus_{\delta^V,\delta^W}:\overline{V\oplus W}\longrightarrow \overline{V}\oplus \overline{W}, \qquad \overline{(v,w)}\longmapsto (\overline{v},\overline{w}),$$
    
    $$\Xi^\otimes_{\delta^V,\delta^W}:\overline{V\otimes W}\longrightarrow \overline{W}\otimes \overline{V}, \qquad \overline{v\otimes w}\longmapsto \overline{w}\otimes \overline{v},$$ 
    
    $$\star:\C\longrightarrow \overline{C}, \qquad \lambda\longmapsto \overline{\lambda^\ast}.$$
\end{proof}

\subsection{Quantum Vector Bundles}

This subsection will be based on reference \cite{18}. 

In Differential Geometry, the Serre--Swan theorem establishes that a complex (finite--dimensional) vector bundle\footnote{Where $VM$ is the total space, $M$ is the base space which will be taken compact, and $\pi$ is the bundle projection.} $$\pi:VM\longrightarrow M$$ is equivalent to the $C^\infty_\C(M)$--bimodule of its space of global smooth sections $$\Gamma(VM),$$ which is a finitely generated and projective. Here,  $C^\infty_\C(M)$ denotes the space of all $\C$--valued smooth functions of $M$. Moreover, $\pi:VM\longrightarrow M$ is the trivial vector bundle if and only $\Gamma(VM)$ is free. This equivalence is actually categorical, and allows to define the notion of vector bundle in non--commutative geometry \cite{18}. In fact

\begin{Definition}
\label{qvb}
    Let $(B,m,\mathbbm{1},\ast)$ be a quantum space. Then, we define a  (finite--dimensional) quantum vector bundle (qvb) on $B$ as a quantum structure  $\zeta$
formally represented by a  $B$--bimodule
\begin{equation}
\label{ec.2.9}
(E,+,\cdot)
\end{equation}
which is finitely generated and projective as left and right $B$--module. The module $E$ is interpreted as the space of smooth sections of $\zeta$; so we are going to identify $\zeta$ with  $E$. Moreover, we say that a qvb over $B$ is trivial if it is free.

The space of all qvb's over $B$ will be denoted by  $$\Obj(\VB_B).$$  
\end{Definition} 

Also, we have

\begin{Definition}
    \label{qvbmorphism}
    Given two qvb $E_1$, $E_2$ over $B$, we define a qvb morphism between $E_1$ and $E_2$ as a $B$--bimodule morphism $$ A:E_1\longrightarrow E_2.$$ The set of all qvb morphisms between $E_1$ and $E_2$ will be denoted by $$\Mor(E_1,E_2).$$
\end{Definition}

\begin{Remark}
    \label{remaqvb}
    Quantum vector bundles over $B$ form a  category, denoted by $$\VB_B.$$ The notions of monomorphism, epimorphism, and isomorphism of qvb morphisms should be clear.
\end{Remark}

Before continue, let us recall that $\VB_B$  admits direct
sums $\oplus$, tensor products $\otimes$ and  conjugates.

In fact, let $E_1$, $E_2$ be two qvb's over $B$. We define the direct sum of $E_1$ with $E_2$ as the  qvb over $B$ given by  
\begin{equation}
\label{ec.2.11}
E_1\oplus E_2.
\end{equation}
We define the direct sum of qvb morphisms $A_1$ $\in$ $\Mor(E_1,E_2)$, $A_2$ $\in$ $\Mor(E_3,E_4)$ as the element of $\Mor(E_1\oplus E_3,E_2\oplus E_4)$ given by
\begin{equation}
    \label{ec.2.12}
    \oplus(A_1,A_2):=A_1\oplus A_2: E_1\oplus E_3\longrightarrow E_2\oplus E_4,\qquad (x,y)\longmapsto (A_1(x),A_2(y)).
\end{equation}

On the other hand,  we define the tensor product of $E_1$ with $E_2$ as the qvb over $B$ given by  
\begin{equation}
\label{ec.2.13}
E_1\otimes_B E_2.
\end{equation}
We define the tensor product of qvb morphisms $A_1$ $\in$ $\Mor(E_1,E_2)$, $A_2$ $\in$ $\Mor(E_3,E_4)$ as the element of $\Mor(E_1\otimes_B E_3,E_2\otimes_B E_4)$ given by 
\begin{equation}
    \label{ec.2.12.1}
    \otimes(A_1,A_2):=A_1\otimes_B A_2: E_1\otimes_B E_3\longrightarrow E_2\otimes_B E_4,\qquad x\otimes_B y\longmapsto A_1(x)\otimes_B A_2(y).
\end{equation}

Finally, for a qvb $E$, the products $$\overline{\cdot}: B\times E_1,\qquad (b,x)\longmapsto x\,b^\ast \qquad \mbox{and}\qquad  \overline{\cdot}:  E_1\times B,\qquad   (x,b)\longmapsto b^\ast x $$  equip $(E,+)$ with another $B$--bimodule structure which will be denoted by  $\overline{E}$  and it turns out to be a finitely generated projective left and right $B$--module as well. We are going to use the notation $\overline{x}$ for elements of $\overline{E}$. In this way,  we define {\it the conjugate qvb} of $E$ as the qvb over $B$ given by  
\begin{equation}
\label{ec.2.13.1}
\overline{E}.
\end{equation}
In addition, we define the conjugate of a qvb morphism $A$ $\in$ $\Mor(E_1,E_2)$ as the element of $\Mor(\overline{E}_1,\overline{E}_2)$ given by 
\begin{equation}
    \label{ec.2.14}
    \overline{A}:\overline{E}_1\longrightarrow \overline{E}_2,\qquad \overline{A}(\overline{x})=\overline{A(x)}.
\end{equation}

Notice that
\begin{equation}
    \label{ec.2.10.1}
    \Mor(E_1,\overline{E_2})=\{ A:E_1\longrightarrow E_2\mid A \mbox{ is additive and } A(b \,x)=A(x)\,b^\ast,\;  A(x\,b)=b^\ast\,A(x) \}
\end{equation}

\begin{Proposition}
    \label{prop2}
    The category $\VB_B$ is a bar category (see Definition \ref{b.3} in Appendix B).
\end{Proposition}
\begin{proof}
    This is statement has been proven in  Section 2 of reference \cite{barcategories}, where the functor $${\bf{-}}:\VB_B\longrightarrow \VB_B $$ is given by $$E\longmapsto \overline{E},\qquad (A:E_1\longmapsto E_2)\longmapsto (\overline{A}:\overline{E}_1\longrightarrow \overline{E}_2);$$ and $$\mathrm{bb}_{E}: E\longrightarrow \overline{\overline{E}}, \qquad x\longrightarrow \overline{\overline{x}},$$ 
    
    $$\Xi^\oplus_{E_1,E_2}:\overline{E_1\oplus E_2}\longrightarrow \overline{E_1}\oplus \overline{E_2}, \qquad \overline{(x,y)}\longmapsto (\overline{x},\overline{y}),$$
    
    $$\Xi^\otimes_{E_1,E_2}:\overline{E_1\otimes_B E_2}\longrightarrow \overline{E_2}\otimes_B \overline{E_1}, \qquad \overline{x\otimes_B y}\longmapsto \overline{y}\otimes_B \overline{x},$$
    
    $$\star:B\longrightarrow \overline{B}, \qquad b\longmapsto \overline{b^\ast}.$$
\end{proof}

\subsection{Quantum Principal Bundles}

As mentioned in the first section, this paper is developed within M. Durdevich's framework of quantum principal bundles. The reader may consult the reference \cite{stheve}, written by S.~Sontz, for further details. In addition, the original works \cite{micho1,micho2,micho3} provide the foundational presentation of this theory.

\begin{Definition}
\label{qpbdef}
Let $(B,m,\mathbbm{1},\ast)$ be a quantum space and let $H^\infty=(H,m,\mathbbm{1},\Delta,\epsilon,S,\ast)$ be a $\ast$--Hopf algebra. A quantum principal $H$--bundle (qpb) on $B$ is a quantum structure formally represented by the triplet
\begin{equation}
\label{f29}
\zeta=(P,B,\Delta_P),
\end{equation}
where $(P,m,\mathbbm{1},\ast)$ is a quantum space called the {\it quantum total space} with $(B,m,\mathbbm{1},\ast)$ as subspace, which receives the name of {\it quantum base space},  and $$\Delta:P \longrightarrow P\otimes H$$ is a $\ast$--algebra morphism that satisfies
\begin{enumerate}
\item $\Delta_P$ is a $H$--corepresentation.
\item $\Delta(x)=x\otimes \mathbbm{1}$ if and only if $x$ $\in$ $B$.
\item The linear map $\beta:P\otimes P\longrightarrow P\otimes H$ given by $$\beta(x\otimes y):=x\cdot \Delta_P(y):=(x\otimes \mathbbm{1})\cdot \Delta_P(y) $$ is surjective. 
\end{enumerate}
\end{Definition}

On the other hand, we have
\begin{Definition}
\label{def3}
Let $\zeta_i=(P_i,B,\Delta_{P_i})$ be a quantum principal $H_i$--bundle over $B$, for $i=1,2$. We define a qpb morphism between $\zeta_1$ and $\zeta_2$ as a pair $$(h,F),$$ where $$h:H_1\longrightarrow H_2$$ is a $\ast$--Hopf algebra morphism and  $$F:P_1\longrightarrow P_2$$ is a left $B$--module morphism such that $$(F\otimes h)\circ \Delta_{P_1}=\Delta_{P_2}\circ F.$$ 
\end{Definition}

\begin{Definition}
    \label{qpbcat}
        We define $\PB_B$ as the category whose objects are pairs  $(H^\infty,\zeta)$, where $$H^\infty=(H,m,\mathbbm{1},\Delta,\epsilon,S,\ast)$$ is $\ast$--Hopf algebra, and $$\zeta=(P,B,\Delta_P)$$ is a quantum principal $H$--bundle over $B$; and whose morphisms are qpb morphisms. 
\end{Definition}
The notions of monomorphism, epimorphism, and isomorphism of qpb morphisms should be clear.

Let $(H^\infty,\zeta)$ $\in$ $\Obj(\PB_B)$. It is worth emphasizing that, by Theorem \ref{theoremrep}, we can consider the compact quantum group $(\G,\Delta)$ associated with $H$, and we can consider $$\int_H:H\longrightarrow \C,$$ the Haar measure of $\G$ restricted to $H$ \cite{librogroups}. It is well--known that (\cite{librogroups}) $$(\id_H \otimes\int_H) \Delta=\mathbbm{1}\,\int_H\qquad \mbox{ and }\qquad \int_H\mathbbm{1}=1. $$ Then, by Theorem 5.9 of reference \cite{libro} it follows that the quotient map $$\widetilde{\beta}:P\otimes_B P\longrightarrow P\otimes H$$ induced by $\beta$ is always bijective. In this situation, one says that $P$ \emph{is a  Hopf--Galois extension of} $B$. In other words, in this paper, qpb's are always equivalent to Hopf--Galois extensions. Furthermore, Theorem 5.9 of \cite{libro} also guarantees the existence of a universal strong connection on $\zeta$. At this point in our discussion, we are not interested in introducing the notion of differential calculus on $\zeta$ (neither in Durdevich's framework nor in other approaches such as \cite{libro}). This will be addressed in the next section. In this way, in this paper, we are going to use the following definition of a universal strong connection \cite{strongconnections}.
\begin{Definition}
    \label{tonteriainecesaria}
    A universal strong connection is a unital left $B$--module morphism $$s:P\longrightarrow B\otimes P $$ such that $$ s\;\in\;\Mor(\Delta_P,\Delta_{B\otimes P}) \quad \mbox{ and }\quad m' \circ s=\id_P, $$ where $$\Delta_{B\otimes P}:=\id_B\otimes \Delta_P: B\otimes P\longrightarrow B\otimes P\otimes H $$ is a $H$--coaction, and $$m':B\otimes P\longrightarrow P\quad \mbox{ is given by }\quad m'(b\otimes p)=b\,p.$$ It is common to say that a universal strong connection is a left $B$--linear $H$--colinear splitting of $m'$.
\end{Definition}

\begin{Remark}
\label{inecesario}
In summary, in this paper, qpb's are always Hopf--Galois extensions with at least, one universal strong connection.
\end{Remark}

\subsection{Associated Quantum Vector Bundles}

This subsection will be based on references \cite{micho2,michokrein}

In Differential Geometry, given a principal $G$--bundle\footnote{Where $GM$ is the total space (a compact manifold), $M$ is the base space (a compact manifold),  $\pi$ is the bundle projection and $\cdot:GM\times G\longrightarrow GM$ is the smooth action of the (compact) Lie group $G$ on $GM$.} $$\pi:GM\longrightarrow M$$ and $$\alpha: G\longrightarrow \mathrm{GL}(V)$$ a finite--dimensional linear $G$--representation on $V$, the associated vector bundle of $\zeta$ with respect to $\alpha$ is the vector bundle $$\pi_V: V^\alpha M\longrightarrow M, \qquad [x,v]\longmapsto \pi(x),$$ where $$V^\alpha M:=GM\times_G V:=(GM\times V)/G,$$ where the action of $G$ on $GM\times V$ is given by $$(x,v,A)\longmapsto (x\cdot g,\alpha(A^{-1})v) $$ for all $x$ $\in$ $GM$, $v$ $\in$ $V$, $A$ $\in$ $G$ \cite{diff,saldgreg}. By the Serre--Swan theorem, $V^\alpha M$ is equivalent to the $C^\infty_\C(M)$--bimodule of its space of global smooth sections $$\Gamma(V^\alpha M).$$ In addition, it is well--known that as $C^\infty_\C(M)$--bimodules, $\Gamma(V^\alpha M)$ is isomorphic to the $C^\infty_\C(B)$--bimodule of $G$--equivariant maps 
\begin{equation}
    \label{ec.2.20}
    \begin{aligned}
        C^\infty_{\C}(P,V)^G= \{f:P\longrightarrow V \mid f & \mbox{ is smooth and } \\
       & f(xA)=\alpha(A^{-1})f(x) \mbox{ for all } x \,\in \, P,\, A \,\in\, G \}.
    \end{aligned}
\end{equation}

By dualizing via the pull--back the space of  $G$--equivariant maps, it follows that the \emph{non--commutative geometrical counterpart} of  $C^\infty_{\C}(P,V)^G$ is the space 
$$\Mor(\delta^V,\Delta_P),$$ where $\delta^V$ is the pull--back of $\alpha$ (identifying the dual space $V^\#$ of $V$ with $V$), and $$\Delta_P:C^\infty_\C(GM)\longrightarrow C^\infty_\C(GM\times G)\supset C^\infty_\C(GM)\otimes C^\infty_\C(G)$$ is the pull--back of the action $\cdot$ of $G$ on $GM$. In this way, in Durdevich's formulation of qpb's, we have the following definition \cite{micho2,michokrein}.

\begin{Definition}
\label{associatedqvb}
 Let $\zeta=(P,B,\Delta_P)$ be a qpb and $\delta^V$ $\in$ $\Obj(\Rep_H)$. Consider the $\C$--vector space $$E^V:=\Mor(\delta^V,\Delta_P)$$ equipped with the following $B$--bimodule structure: $$(b,T)\longmapsto b\,T,\quad \mbox{ where }\quad b\,T: V\longrightarrow P $$ is given by $(b\,T)(v)=b\,T(v)$ for all $v$ $\in$ $V$,  and 
$$(T,b)\longmapsto T\,b,\quad \mbox{ where }\quad T\,b: V\longrightarrow P $$ is given by $(T\,b)(v)=T(v)\,b$ for all $v$ $\in$ $V$. We define the associated qvb of $\zeta$ with respect to $\delta^V$ as the $B$--bimodule $E^V$.
\end{Definition}

To be precise, at this point we cannot yet consider $E^V$ as a qvb in the sense of Definition \ref{qvb}. Considering Remark \ref{inecesario} and in light of Corollary 2.6 of reference \cite{strongconnections}, we have 

\begin{Proposition}
    \label{prop2.4.1}
    Let $\delta^V$ $\in$ $\Obj(\Rep_H)$. Then $\Mor(\delta^V,\Delta_P)$ is finitely generated and projective, as left $B$--module.
\end{Proposition}

\noindent It is worth mentioning that, in the proof of the previous proposition, the hypothesis of the existence of a universal strong connection is required. 

In addition, we have

\begin{Proposition}
    \label{prop2.4.2}
    Let $\delta^V$ $\in$ $\Obj(\Rep_H)$. Then $E^V$ $\in$ $\Obj(\VB_B)$.
\end{Proposition}
\begin{proof}
    Let $\delta^V$ $\in$ $\Obj(\Rep_H)$. Then $\delta^{\overline{V}}$ $\in$ $\Obj(\Rep_{\G})$ and hence $\Mor(\delta^{\overline{V}},\Delta_P)$ is a finitely generated projective left $B$--module. Endowing $\Mor(\delta^{\overline{V}},\Delta_P)$  with the right 
$B$--module structure given by $$T\cdot b:=b^\ast\, T\qquad \mbox{ where } \qquad b^\ast\, T:\overline{V}\longrightarrow P $$ is given by $(b^\ast\, T)(\overline{v})=b^\ast\, T(\overline{v})$ for all $\overline{v}$ $\in$ $\overline{V}$, it becomes a finitely generated projective right $B$--module. Furthermore, the map 
\begin{equation*}
\ast:\Mor(\delta^{V},\Delta_P) \longrightarrow \Mor(\delta^{\overline{V}},\Delta_P),\qquad
T  \longmapsto T^\ast
\end{equation*}
is a right $B$--module isomorphism, where $T^\ast$ is defined as $(T^\ast)(\overline{v}):=T(v)^\ast$ for all $\overline{v}$ $\in$ $\overline{V}$. In this way, 
$\Mor(\delta^{V},\Delta_P)$ is finitely generated and projective, as right $B$--module. 
\end{proof}

In accordance with Proposition 6.1 of reference \cite{br2}, $$P\,\square_{H}\, V^{\#} \cong \Mor(\delta^V,\Delta_P)$$ (for the natural left $H$--coaction on $V^{\#}$, the dual space of $V$). Here, $\square_{H}$ denotes the cotensor product of $P$ and $V^{\#}$. This construction is the common one for associated qvb's in the formulation of qpb's presented in, for example \cite{libro}. Nevertheless, we have chosen to work with $E^V:=\Mor(\delta^V,\Delta_P)$ because
\begin{enumerate}
    \item This is the traditional definition of associated qvb in Durdevich's formulation (\cite{micho2}), which is the formulation we are following.
    \item The definition of the induced quantum linear connection  becomes completely analogous to their \emph{classical}  counterpart. This will be addressed in the next section (see Remark \ref{conection}).
    \item In Section 3 of reference \cite{sald2} there is a formulation of Hermitian structures on $E^V$ for which induced quantum linear connections are Hermitian. This is the no--commutative geometrical counterpart of an important result in Differential Geometry \cite{diff}. 
\end{enumerate}

\subsection{The functor $\A_{\zeta}$ (without Quantum Connections)}

This subsection will be based on \cite{michokrein,libro}.

A direct calculation proves the following proposition
\begin{Proposition}
\label{prop2.5.1}
Let $\delta^{V}$, $\delta^{W}$ $\in$ $\Obj(\Rep_H)$ and consider the associated qvb's $E^V$, $E^W$ of $\delta^{V}$, $\delta^{W}$, respectively, for a given qpb $\zeta$. Let $f$ $\in$ $\Mor(\delta^{V},\delta^{W})$. Then, the map
\begin{equation*}
A_f:E^{W} \longrightarrow E^{V},\qquad T \longmapsto T \circ f
\end{equation*}
is an element of $\Mor(E^{W},E^{V})$.
\end{Proposition}

In this way, in accordance with \cite{michokrein,libro}, we have

\begin{Definition}
\label{qfunctorq}
    Let $(H^\infty,\zeta)$ $\in$ $\Obj(\PB_B)$. We define the quantum association functor $$ \A_{\zeta}:\Rep_H\longrightarrow \VB_B$$ as the contravariant functor such that on objects is given by $$ \A_{\zeta}(\delta^V):=E^V$$ and for a morphism $f$ $\in$ $\Mor(\delta^V,\delta^W)$, we  define 
    $$\A_\zeta(f):=A_f.$$
\end{Definition}

Notice that for an element $$f\,\in\, \Mor(\delta^{\overline{V}},\delta^W)=\{ T:V\longrightarrow W\mid T \mbox{ is antilinear and } (T\otimes \ast)\circ \delta^V=\delta^W\circ T \},$$ the map $\A_\zeta(f):=A_f$ $\in$ $\Mor(E^W,E^{\overline{V}})$ can be viewed as
\begin{equation}
    \label{ec.2.21.1}
    A^\ast_f:E^W\longrightarrow E^V,\qquad T\longrightarrow T^\ast\circ f.
\end{equation}

The next proposition follows from Theorem 2.3 of reference \cite{11} and the fact that $P\,\square_{H}\, V^{\#}$ is isomorphic to $E^V$ (\cite{br2}).
\begin{Proposition}
    \label{product1}
    For every qpb $\zeta$, the quantum association functor $$\A_\zeta: \Rep_H\longrightarrow \VB_B$$ is a strict monoidal contravariant functor (see Definition \ref{b.5} in Appendix B).
\end{Proposition}

It is worth mentioning that, since (See Definition \ref{b.1} in Appendix B and reference \cite{barcategories})
$$\mathbbm{1}_{\Rep_H}=\delta^\C_\triv: \C\longrightarrow \C\otimes H,\qquad \lambda\longmapsto \lambda\otimes \mathbbm{1},$$ the  two natural isomorphisms associated with $\A_\zeta$ as a strict monoidal contravariant functor are given by  

\begin{equation}
    \label{ec.2.21}
    \phi_1: B\longrightarrow E^\C_\triv=:\Mor(\delta^\C_\triv,\Delta_P),\qquad b\longrightarrow T_b, 
\end{equation}
\begin{equation}
    \label{ec.2.22}
    \phi_2(\delta^V,\delta^W): E^V\otimes_B E^W\longrightarrow E^{V\otimes W}\quad \mbox{ is such that }\quad T^V\otimes_B T^W\longrightarrow T^V\,T^W,
\end{equation}
where $$T_b:\C\longrightarrow P,\qquad \lambda\longmapsto \lambda\,b $$ and $$T^V\,T^W:V\otimes W\longrightarrow P\quad \mbox{ is given by }\quad T^V\,T^W(v\otimes w) =T^V(v) \, T^W(w).$$

\begin{Theorem}
\label{productbar1}
    For every qpb $\zeta$, the quantum association functor $\A_\zeta$ is a contravariant bar functor (see Definition \ref{b.6} in Appendix B)
\end{Theorem}
\begin{proof}
    Since 
    \begin{equation}
        \Mor(\delta^V\oplus \delta^W,\Delta_P)\cong \Mor(\delta^V ,\Delta_P)\oplus  \Mor(\delta^W ,\Delta_P)
    \end{equation}
    is straightforward to check that $\A_\zeta$ is an additive functor.  

    First, we denote the bar category structure of $\VB_B$ (see Proposition \ref{prop2}) by an apostrophe in order to distinguish it from that of $\Rep_H$ (see Proposition \ref{prop1}).

    Let $\delta^V$ $\in$ $\Obj(\Rep_H)$. By equation (\ref{ec.2.10}) for $\delta^W=\Delta_P$ we have
    $$E^{\overline{V}}=\{ T:V\longrightarrow P\mid T \mbox{ is antilinear and } (T\otimes \ast)\circ \delta^V=\Delta_P\circ T \}.$$ Thus, consider the $B$--bimodule isomorphism 
    \begin{equation}
        \label{ec.2.26.1}
        \mathrm{bf}_{\delta^V}:\overline{E^V}\longrightarrow E^{\overline{V}},\qquad \overline{T}\longmapsto T^\ast,
    \end{equation}
     where $$T^\ast: \overline{V}\longrightarrow P, \qquad \overline{v}\longrightarrow T^\ast(\overline{v}):=T(v)^\ast. $$
    \begin{enumerate}
        \item We have $$\A_\zeta(\star^{-1})\circ \phi_1: B\longrightarrow E^{\overline{\C}}_\triv=:\Mor(\delta^{\overline{\C}},\Delta_P),\qquad b\longmapsto T_b\circ \star^{-1},$$ with $$T_b\circ \star^{-1}: \overline{\C}\longrightarrow P$$ the linear map (linear with respect to $\overline{\C}$, i.e., antilinear with respect to $\C$) given by
        $(T_b\circ \star^{-1})(\overline{1})=b.$ In addition,
         $$ \mathrm{bf}_{\delta^\C_\triv} \circ \overline{\phi_1}\circ \star':B\longrightarrow E^{\overline{\C}}_\triv,\qquad b\longrightarrow T^\ast_{{b^\ast}},$$ where $T^\ast_{{b^\ast}}$ is the linear map  (linear with respect to $\overline{\C}$, i.e., antilinear with respect to $\C$) given by $T^\ast_{{b^\ast}}(\overline{1})=(b^\ast)^\ast=b.$ Hence  $$\A_\zeta(\star^{-1})\circ \phi_1 = \mathrm{bf}_{\delta^{\mathbf{C}}_\triv}\circ \overline{\phi_1}\circ \star'.$$
         \item Let $\delta^V$ $\in$ $\Obj(\Rep_H)$. Then $$\A_\zeta(\mathrm{bb}^{-1}_{\delta^V}):=A_{\mathrm{bb}^{-1}_{\delta^V}}:E^V\longrightarrow E^{\overline{\overline{V}}},\qquad T\longmapsto A_{\mathrm{bb}^{-1}_{\delta^V}}(T):=T\circ \mathrm{bb}^{-1}_{\delta^V}$$ with $$A_{\mathrm{bb}^{-1}_{\delta^V}}(T):\overline{\overline{V}}\longmapsto P,\qquad \overline{\overline{v}}\longmapsto T(\mathrm{bb}^{-1}_{\delta^V}(\overline{\overline{v}}))=T(v);$$ and for $\A_\zeta(\delta^V)=E^V$, we have $$E^V
\xrightarrow{\mathrm{bb}'_{E^V}}
\overline{\overline{E^V}}
\xrightarrow{\overline{\mathrm{fb}_{\delta^V}}}
\overline{E^{\overline{V}}} \xrightarrow{\mathrm{fb}_{\delta^{\overline{V}}}}E^{\overline{\overline{V}}},\qquad T\longmapsto \overline{\overline{T}}\longmapsto \overline{T^\ast}\longmapsto T^{\ast\ast}$$
with $$T^{\ast\ast}:\overline{\overline{V}}\longmapsto P,\qquad \overline{\overline{v}}\longmapsto T(\overline{\overline{v}})^{\ast\ast}= T(v).$$  Therefore $$ \A_\zeta(\mathrm{bb}^{-1}_{\delta^V})= \mathrm{fb}_{\delta^{\overline{V}}}\circ\overline{\mathrm{fb}_{\delta^V}}\circ \mathrm{bb}'_{E^V}.$$
\item Let $\delta^V$, $\delta^W$ $\in$ $\Obj(\Rep_H)$. Then $$\overline{\phi^{-1}_2(\delta^V,\delta^W)}:E^{\overline{V\otimes W}}\longrightarrow \overline{E^V\otimes_B E^W}, \qquad \overline{T}^\otimes\longmapsto \sum_i \overline{T^V_i\otimes_B T^W_i}$$ such that $$\sum_i T^V_i\,T^W_i=T^\otimes.$$ Thus  $$(\Xi^{\otimes'}_{E^V,E^W}\circ\overline{\phi^{-1}_2(\delta^V,\delta^W)})(T^\otimes)= \sum_i \Xi^{\otimes'}_{E^V,E^W}(\overline{T^V_i\otimes_B T^W_i})=\sum_i\overline{T^W_i}\otimes_B \overline{T^V_i} $$ and $$(\mathrm{bf}_{\delta^W}\otimes \mathrm{fb}_{\delta^V})(\sum_i\overline{T^W_i}\otimes_B \overline{T^V_i})=\sum_i T^{W\ast}_i\otimes_B T^{V\ast}_i.$$ In this way, we get $$(\phi_2(\delta^{\overline{W}},\delta^{\overline{V}})\circ (\mathrm{bf}_{\delta^W}\otimes \mathrm{fb}_{\delta^V})\circ\Xi^{\otimes'}_{E^V,E^W}\circ\overline{\phi^{-1}_2(\delta^V,\delta^W)})(T^\otimes)=\sum_i T^{W\ast}_i\,T^{V\ast}_i,$$ where $$\sum_i T^{W\ast}_i\,T^{V\ast}_i: \overline{W}\otimes \overline{V}\longrightarrow P$$ is such that  $$\sum_i T^{W\ast}_i\,T^{V\ast}_i(\overline{w}\otimes \overline{v})=\sum_i T^{W}_i(w)^\ast\,T^{V}_i(v)^\ast=\sum_i T^\otimes(v\otimes w)^\ast.$$ On the other hand, $$\overline{E^{V\otimes W}}
\xrightarrow{\mathrm{bf}_{\delta^V\otimes \delta^W}}
E^{\overline{V\otimes W}}
\xrightarrow{\A(\Xi^{\otimes -1}_{\delta^V,\delta^W})}
E^{\overline{W}\otimes \overline{V}},\qquad \overline{T^\otimes}\longmapsto T^{\otimes \ast}\longmapsto T^{\otimes \ast}\circ \Xi^{\otimes -1}_{\delta^V,\delta^W},$$ where $$T^{\otimes \ast}\circ \Xi^{\otimes -1}_{\delta^V,\delta^W}: \overline{W}\otimes \overline{V}\longmapsto P$$ is such that $$(T^{\otimes \ast}\circ \Xi^{\otimes -1}_{\delta^V,\delta^W})(\overline{w}\otimes \overline{v})=T^{\otimes\ast}(\overline{v\otimes w})=T^\otimes(v\otimes w)^\ast.$$ Hence $$\phi_2(\delta^{\overline{W}},\delta^{\overline{V}})\circ (\mathrm{bf}_{\delta^W}\otimes \mathrm{fb}_{\delta^V})\circ\Xi^{\otimes'}_{E^V,E^W}\circ\overline{\phi^{-1}_2(\delta^V,\delta^W)}= \A(\Xi^{\otimes -1}_{\delta^V,\delta^W})\circ  \mathrm{bf}_{\delta^V\otimes \delta^W}$$
    \end{enumerate}
and the theorem follows.
\end{proof}

Let $(H^\infty,\zeta)$ $\in$ $\Obj(\PB_B)$ and $\delta^V$ $\in$ $\T$ (see equation (\ref{ec.irre})). By consider the matrix coefficients $\{g^V_{ij}\}^{n_V}_{i,j}$ of $\delta^V$ (see equation (\ref{ec.2.3})), we define the multiple irreducible subspace
\begin{equation}
    \label{ec.2.25}
    P^{V}:=\{x\in P\mid \Delta_P(x)\in P\otimes \mathrm{span}_\C\{\{g^V_{ij}\}^{n_V}_{i,j} \} \} \subseteq P
\end{equation}
associated with $\delta^V$. Each $P^{V}$ is a $B$--bimodule and the map
\begin{equation}
    \label{ec.2.26}
    E^V\otimes V\longrightarrow P^V,\qquad T\otimes v\longmapsto T(v) 
\end{equation}
is a $B$--bimodule isomorphism \cite{michokrein}. Furthermore, the following relation holds
\begin{equation}
\label{f44}
P\;\cong \bigoplus_{\delta^{V}\,\in\, \T}P^{V}\cong \bigoplus_{\delta^{V}\,\in\, \T} E^V\otimes V
\end{equation}
as $B$--bimodules. Equation (\ref{ec.2.26}) is a  $H$--corepresentation isomorphism between (\cite{micho2}) $$\Delta_P|_{P^{V}}\qquad \mbox{ and } \qquad \id_{E^{V}}\otimes \delta^{V}.$$  

 According to equation (\ref{ec.2.21}), there is a \emph{canonical} inclusion of $B$ on the right--hand side of equation (\ref{f44}) since $\delta^\C_\triv$ $\in$ $\T$. Furthermore, since every $\delta^V$ $\in$ $\Obj(\Rep_H)$ is the finite direct sum of elements of $\T$, every finite--dimensional $H$--corepresentation appears in the right--hand side of equation (\ref{f44}). In particular, this happens for $\delta^V\otimes \delta^W$ and $\delta^{\overline{V}}$ with $\delta^V$, $\delta^W$ $\in$ $\T$. Thus,  we can get a algebra structure on the right--hand side of equation (\ref{f44}) by means of 
\begin{equation}
    \label{ec.2.29}
     (T^V\otimes v)\cdot (T^W\otimes w):=\phi_2 (\delta^V,\delta^W)(T^V\otimes_B T^W)\otimes (v\otimes w);
\end{equation}
and a $\ast$ operation by means of 
\begin{equation}
    \label{ec.2.30}
    (T^V\otimes v)^\ast:=\mathrm{bf}_{\delta^V}(\overline{T^V})\otimes \overline{v}.
\end{equation}
These operations equip the right--hand side of equation (\ref{f44}) with the structure of a $\ast$--algebra. The consistency of these definitions follows directly from the functoriality of our constructions.

On the other hand, $H$ acts on the right--hand side of equation (\ref{f44}) with $$\bigoplus_{\delta^{V}\,\in\, \T}(\id_{E^{V}}\otimes \delta^{V})$$ and in accordance with Proposition 3.2 of reference \cite{michokrein},  this $H$--corepresentation is a $\ast$--algebra morphism. Hence, equation (\ref{f44}) holds as qpb's. In other words, equation (\ref{f44}) induces an isomorphism in $\PB_B$ between 
\begin{equation}
    \label{ec.2.28}
    (H^\infty,\zeta) \quad\mbox{ and }\quad (H^\infty,(\bigoplus_{\delta^{V}\,\in\, \T}P^{V},B,\bigoplus_{\delta^{V}\,\in\, \T}(\id_{E^{V}}\otimes \delta^{V}))).
\end{equation}
This shows that given $\A_{\zeta}$, we can recreate $\zeta$. The reader is encouraged to consult Section 3 of reference \cite{michokrein} for further details. 

It is worth mentioning that, in order to reconstruct $\zeta$ from $\A_{\zeta}$, it is only necessary that equation (\ref{f44}), Proposition \ref{product1} (which defines the product as in equation~(\ref{ec.2.29})), and Theorem \ref{productbar1} (which defines the $\ast$ operation as in equation (\ref{ec.2.30})) hold. For this purpose, any universal strong connection suffices. However, in~\cite{michokrein} the author uses a particular choice of universal strong connection, which we will discuss in Proposition~\ref{strongrema} of the following section.

\section{Differential Structures}

The next step is to introduce differential structures on all the spaces mentioned in the previous section. In this section, we address this subject.

\subsection{On $\ast$--Hopf Algebras}
This subsection will be based on references \cite{micho1,micho2,stheve,tomas,appendix}. Let 
\begin{equation}
    \label{ec.3.1}
    (\Gamma,d)
\end{equation}
be a bicovariant $\ast$--FODC on $H$ (see Section 6 of reference \cite{stheve}). Then, it is well--known that
\begin{equation}
    \label{ec.3.2}
    (\Gamma,d)\cong (H\otimes {\Ker(\epsilon)\over \mathcal{R}},d )
\end{equation}
for some right $H$--ideal $\mathcal{R}$ $\subseteq$ $\Ker(\epsilon)$ that satisfies $$\Ad(\mathcal{R})\subseteq \mathcal{R}\otimes H \qquad \mbox{ and }\qquad S(\mathcal{R})^\ast\subseteq \mathcal{R},$$ where $$\Ad:H\longrightarrow H\otimes H, \qquad g\longmapsto g^{(2)}\otimes S(g^{(1)})g^{(3)}$$ is the (right) adjoint coaction of $H$. Let us define  
\begin{equation}
    \label{ec.3.3}
    \mathfrak{qg}^\#:= {\Ker(\epsilon)\over \mathcal{R}}
\end{equation}
and consider the quantum germs map (see Section 6.4 of reference \cite{stheve})
\begin{equation}
\label{ec.3.4}
\begin{aligned}
\pi:H &\longrightarrow \mathfrak{qg}^\#\\
g &\longmapsto S(g^{1})dg^{2}.
\end{aligned}
\end{equation} 
Furthermore, the bicovariance of $(\Gamma,d)$ implies the existence of $\ast$--preserving linear maps
\begin{equation}
    \label{ec.3.5}
    {_\Gamma}\Phi:\Gamma \longrightarrow \Gamma\otimes H, \qquad \Phi_\Gamma:\Gamma\longrightarrow H\otimes \Gamma  
\end{equation}
such that 
\begin{equation*}
    \begin{aligned}
        {_\Gamma}\Phi(\vartheta\,g)={_\Gamma}\Phi(\vartheta)\Delta(g),&\qquad \Phi_\Gamma(g\, \vartheta)=\Delta(g)\Phi_\Gamma(\vartheta),\\
        ({_\Gamma}\Phi\otimes \id_H)\circ {_\Gamma}\Phi=(\id_H\otimes \Delta)\circ {_\Gamma}\Phi,&\qquad (\id_H\otimes \Phi_\Gamma)\circ \Phi_\Gamma=(\Delta\otimes \id_H)\circ \Phi_\Gamma,\\ (\id_H\otimes \epsilon)\circ  {_\Gamma}\Phi=\id_{H}
        ,&\qquad (\epsilon\otimes \id_H)\circ  \Phi_\Gamma=\id_{H},\\{_\Gamma}\Phi\circ d=(d\otimes \id_H)\circ \Delta,&\qquad \Phi_\Gamma\circ d=(\id_H\otimes d)\circ \Delta,
    \end{aligned}
\end{equation*}
and it is worth mentioning that ${_\Gamma}\Phi(\mathfrak{qg}^\#)\subseteq \mathfrak{qg}^\#\otimes H$; so 
\begin{equation}
\label{ec.3.6}
\ad:={_\Gamma}\Phi|_{\mathfrak{qg}^\#}: \mathfrak{qg}^\#\longrightarrow \mathfrak{qg}^\#\otimes H
\end{equation}
is a $H$--corepresentation and it fulfills (\cite{stheve})
\begin{equation}
\label{ec.3.7}
\ad \circ \pi= (\pi\otimes \id_H) \circ \Ad.
\end{equation}

On the other hand, there is a right $H$--module structure on $\mathfrak{qg}^\#$ given by 
\begin{equation}
\label{ec.3.8}
\theta \diamondsuit g:=\pi(hg-\epsilon(h)g)
\end{equation}
for every $\theta=\pi(h)$ $\in$ $\mathfrak{qg}^\#$. In particular, we have $(\theta\diamondsuit g)^\ast=\theta^{\ast}\diamondsuit S(g)^\ast$ (\cite{stheve}).

Consider the graded vector space
$$\bigotimes^\bullet_A\Gamma:=\bigoplus_k (\otimes^k_H\Gamma)\qquad \mbox{ with } \qquad \otimes^0_H \Gamma=H,\qquad  \otimes^k_H\Gamma:=\underbrace{\Gamma\otimes_H\cdots\otimes_H \Gamma}_{k\; times}$$ ($k\in \N$) endowed with its canonical graded $\ast$--algebra structure, which is given by
$$(\vartheta_1\otimes_{H}\cdots\otimes_{H}\vartheta_k)\cdot(\vartheta'_1\otimes_{H}\cdots\otimes_{H}\vartheta'_l):=\vartheta_1\otimes_{H}\cdots\otimes_{H}\vartheta_k\otimes_{H}\vartheta'_1\otimes_{H}\cdots\otimes_{H}\vartheta'_l,$$ $$(\vartheta_1\otimes_{H}\cdots\otimes_{H}\vartheta_k)^{\ast}:=(-1)^{k(k-1)\over 2}\,\vartheta^{\ast}_k\otimes_{H}\cdots\otimes_{H}\vartheta^{\ast}_1, $$ for $\vartheta_1\otimes_{H}\cdots\otimes_{H}\vartheta_k$ $\in$ $\otimes^{k}_H\Gamma$ and $\vartheta'_1\otimes_{H}\cdots\otimes_{H}\vartheta'_l$ $\in$ $\otimes^{l}_H\Gamma$. Now, let us consider the quotient graded space
  \begin{equation}
      \label{ec.3.9}     \Gamma^\wedge:=\otimes^\bullet_H\Gamma/\mathcal{Q},
  \end{equation}
  where $\mathcal{Q}$ is the two--side ideal of $\otimes^\bullet_H\Gamma$ generated by elements
  \begin{equation}
      \label{ec.3.10}
     \sum_i dg_i\otimes_H dh_i \quad \mbox{ such that } \quad \sum_i g_i\,dh_i=0,
  \end{equation}
  for all $g_i$, $h_i$ $\in$ $H$. According to \cite{micho1,stheve}, the graded $\ast$--algebra structure of $\otimes^\bullet_H \Gamma$ endows $\Gamma^\wedge$ with structure of graded $\ast$--algebra. The product in $\Gamma^\wedge$ will be denoted simply by juxtaposition of elements. On the other hand, for a given $t=\vartheta_1\cdots \vartheta_n$ $\in$ $\Gamma^{\wedge\,n }$ with $\vartheta_1$,..., $\vartheta_n$ $\in$ $\Gamma$, the linear map
  \begin{equation}
      \label{ec.3.11}
     d:\Gamma^\wedge\longrightarrow \Gamma^\wedge 
  \end{equation}
  given by $$d(t)=d(\vartheta_1\cdots \vartheta_n)=\displaystyle \sum^n_{j=1}(-1)^{j-1}\vartheta_1\cdots \vartheta_{j-1}\cdot d\vartheta_j\cdot \vartheta_{j+1}\cdots \vartheta_n \; \in \; \Gamma^{\wedge\,n+1 },$$ where $d\vartheta_j=\displaystyle \sum_l dg_l\,dh_l$ if $\vartheta_j=\displaystyle \sum_l g_l\,(dh_l)$ is well--defined, satisfies the graded Leibniz rule, $d^2=0$ and $d(t^\ast)=(dt)^\ast$ \cite{micho1,stheve}. In this way 
  \begin{equation}
    \label{ec.3.12}
     (\Gamma^\wedge,d,\ast)
\end{equation}
  is a graded differential $\ast$--algebra generated by its degree 0 elements $\Gamma^{\wedge\,0}= H$  and in references \cite{micho1,stheve}, it is called {\it the universal differential envelope $\ast$--calculus} of $(\Gamma,d)$.  It is worth mentioning that the previous construction holds for any $\ast$--FODC. In other words, at this point, the bicovariance of $(\Gamma,d)$ is not required. Moreover, as shown in Appendix B of reference \cite{micho1}, the following proposition holds for any $\ast$--FODC.

\begin{Proposition}
\label{prop3.1}
Suppose $(\Omega^\bullet=\bigoplus_k \Omega^k,d,\ast)$ is a graded differential $\ast$--algebra and $(\Gamma,d)$ is a $\ast$--FODC over $H$.

Let $$\phi^0: \Gamma^{\wedge 0}=H\longrightarrow \Omega^0$$ be a $\ast$--algebra morphism and $$\phi^1:\Gamma^{\wedge 1}=\Gamma\longrightarrow \Omega^1$$ be a linear map such that $$\phi^1(a \,db)=\phi^0(a)\,d(\phi^0(b))$$ for all $a$, $b$ $\in$ $A$. Then, there exist unique linear maps $$\phi^k: \Gamma^{\wedge k}\longrightarrow \Omega^k$$ for all $k\geq 2$ such that $$\phi:=\bigoplus_k\phi^k: \Gamma^\wedge\longrightarrow \Omega^\bullet$$ is a graded differential $\ast$--algebra morphism.
\end{Proposition}

In light of \cite{micho2,stheve}, it can be proven that for a given $\ast$--FODC $(\Gamma, d)$ on $H$, its maximal prolongation, i.e., the biggest graded differential $\ast$–algebra generated by its degree--zero elements (elements of $H$) and whose degree--one component is $\Gamma$, is $(\Gamma^\wedge, d, \ast)$.

Let $(\Gamma,d)$ be a bicovariant $\ast$--FODC on $H$ and consider its  universal differential envelope $\ast$--calculus $(\Gamma^\wedge,d,\ast)$. In this way,  we can consider the following tensor product of graded differential $\ast$--algebras  $$(\Gamma^\wedge\otimes H,d_\otimes,\ast),\qquad (H\otimes \Gamma^\wedge,d_\otimes,\ast),$$ where the structure of graded differential $\ast$--algebra on $H$ is the trivial one, i.e., $d=0$. Define ${_{\Gamma^\wedge}}\Phi^{0}=\Phi^{0}_{\Gamma^\wedge}=\Delta$, ${_{\Gamma^\wedge}}\Phi^{1}={_\Gamma}\Phi$  and $\Phi^{1}_{\Gamma^\wedge}=\Phi_\Gamma$. Thus, by Proposition \ref{prop3.1} we obtain graded differential $\ast$--algebra morphisms 
\begin{equation}
   \label{ec.3.13}
   {_{\Gamma^\wedge}}\Phi: \Gamma^\wedge \longrightarrow \Gamma^\wedge \otimes H,\qquad \Phi{_{\Gamma^\wedge}}: \Gamma^\wedge \longrightarrow H \otimes \Gamma^\wedge. 
\end{equation}
 Similarly, consider now the tensor product of $(\Gamma^\wedge,d,\ast)$ with itself $$ (\Gamma^\wedge\otimes \Gamma^\wedge,d_\otimes,\ast)$$ and by setting $\Delta^{0}=\Delta$ and $\Delta^{1}={_\Gamma}\Phi+ \Phi_{\Gamma}$, we can use Proposition \ref{prop3.1} to extend the coproduct to a graded differential $\ast$--algebra morphism
 \begin{equation}
\label{ec.3.14}
\Delta: \Gamma^\wedge \longrightarrow \Gamma^\wedge \otimes  \Gamma^\wedge.
\end{equation}
In particular, in accordance with reference \cite{micho1,stheve}, we have 
\begin{equation}
    \label{ec.3.15}
    \Delta(\theta)=\mathbbm{1}\otimes \theta+\ad(\theta)
\end{equation}
for all $\theta$ $\in$ $\mathfrak{qg}^\#$. Note that $\Delta$ coincides with ${_{\Gamma^\wedge}}\Phi + \Phi_{{\Gamma^\wedge}}$ only in degree $1$.

According to Appendix B of reference \cite{micho1}, the counit   $\epsilon$ and the antipode $S$ can also be extended to $\Gamma^\wedge$
\begin{equation}
\label{ec.3.16}
\epsilon: \Gamma^\wedge \longrightarrow \C, \qquad S: \Gamma^\wedge \longrightarrow \Gamma^\wedge;
\end{equation}
however, this part of the theory is not necessary for the purposes of this paper. If the reader is interested in these extensions, see reference \cite{micho1}.

Equation (\ref{ec.3.2}) can be extended to every degree. In fact, 
 let $(\Gamma,d)$ be a bicovariant $\ast$--FODC. Now, let us take
\begin{equation}
    \label{ec.3.17}
    \begin{aligned}
        \mathfrak{qg}^{\#\wedge}=\otimes^\bullet \mathfrak{qg}^\#/S^\wedge, \qquad \otimes^\bullet\mathfrak{qg}^\#:=\bigoplus_k (\otimes^k\mathfrak{qg}^\#),\\
    \otimes^0\mathfrak{qg}^\#=\C\mathbbm{1},\qquad \otimes^k\mathfrak{qg}^\#:=\underbrace{\mathfrak{qg}^\#\otimes\cdots\otimes \mathfrak{qg}^\#}_{k\; times}
    \end{aligned}
\end{equation}
($k\in \N$), where $S^\wedge$ is the graded two--side ideal of $\otimes^\bullet\mathfrak{qg}^\#$ generated by elements 
\begin{equation}
    \label{ec.3.18}
    \pi(g^{(1)})\otimes \pi(g^{(2)})\qquad \mbox{ for all }\qquad g \,\in\, \mathcal{R}.
\end{equation}
Then, in light of \cite{micho1,tomas}, we have
\begin{equation}
    \label{ec.3.19}
        (\Gamma^\wedge,d,\ast)\cong (H\otimes \mathfrak{qg}^{\#\wedge},d,\ast).
\end{equation}
For degree 0, the previous isomorphism is given by $H\cong H\otimes \mathbbm{1}$ in the canonical way and for degree $1$, the previous isomorphism coincides with the morphism of equation (\ref{ec.3.2}) (\cite{micho1,tomas}). Furthermore, there is a Maurer--Cartan formula (\cite{micho1,stheve})
\begin{equation}
    \label{ec.3.20}
    d\pi(g)=-\pi(g^{(1)})\pi(g^{(2)})
\end{equation}
for all $g$ $\in$ $H$. Moreover, it is worth mentioning that $\mathfrak{qg}^{\#\wedge}$ is a graded differential $\ast$--subalgebra of $\Gamma^\wedge$ (\cite{micho2}) and  we can extend the right $H$--module structure of $\mathfrak{qg}^\#$ (see equation (\ref{ec.3.8})) to $\mathfrak{qg}^{\#\wedge}$ by means of 
\begin{equation}
    \label{ec.3.21.1}
    1\diamondsuit g=\epsilon(g),\quad (\theta_1\theta_2)\diamondsuit g=(\theta_1\diamondsuit g^{(1)})(\theta_2\diamondsuit g^{(2)}).
\end{equation}
for $\theta_1$, $\theta_2$ $\in$ $\mathfrak{qg}^\#$.

Let $G$ $\subset$ $M_n(\C)$  be a compact matrix Lie group  and $H^\infty=(H,m,\mathbbm{1},\Delta,\epsilon,S,\ast)$ its canonical associated $\ast$--Hopf algebra, i.e., $H$ is the space of polynomical functions on $G$ \cite{woro1}. If a bicovariant $\ast$--FODC of $H$ is defined by $$\mathcal{R}=\Ker^2(\epsilon)=\displaystyle\{\sum^n_{i=1} a_i\,b_i\mid a_i,\,b_i\,\in\,\Ker(\epsilon),\;  n\, \in\, \N \},$$ then, according to \cite{woro2,appendix} we have 
\begin{equation}
    \label{ec.3.21}
    \Gamma\cong H\otimes \mathfrak{g}^\#_\C,\qquad \mbox{ where }\qquad \mathfrak{qg^\#}={\Ker(\epsilon)\over \mathcal{R}}={\Ker(\epsilon)\over \Ker^2(\epsilon)}=\mathfrak{g}^\#_\C
\end{equation}
is the complexification of the dual space of the Lie algebra $\mathfrak{g}$ of $G$. Moreover, the  $H$--corepresentation $\ad$ on $\mathfrak{qg^\#}=\mathfrak{g}^\#_\C$ of equation (\ref{ec.3.6}) coincides with the complexification of the pull--back of the right adjoint action of $G$ on $\mathfrak{g}$ \cite{appendix}. In addition,  $$\mathfrak{qg}^{\#\wedge}=\bigwedge \mathfrak{g}^\#_\C$$ is the exterior algebra of $\mathfrak{g}^\#_\C$ and by equation (\ref{ec.3.19}) we obtain $\Gamma^\wedge\cong H \otimes \bigwedge \mathfrak{g}^\#_\C$ \cite{appendix}. In other words, $$(\Gamma^\wedge,d,\ast)$$ is a $\ast$--subalgebra of the algebra of $\C$--valued differential forms $$(C^\infty_\C(G)\otimes \bigwedge \mathfrak{g}^\#_\C,d ) $$ of $G$, and by considering convergent sequences, one can recover the full algebra. For further details, see reference \cite{appendix}.

Therefore, we can conclude that the universal differential envelope $\ast$--calculus is a proper generalization of the algebra of $\C$--valued differential forms of $G$ in non--commutative geometry.  In this way, for a given $\ast$--Hopf algebra $H^\infty=(H,m,\mathbbm{1},\Delta,\epsilon,S,\ast)$ and a bicovariant $\ast$--FODC $(\Gamma,d)$ on $H$, the triplet $$(\Gamma^\wedge,d,\ast)$$ will be interpreted as the $\ast$--algebra of {\it quantum differential forms} of $H$. In this sense,  the space 
$$\mathfrak{qg}^\#={\Ker(\epsilon)\over \mathcal{R}}$$ plays the role of the {\it quantum dual Lie algebra} and the $H$--corepresentation $\ad$ on $\mathfrak{qg}^\#$ of equation (\ref{ec.3.6}) plays the role of the {\it dualization} of the right adjoint action of $G$ on $\mathfrak{g}$.

\subsection{On Quantum Vector Bundles}
This subsection will be based on reference \cite{18,starconnections}.

Let $E$ be a qvb over $B$. 
We say that a graded differential $\ast$--algebra
$$(\Omega^\bullet(B),d,\ast),\;\;\;\;\; \Omega^\bullet(B):=\bigoplus_{k\geq 0}\,\Omega^{k}(B)$$ is an admissible differential $\ast$--calculus if $\Omega^{0}(B)=B$ and there exists a graded $B$--bimodule isomorphism
\begin{equation}
\label{EC.3.22}
\sigma:\Omega^\bullet(B)\otimes_{B} E\longrightarrow E\otimes_{B}\Omega^\bullet(B).
\end{equation}

Whenever we work with an admissible differential $\ast$--calculus for a qvb, we assume that the morphism $\sigma$ is fixed. The graded differential $\ast$--algebra $(\Omega^\bullet(B),d,\ast)$ will be interpreted as \emph{quantum differential forms of $B$}.

\begin{Definition}
    \label{def3.1}
    Let $E$ $\in$ $\Obj(\VB_B)$ and let $(\Omega^\bullet(B),d,\ast)$ be an admissible differential $\ast$--calculus. A quantum linear connection (qlc) on $E$ is a linear map 
    \begin{equation}
\label{ec.3.23}
\nabla: E \longrightarrow \Omega^{1}(B)\otimes_{B} E 
\end{equation}
that satisfies the left and right Leibniz rule: 
\begin{equation}
    \label{ec.2.23.1}
    \nabla(bx)=b\, \nabla(x)+db\otimes_{B}x,\quad \nabla(xb)=\nabla(x)\,b+\sigma^{-1}(x\otimes_{B}db)
\end{equation}
for all $b$ $\in$ $B$ and all $x$ $\in$ $E$. 

Fix a graded differential $\ast$--algebra $(\Omega^\bullet(B),d,\ast)$ such that $\Omega^0(B)=B$. The set of  all qvb's $E$ over $B$ with qlc's $\nabla:E\longrightarrow \Omega^1(B)\otimes_B E$ will be denoted by $$\Obj(\VB^\nabla_{\Omega^\bullet(B)}).$$ Elements of $\Obj(\VB^\nabla_{\Omega^\bullet(B)})$ will be denoted by  the pair  $(E,\nabla)$.
\end{Definition}

 It is worth mentioning that qlc's depend on the choice of the admissible differential $\ast$--calculus.

Inspired by the \emph{classical} case, one may regard $E$ and $\Omega^1(B)\otimes_B E$ as the spaces of $E$--valued $0$--forms and $E$--valued $1$--forms on $B$, respectively. Moreover, via $\sigma$, the space $E\otimes_B \Omega^1(B)$ can also be interpreted as the space of $E$--valued $1$--forms on $B$. In this way, a qlc $\nabla$ can be extended to $E$--valued differential forms on $B$
\begin{equation} 
\label{def.3.2}
d^{\,\nabla}:\Omega^\bullet(B)\otimes_B E \longrightarrow \Omega^\bullet(B)\otimes_{B} E
\end{equation}
by means of $$d^{\,\nabla}(\mu\otimes_{B}x)=d\mu\otimes_{B}x+(-1)^{k}\mu\,\nabla x,$$ if $\mu$ $\in$ $\Omega^k(B)$.  Given $(E,\nabla)$, we define the curvature of $\nabla$ as the linear map
\begin{equation}
    \label{ec.3.26}
    R^{\,\nabla}:=d^{\,\nabla}\circ \nabla.
\end{equation}

\begin{Definition}
\label{def3.2}
    Let $(E_1,\nabla_1)$, $(E_2,\nabla_2)$ $\in$ $\Obj(\VB^\nabla_{\Omega^\bullet(B)})$. We say that an element of $A$ $\in$ $\Mor(E_1,E_2)$ is parallel  if
\begin{equation}
\label{ec.3.27}
(\id_{\Omega^\bullet(B)}\otimes_{B}A)\circ \nabla_1=\nabla_2\circ A,\qquad (A\otimes_{B}\id_{\Omega^\bullet(B)})\circ \sigma_1=\sigma_2\circ (\id_{\Omega^\bullet(B)}\otimes_{B}A)).
\end{equation}
The set of all parallel qvb morphisms between $E_1$ and $E_2$ will be denoted by $$\Mor((E_1,\nabla_1),(E_2,\nabla_2)).$$
\end{Definition}

 By using Definition \ref{def3.2}, it is straightforward to get relations for the map $d^{\,\nabla}$. For example $$ (\id_{\Omega^\bullet(B)}\otimes_{B}A)\circ d^{\,\nabla}_{1}=d^{\,\nabla}_{2}\circ (\id_{\Omega^\bullet(B)}\otimes_{B}A)$$ for all $A$ $\in$ $\Mor((E_1,\nabla_1),(E_2,\nabla_2))$.

\begin{Remark}
    \label{remaqvbqlc}
    Fix a graded differential $\ast$--algebra $(\Omega^\bullet(B),d,\ast)$ such that $\Omega^0(B)=B$. Quantum vector bundles over $B$ with quantum linear connections $\nabla$ (with respect to $\Omega^1(B)$) and parallel qvb morphisms form a category, denoted by $$\VB^\nabla_{\Omega^\bullet(B)}.$$ The notions of monomorphism, epimorphism, and isomorphism of parallel qvb morphisms should be clear.
\end{Remark}

Before continue, let us recall that $\VB^\nabla_{\Omega^\bullet(B)}$ admit direct
sums $\oplus$, tensor products $\otimes$ and conjugates.

In fact, let $(E_1,\nabla_1)$, $(E_2,\nabla_2)$ $\in$ $\Obj(\VB^\nabla_{\Omega^\bullet(B)})$. Then $(\Omega^\bullet(B),d,\ast)$ is an admissible differential calculus for $E_1\oplus E_2$ by means of the map
\begin{equation}
    \label{ec.3.30}
    \sigma_1\oplus \sigma_2: \Omega^\bullet(B)\otimes_B (E_1\oplus E_2)\longrightarrow (E_1\oplus E_2)\otimes_B\Omega^\bullet(B)
\end{equation}
given by $$\sigma^\oplus(\mu\otimes_B (x_1,x_2))=(\sigma_1(x_1),\sigma_2(x_2)),$$ where we have considered that $$ \Omega^\bullet(B)\otimes_B (E_1\oplus E_2)\cong (\Omega^\bullet(B)\otimes_B E_1) \oplus (\Omega^\bullet(B)\otimes_B E_2)$$ and $$ (E_1\oplus E_2)\otimes_B \Omega^\bullet(B) \cong (E_1\otimes_B\Omega^\bullet(B))  \oplus (E_2\otimes_B\Omega^\bullet(B)).$$ Moreover,  we define the direct sum of $\nabla_1$ with $\nabla_2$ as the  qlc 
\begin{equation}
\label{ec.3.31}
\nabla_1\oplus \nabla_2:E_1\oplus E_2\longrightarrow \Omega^1(B)\otimes_B (E_1\oplus E_2)\cong (\Omega^1(B)\otimes_B E_1)\oplus (\Omega^1(B)\otimes_B E_2)
\end{equation}
given by $$(\nabla_1\oplus \nabla_2)(x_1,x_2)=(\nabla_1(x_1),\nabla_2(x_2)).$$

Furthermore, $(\Omega^\bullet(B),d,\ast)$ is also an admissible differential $\ast$--calculus for $E_1\otimes_B E_2$ by means of 
\begin{equation}
\label{ec.3.32}
\sigma_1\otimes\sigma_2:=(\id_{E_1}\otimes_{B} \sigma_2)\circ (\sigma_1\otimes_{B}\id_{E_2}): \Omega^\bullet(B)\otimes_B(E_1\otimes_B E_2)\longrightarrow (E_1\otimes_B E_2)\otimes_B \Omega^\bullet(B).
\end{equation}
Also, we define the tensor product of $\nabla_1$ with $\nabla_2$ as the  qlc 
\begin{equation}
\label{ec.3.33}
\nabla_1\otimes\nabla_2:E_1\otimes_{B}E_2\longrightarrow \Omega^{1}(B)\otimes_{B}(E_1\otimes_{B}E_2) 
\end{equation}
such that $$(\nabla_1\otimes\nabla_2)(x_1\otimes_{B}x_2)=\nabla_1(x_1)\otimes_{B}x_2+(\sigma^{-1}_1\otimes_{B}\id_{E_2})(x_1\otimes_{B}\nabla_2(x_2)).$$

Now, consider the linear map $$\star_B: \Omega^1(B)\longrightarrow \overline{\Omega^1(B)},\qquad a\,db\longmapsto \overline{db^\ast\,a^\ast}, $$ where $\overline{\Omega^1(B)}$ is the conjugate of the $B$--bimodule $\Omega^1(B)$. For a qvb $E$, the map 
\begin{equation}
\label{ec.3.34}
\overline{\sigma}^{-1}:= (\star^{-1}_B\otimes \id_{\overline{E}})\circ\Xi^\otimes \circ \overline{\sigma}\circ \Xi^{\otimes-1}\circ (\id_{\overline{E}}\otimes \star_B): \overline{E}
\otimes_B \Omega^\bullet(B)\longrightarrow \Omega^\bullet(B)\otimes_B \overline{E} 
\end{equation}
tells us that $(\Omega^\bullet(B),d,\ast)$ is an admissible differential $\ast$--calculus for $\overline{E}$ as well \cite{starconnections}. Even more, for every qlc $\nabla$ on $E$, we define the conjugate qlc of $\nabla$ as the qlc (\cite{starconnections})
\begin{equation}
\label{ec.3.35}
\overline{\nabla}:= (\star_B\otimes \id_{\overline{E}})\circ \Xi^{\otimes}\circ \overline{\sigma\circ \nabla}: \overline{E}\longrightarrow \Omega^1(B)\otimes_B \overline{E}.
\end{equation}

In Theorem 2.4.2 of reference \cite{starconnections}, the reader can find a proof of the following proposition. It is worth mentioning that this proof consists of verifying that the natural isomorphisms proposed in Proposition~\ref{prop2} for $\VB_B$ are also isomorphisms for $\VB^\nabla_{\Omega^\bullet(B)}$.
\begin{Proposition}
    \label{prop3.3}
    The category $\VB^\nabla_{\Omega^\bullet(B)}$ is a bar category (see Definition \ref{b.3} in Appendix B).
\end{Proposition}

\subsection{On Quantum Principal Bundles}

As Section 2.3, this subsection will be based on Durdevich’s formulation of quantum principal bundles. The reader is encourage to consult references \cite{micho1,micho2,micho3,stheve} for further details. In addition, in reference \cite{tomas} the authors show the relation between Durdevich's formulation and the usual Brzeziński--Majid formulation, for example, presented in \cite{libro}.

\begin{Definition}
    \label{calculusqpb}
     Given $\zeta=(P,B,\Delta_P)$ a quantum principal $H$--bundle over $B$ (see Definition \ref{qpbdef}), a {\it differential calculus} on it is:
 \begin{enumerate}
 \item A graded differential $\ast$--algebra $(\Omega^\bullet(P),d,\ast)$ generated by $\Omega^0(P)=P$ ({\it quantum differential forms of $P$}).
 \item  A bicovariant $\ast$--FODC $(\Gamma,d)$ on $H$.
 \item An extension of the $\Delta_P$ to a graded differential $\ast$--algebra morphism $$\Delta_{\Omega^\bullet(P)}:\Omega^\bullet(P)\longrightarrow \Omega^\bullet(P)\otimes \Gamma^{\wedge}.$$ Here, we have considered that $\otimes$ is the tensor product of graded differential $\ast$--algebras.
 \end{enumerate}
\end{Definition}

Notice that if $\Delta_{\Omega^\bullet(P)}$ exists, then it is unique because all our graded differential $\ast$--algebras are generated by their degree $0$ elements. Furthermore,  $\Delta_{\Omega^\bullet(P)}$ is a graded differential $\Gamma^\wedge$--corepresentation on $\Omega^\bullet(P)$ \cite{micho2}.

\begin{Definition}
\label{def.horizontal}
    The space of horizontal forms is defined as 
\begin{equation}
\label{2.f18}
\Hor^\bullet P\,:=\{\varphi \in \Omega^\bullet(P)\mid \Delta_{\Omega^\bullet(P)}(\varphi)\, \in \, \Omega^\bullet(P)\otimes H \},
\end{equation}
and it is a graded  $\ast$--subalgebra of $\Omega^\bullet(P)$ (\cite{stheve}). Furthermore, by definition $$\Delta_{\Omega^\bullet(P)}(\Hor^\bullet P)\subseteq \Hor^\bullet P\otimes H,$$ so the map 
\begin{equation}
\label{2.f19}
\Delta_\Hor:=\Delta_{\Omega^\bullet(P)}|_{\Hor^\bullet P}: \Hor^\bullet P \longrightarrow \Hor^\bullet P\otimes H
\end{equation}
is a $H$--corepresentation on $\Hor^\bullet P$.
\end{Definition}

Also we have

\begin{Definition}
    \label{def.base}
    The space of {\it base} forms ({\it quantum differential forms of $B$}) is defined as 
\begin{equation}
\label{2.f20}
\Omega^\bullet(B):=\{\mu \in \Omega^\bullet(P)\mid \Delta_{\Omega^\bullet(P)}(\mu)=\mu\otimes \mathbbm{1}\}.
\end{equation}
The space of base forms is a graded differential $\ast$--subalgebra of $(\Omega^\bullet(P),d,\ast)$.
\end{Definition}

It is worth mentioning that in general, the space of base forms is not generated by $\Omega^0(B)=B$. An explicit example of this fact can be found in reference \cite{appendix}.

\begin{Definition}
    The space of vertical forms is introduced as the graded differential $\ast$--algebra (see equation (\ref{ec.3.17})) 
\begin{equation}
    \label{2.fver1}
    \Vert^\bullet P:= P\otimes \mathfrak{qg}^{\#\wedge} 
\end{equation}
(here, $\otimes$ is the tensor product of graded vector spaces) with the operations
\begin{equation}
    \label{2.fver2}
    \begin{aligned}
        (x\otimes\theta)(y\otimes \vartheta):=xy^{(0)}\otimes (\theta\diamondsuit y^{(1)})\vartheta,\\
        (x\otimes\theta)^{\ast}:=x^{(0)\ast}\otimes(\theta^{\ast}\diamondsuit x^{(1)\ast}),\\
        d_v(x\otimes \theta)=x\otimes d\theta+x^{(0)}\otimes\pi(x^{(1)})\theta,
    \end{aligned}
\end{equation}
where $x$, $y$ $\in$ $P$, $\theta$, $\vartheta$ $\in$ $\mathfrak{qg}^{\#\wedge}$, $\Delta_P(x)=x^{(0)}\otimes x^{(1)}$, $\Delta_P(y)=y^{(0)}\otimes y^{(1)}$ and $\theta \diamondsuit g$ for $g$ $\in$ $H$ is defined in equations (\ref{ec.3.8}), (\ref{ec.3.21.1}).
\end{Definition}

According to Lemma 3.1 of reference \cite{micho2}, $\Vert^\bullet P$ is generated by its degree $0$ elements $\Vert^0 P=P$. Furthermore, in accordance with Lemma 3.2 of reference \cite{micho2}, the map 
$$\Delta_\Vert: \Vert^\bullet P: \longrightarrow \Vert^\bullet P\otimes \Gamma^\wedge  $$ defined by $\Delta_\Vert(x \otimes \theta)=x^{(0)}\otimes \theta^{(1)} \otimes x^{(0)}\theta^{(2)}$ (with $\Delta(\theta)= \theta^{(1)} \otimes \theta^{(2)})$ is the unique graded differential $\ast$--algebra that is also a $\Gamma^\wedge$--corepresentation and is $\Delta_P$ in the degree $0$ case. The reader can find proofs of the two following propositions in Proposition 3.6 and Lemma 3.7 of reference \cite{micho2}, 

\begin{Proposition}
\label{piv}
The map $$\pi_\V:\Omega^\bullet(P)\longrightarrow \Vert^\bullet P $$ given by $$\pi_\V=(\id_{P}\otimes (\epsilon \otimes \id_{\mathfrak{qg}^{\#\wedge}}))\circ (\id_{\Omega^\bullet(P)}\otimes \rho_k) \circ \Delta_{\Omega^\bullet(P)}$$  is the unique graded differential $\ast$--algebra morphism  such that  $\pi_\V=\id_{P}$ in degree $0$, and $$\Delta_\Vert\circ \pi_\V  =(\pi_\V\otimes\id_{\Gamma^\wedge})\circ \Delta_{\Omega^\bullet(P)}.$$ Moreover, $\pi_\V$ is surjective. Here, $\rho_k:\Gamma^\wedge\longrightarrow \Gamma^{\wedge\, k}$ is the canonical projection onto the degree $k$ elements and we have considered that $\Gamma^\wedge=H \otimes \mathfrak{qg}^{\#\wedge}$ (see equation (\ref{ec.3.19})).
\end{Proposition}

\begin{Proposition}
    \label{seq}
    The Atiyah sequence is exact. In other words, the following sequence of $\ast$--$P$--bimodules
\begin{equation}
\label{3.f1.4}
0\longrightarrow  \Hor^1 P \lhook\joinrel\relbar\joinrel\rightarrow  \Omega^1(P) \xlongrightarrow{\pi_\V} \Vert^1 P \longrightarrow 0
\end{equation}
is always exact.
\end{Proposition}

It is worth mentioning that in Durdevich's formulation of qpb's, the exactness of the Atiyah sequence is a result of the theory, whereas in the Brzeziński--Majid formulation of qpb's, the exactness of the Atiyah sequence is imposed as a condition of the theory (see Section 5 of reference \cite{libro}).  As mentioned at the beginning of this subsection, the interested reader may consult \cite{tomas} for a discussion of the relation between these two formulations of qpb's.

By {\it dualizing} the notion of principal connections in Differential Geometry (\cite{diff}), we have 
\begin{Definition}
    \label{qpc's}
    Let $\zeta$ be a qpb with a differential calculus. A quantum principal connection (qpc) on $\zeta$ is a linear map 
    $$\omega:\mathfrak{qg}^\#\longrightarrow \Omega^{1}(P)$$ such that
    \begin{equation}
    \label{qpc1}
    \Delta_{\Omega^\bullet(P)}(\omega(\theta))=(\omega\otimes \id_G)\ad(\theta)+\mathbbm{1}\otimes\theta,
\end{equation}
   \begin{equation}
   \label{qpc2}
   \omega(\theta^\ast)=\omega(\theta)^\ast,
\end{equation}
for all $\theta$ $\in$ $\mathfrak{qg}^\#$, where $\ad$ is the $H$--corepresentation given in equation (\ref{ec.3.6}). A qpb with a qpc will be denoted by the pair  $(\zeta,\omega)$. 
\end{Definition}

A qpc is called {\it regular} if for all $\varphi$ $\in$ $\Hor^{k}P$ and $\theta$ $\in$ $\mathfrak{qg}^\#$, we have 
\begin{equation}
\label{2.f25}
\omega(\theta)\,\varphi=(-1)^{k}\varphi^{(0)}\omega(\theta\diamondsuit\varphi^{(1)}), 
\end{equation}
where $\Delta_\Hor(\varphi)=\varphi^{(0)}\otimes\varphi^{(1)}.$ A qpc $\omega$ is called {\it multiplicative} if 
\begin{equation}
\label{2.f26}
\omega(\pi(g^{(1)}))\omega(\pi(g^{(2)}))=0
\end{equation}
for all $g$ $\in$ $\mathcal{R}$, where $\Delta(g)=g^{(1)}\otimes g^{(2)}$.

In Theorems 12.8 and 12.10 of reference \cite{stheve}, the reader can find a proof of the following statement. 

\begin{Theorem}
\label{seq1}
Let $(\zeta,\omega)$ be a qpb with a qpc. Define $$\mu_{\omega}:\Vert^1 P\longrightarrow \Omega^1(P) $$ by means of $\mu_{\omega}(x\otimes \theta)=x\,\omega(\theta).$ Then 
\begin{enumerate}
\item $\mu_{\omega}$ splits the Atiyah sequence as left $P$--modules. In particular, $\mu_{\omega}$ is injective.
\item $\Delta_{\Omega^\bullet(P)} \circ \mu_{\omega}=(\mu_{\omega}\otimes \id_{\Gamma^\wedge})\circ \Delta_\Vert.$
\item $\mu_\omega(\mathbbm{1}\otimes \theta^\ast)=(\mu_\omega(\mathbbm{1}\otimes \theta))^\ast$ for all $\theta$ $\in$ $\mathfrak{qg}^\#$.
\end{enumerate}
Reciprocally, if a left $P$--module morphism $$\mu:\Vert^1 P\longrightarrow \Omega^1(P)$$ satisfies properties $1$, $2$ and $3$, then it defines a unique  qpc $\omega$ on $\zeta$ by means of $\omega(\theta):=\mu(\mathbbm{1}\otimes \theta)$.

In addition, if $\omega$ is a regular qpc, then the map $\mu_\omega$ splits the Atiyah sequence as $\ast$--$P$--bimodules.
\end{Theorem}

By {\it dualizing} the notion of covariant derivative of a principal connection in Differential Geometry (\cite{diff}), we have 

\begin{Definition}
    \label{covder}
    For a given qpc $\omega$, we define its  covariant derivative as the first--order linear map (\cite{micho3})
\begin{equation}
\label{2.f30}
D^{\omega}: \Hor^\bullet P \longrightarrow \Hor^\bullet P
\end{equation}
such that for every $\varphi$ $\in$ $\Hor^k P$ $$
D^{\omega}(\varphi)=  d\varphi-(-1)^{k}\varphi^{(0)}\omega(\pi(\varphi^{(1)})),$$ where $\Delta_\Hor(\varphi)=\varphi^{(0)}\otimes \varphi^{(1)}$. 
\end{Definition}

Direct calculations prove that (\cite{micho2,micho3})
\begin{equation}
\label{2.f31}
D^{\omega}\,  \in\,\Mor(\Delta_\Hor,\Delta_\Hor),
\qquad
D^{\omega}|_{\Omega^\bullet(B)}=d|_{\Omega^\bullet(B)}
\end{equation}
and
\begin{equation}
\label{2.f32}
\begin{aligned}
    D^{\omega}(\varphi\psi)=D^{\omega}(\varphi)\psi+(-1)^k\varphi D^{\omega}(\psi) + (-1)^k \varphi^{(0)}\ell^{\omega}(\pi(\varphi^{(1)}),\psi)
\end{aligned}
\end{equation}
for $\varphi$ $\in$ $\Hor^k\,P$, where 
\begin{equation*}
\begin{aligned}
\ell^{\omega}:\mathfrak{qg}^\#\times \Hor^\bullet P &\longrightarrow \Hor^\bullet P \\
(\theta\;,\;\varphi)\qquad&\longmapsto \omega(\theta)\varphi-(-1)^k \varphi^{(0)}\omega(\theta\diamondsuit \varphi^{(1)}).
\end{aligned}
\end{equation*}
The map $\ell^{\omega}$ {\it measures the degree of non--regularity} of $\omega$, in the sense of $\ell^{\omega}=0$ if and only if $\omega$ is regular. In other words, the covariant derivative of a qpc $D^\omega$ satisfies the graded Leibniz rule if and only if $\omega$ is regular. This is the main reason to study regular qpc's.

Moreover, for regular qpc's, we have (see Proposition 4.6 of reference \cite{micho2})
\begin{equation}
    \label{ec.5.53}
    D^\omega\circ \ast=\ast\circ D^\omega.
\end{equation}

\begin{Definition}
    \label{curvatureqpc}
    An embedded differential is a linear map 
\begin{equation}
    \label{ec.algo0}
    \Theta:\mathfrak{qg}^\#\longrightarrow \mathfrak{qg}^\#\otimes \mathfrak{qg}^\# 
\end{equation} 
such that
\begin{enumerate}
    \item  $\Theta$ $\in$ $\Mor(\ad,\ad\otimes \ad)$, where $\ad\otimes \ad$ is the $H$--corepresentation tensor product  of $\ad$ with itself. 
    \item If $\Theta(\theta)=\displaystyle\sum^m_{i,j=1}\theta_i\otimes \theta'_j$, then $d\theta=\displaystyle\sum^m_{i,j=1}\theta_i \theta'_j$ and $\Theta(\theta^\ast)=-\displaystyle\sum^m_{i,j=1}\theta'^\ast_j \otimes \theta^\ast_i$.
\end{enumerate}
Fix an embedded differential $\Theta$ and let $\omega$ be a qpc. We define the curvature of $\omega$ as the linear map
\begin{equation}
    \label{ec.curvature}
    R^\omega:\mathfrak{qg}^\#\longrightarrow \Omega^2(P)
\end{equation}
given by
$$R^\omega(\theta)=d\omega(\theta)-m_\Omega(\omega\otimes \omega)\Theta(\theta),$$ where $m_\Omega:\Omega^\bullet(P)\otimes \Omega^\bullet(P)\longrightarrow \Omega^\bullet(P)$ is the product map.
\end{Definition}

In light of Section 12.8 of reference \cite{stheve}, we have that 
\begin{equation}
    \label{ec.curva}
    \Im(R^\omega)\subseteq \Hor^2\,P\qquad \mbox{ and }\qquad R^\omega\,\in\,\Mor(\ad,\Delta_\Hor).
\end{equation}

The reason for using an embedded differential $\Theta$ in the definition of $R^\omega$ is that, in Durdevich's formulation of qpb's, the aim is for $R^\omega$ to have domain in the quantum dual Lie algebra $\mathfrak{qg}^\#$. This reflects the \emph{classical} situation in Differential Geometry, where the curvature of a principal connection takes values in the Lie algebra $\mathfrak{g}$ of the structure group $G$ of the bundle. In concrete, for a given principal connection of a principal bundle, the pull--back of the curvature map coincides with the map of equation (\ref{ec.curvature}) for $\Theta=-\displaystyle{1\over 2}c^T$, where $c^T=(\id\otimes \pi)\circ \ad$ is the transpose commutator. In the \emph{quantum} case, the properties of $\Theta$ guarantee equation (\ref{ec.curva}). For further details, see references \cite{micho2,micho3,stheve,saldym,saldcon}.

It is worth mentioning that, according to Proposition 12.14 of reference \cite{stheve}, the curvature $R^\omega$ is independent of the choice of the embedded differential $\Theta$ if  $\omega$ is multiplicative. This is the main reason for studying multiplicative qpc's. Additionally, for multiplicative qpc's we have (\cite{micho2}) 
\begin{equation}
    \label{cuadrado}
    D^{\omega\,2}(\varphi)=-\varphi^{(0)}\,R^{\omega}(\pi(\varphi^{(1)}))
\end{equation}
for all $\varphi$ $\in$ $\Hor^\bullet\,P$ and this formula determines the curvature via
\begin{equation}
    \label{cuadrado1}
    R^\omega(\pi(g))=-\sum_{k}x_k\,D^{\omega\,2}(y_k),
\end{equation}
where $x_k, y_k \in P$ satisfy $\beta(\displaystyle \sum_k x_k\otimes y_k)=\mathbbm{1}\otimes g$ (\cite{micho2}).

\subsection{An Important Restriction}

Unfortunately, the theory is very general for a categorical point of view: we have the freedom to choose too many structures that satisfy too many conditions. Hence, it is necessary to restrict it, in some way, in order to reconstructing the results presented in Section 2.5 with quantum principal connections \cite{libro}.

Let $\zeta=(P,B,\Delta_P)$ be a quantum principal $H$--bundle over $B$ and consider differential calculus on it.

\begin{Proposition}
    \label{galoishor}
    The triple $$\zeta_\Hor=(\Hor^\bullet\,P,\Omega^\bullet(B),\Delta_\Hor)$$ is a quantum principal $H$--bundle over $\Omega^\bullet(B)$.
\end{Proposition}
\begin{proof}
    Notice that we only need to prove that the map $$\beta_\Hor:\Hor^\bullet\,P\otimes \Hor^\bullet\,P\longrightarrow \Hor^\bullet\,P\otimes H$$ given by $$\beta_\Hor(\varphi\otimes \psi):=\varphi\cdot \Delta_\Hor(\psi):=(\varphi\otimes \mathbbm{1})\cdot \Delta_\Hor(\psi) $$ is surjective. In addition, since $\beta_\Hor$ is a left $(\Hor^\bullet\,P)$--bimodule, the subjectivity of $\beta_\Hor$ is equivalent to prove that $\mathbbm{1}\otimes g$ $\in$ $\Im(\beta_\Hor)$ for every $g$ $\in$ $H$.

    Let $g$ $\in$ $H$ and consider $\beta$  the map associated with $\zeta$ in Definition \ref{qpbdef}. Then $\mathbbm{1}\otimes g$ $\in$ $\Im(\beta)$ and the proposition follows from the fact that $\beta_\Hor|_P=\beta$.
\end{proof}

As we check in Section 2.3, there exists a universal strong connection 
\begin{equation}
    \label{ec.3.55}
    s_\Hor: \Hor^\bullet\,P\longrightarrow \Omega^\bullet(B)\otimes \Hor^\bullet\,P.
\end{equation}
Here, the tensor product $\otimes$ is the algebraic one. In addition, $s_\Hor|_P$ is a universal strong connection of $\zeta$. 

Therefore, all the results presented in Section 2.5 hold for the \emph{horizontal} quantum association functor
\begin{equation}
    \label{ec.3.56}
    A_{\zeta_\Hor}:\Rep_H\longrightarrow \VB_{\Omega^\bullet(B)}
\end{equation}
such that
\begin{equation}
    \label{ec.3.57}
    A_{\zeta_\Hor}(\delta_V):=\E^V:=\Mor(\delta^V,\Delta_\Hor).
\end{equation}
for every $\delta^V$ $\in$ $\Obj(\Rep_H)$, and
\begin{equation}
    \label{ec.3.58}
    A_{\zeta_\Hor}(f):A_f: \E^W\longrightarrow \E^V,\qquad \tau\longmapsto \tau\circ f
\end{equation}
for every $f$ $\in$ $\Mor(\delta^V,\delta^W)$. Of course, the $\Omega^\bullet(B)$--bimodule structure of $\E^V$ is similar to the $B$--bimodule structure of $E^V$ (see Definition \ref{associatedqvb}).

In particular we have the following construction. By consider the matrix coefficients $\{g^V_{ij}\}^{n_V}_{i,j}$ of $\delta^V$ $\in$ $\T$ (see equation (\ref{ec.2.3})), we define the multiple irreducible subspace
\begin{equation}
    \label{ec.3.59}
    \Hor^\bullet\,P^{V}:=\{\varphi\in \Hor^\bullet\,P\mid \Delta_\Hor(x)\in \Hor^\bullet\,P\otimes \mathrm{span}_\C\{\{g^V_{ij}\}^{n_V}_{i,j} \} \} \subseteq \Hor^\bullet\,P
\end{equation}
associated with $\delta^V$. Each $\Hor^\bullet\,P^{V}$ is a $\Omega^\bullet(B)$--bimodule and the map
\begin{equation}
    \label{ec.3.60}
    \E^V\otimes V\longrightarrow \Hor^\bullet\,P^{V},\qquad \tau\otimes v\longmapsto \tau(v) 
\end{equation}
is a $\Omega^\bullet(B)$--bimodule isomorphism \cite{michokrein}. Furthermore, the following relation holds
\begin{equation}
\label{ec.3.61}
\Hor^\bullet\,P\;\cong \bigoplus_{\delta^{V}\,\in\, \T}\Hor^\bullet\,P^{V}\cong \bigoplus_{\delta^{V}\,\in\, \T}\E^V\otimes V
\end{equation}
as $\Omega^\bullet(B)$--bimodules. Equation (\ref{ec.3.60}) is a  $H$--corepresentation isomorphism between (\cite{micho2}) $$\Delta_\Hor|_{\Hor^\bullet\,P^{V}}\qquad \mbox{ and } \qquad \id_{\E^{V}}\otimes \delta^{V}.$$  There is a \emph{canonical} inclusion of $\Omega^\bullet(B)$ on the right--hand side of equation (\ref{ec.3.61}) since $\delta^\C_\triv$ $\in$ $\T$. Furthermore, by using  Proposition \ref{product1} and Theorem \ref{productbar1} for $\A_{\zeta_\Hor}$, we can get a $\ast$--algebra structure on the right--hand side of equation (\ref{ec.3.61}), where the product is given by  
\begin{equation}
    \label{ec.3.62}
    (\tau^V\otimes v)\cdot (\tau^W\otimes w):=\phi_2 (\delta^V,\delta^W)(\tau^V\otimes_B \tau^W)\otimes (v\otimes w);
\end{equation}
and the $\ast$ operation is given by
\begin{equation}
    \label{ec.3.63}
    (\tau^V\otimes v)^\ast:={\mathrm{bf}}_{\delta^V}(\overline{\tau^V})\otimes \overline{v}.
\end{equation}
In other words, equation (\ref{ec.3.61}) induces an isomorphism in $\PB_{\Omega^\bullet(B)}$ between 
\begin{equation}
    \label{ec.3.64}
    (H^\infty,\zeta_\Hor) \quad\mbox{ and }\quad (H^\infty,(\bigoplus_{\delta^{V}\,\in\, \T}\Hor^\bullet\,P^{V},\Omega^\bullet(B),\bigoplus_{\delta^{V}\,\in\, \T}(\id_{\E^{V}}\otimes \delta^{V}))).
\end{equation}

Incorporating qpc's $\omega$ into the functor $\A_\zeta$ requires extending the category $\VB_B$ to the category $\VB^\nabla_{\Omega^\bullet(B)}$. However, the canonical linear map from $E^V$ to $\Omega^1(B)\otimes_B E^V$ induced by $\omega$, in general, does not satisfy the left and right Leibniz rule (see equation (\ref{ec.2.23.1})); so it cannot be considered as a quantum linear connection in the sense of Section 3.2 \cite{libro,sald2}. 

In order to ensure this and turn the quantum association functor with $\omega$ into a bar functor, the authors of~\cite{libro} impose two conditions in Section 5.4.2:
\begin{enumerate}
\item The qpb is strong.
\item Every qpc considered induces a $\ast$--bimodule splitting of the Atiyah sequence.
\end{enumerate}

In this paper, we will also impose these two conditions, but within Durdevich's framework and from a categorical point of view. The first condition is addressed in Proposition~\ref{horizontalstrong} through Remark~\ref{rema} (or to be more specific through Theorem \ref{gen}), while the second will be addressed in Remark~\ref{calculus} of the next subsection.

\begin{Remark}
\label{rema}
From this point onward until the end of the paper, we shall restrict our attention exclusively to qpb's for which the quantum base space $(B,m,\mathbbm{1},\ast)$ is a $\ast$--algebra stable under holomorphic calculus (\cite{con}).
\end{Remark}

It is worth mentioning that the previous assumption is easily satisfied and constitutes a common framework in Non--Commutative Geometry. For example, the Yang--Mills theory formulated by A. Connes in Chapter 6 of reference \cite{con} holds only for quantum spaces that are stable under holomorphic calculus. In Connes' words, only in this situation \emph{all possible notions of positivity coincide} (\cite{con}).

The hypothesis of Remark \ref{rema} was introduced exclusively in order to provide a common framework within which the following algebraic result can be guaranteed. A proof of this result can be found in Appendix B of reference \cite{micho3}.
\begin{Theorem}
    \label{gen}
    Let $\zeta=(P,B,\Delta_P)$.  Then, for every  $\delta^V$ $\in$ $\T$ (see equation (\ref{ec.irre})) there exists a set  $$\{T^V_k \}^{d_{V}}_{k=1} \subseteq \Mor(\delta^V,\Delta_P)$$ for some $d_{V}$ $\in$ $\N$ such that
\begin{equation}
    \label{generators}
\sum^{d_{V}}_{k=1}x^{V\,\ast}_{ki}x^{V}_{kj}=\delta_{ij}\mathbbm{1},
\end{equation}
with $x^{V}_{ki}:=T^V_k(e_i)$. Here $\{e_i\}^{n_{V}}_{i=1}$ is the  orthonormal basis of $V$ of equation (\ref{ec.2.3}).
\end{Theorem}

In Proposition 2.7 of reference \cite{sald2}, we determine the explicit form of the maps $\{T^V_k\}^{d_V}_{k=1}$ in the setting of Differential Geometry for a \emph{classical} principal $G$--bundle. Furthermore, in Section 5 of \cite{sald2} we determine the explicit form of the maps $\{T^V_k\}^{d_V}_{k=1}$ for trivial quantum principal bundles and for homogeneous quantum principal bundles.

\begin{Proposition}
    \label{strongrema}
    Let $\zeta=(P,B,\Delta_P)$ be a qpb. Then, the choice of the maps $\{T^V_k\}^{d_V}_{k=1}$ for each $\delta^V$ $\in$ $\T$ induces a particular choice of a universal strong connection of $\zeta$ (see Definition \ref{tonteriainecesaria}). Moreover, if $\zeta$ is equipped with a differential calculus, this universal strong connection can be extended to a universal strong connection of $\zeta_\Hor$.
\end{Proposition}

\begin{proof}
    Assume that for every $\delta^V$ $\in$ $\T$ we have chosen the maps $\{T^V_k \}^{d_{V}}_{k=1}$. According to Proposition 3.2 of reference \cite{sald2}, 
    every element $T$ of $E^V=\Mor(\delta^V,\Delta_P)$ is of the form 
    \begin{equation}
        \label{3.ec.66}
        T=\sum^{d_V}_{k=1}b^{_T}_k\,T^V_k\qquad \mbox{ with }\qquad b^{_T}_k=\sum^{n_V}_{i=1}T(e_i)\,x^{V\ast}_{ki}\,\in\,B.
    \end{equation}
    
   Consider the left $B$--module morphism $$s_{E^V}: E^V\longrightarrow B\otimes E^V, \qquad T\longmapsto \sum_k b^{_T}_k\otimes T^V_k;$$ which induces the left $B$--module morphism $$\widetilde{s}_V: E^V\otimes V\longrightarrow B\otimes E^V\otimes V \qquad \mbox{ given by }\qquad T\otimes v\longmapsto \sum_k b^{_T}_k\otimes T^V_k\otimes v.$$
   
   We claim that $\widetilde{s}_V$ is a left $B$--linear $H$--colinear splitting of $$\widetilde{m}: B\otimes E^V\otimes V\longrightarrow E^V\otimes V \qquad \mbox{ such that }\qquad \widetilde{m}(b\otimes T\otimes v)=b\,T\otimes v,$$ with respect to the $H$--corepresentations $$\id_{E^V}\otimes \delta^V\qquad \mbox{ and }\qquad \id_B\otimes \id_{E^V}\otimes \delta^V.$$ Indeed, by equation (\ref{3.ec.66}) we get 
   \begin{eqnarray*}
       \widetilde{m}(\widetilde{s}_V (T\otimes v ))=\sum_{k}\widetilde{m}(b^{_{T}}_k\otimes T^V_k\otimes v)=\sum_{k}b^{_{T}}_k\,T^V_k\otimes v
       =
       T\otimes v.
   \end{eqnarray*}
    Furthermore
   \begin{eqnarray*}
      (\id_B\otimes \id_{E^V}\otimes \delta^V)(\widetilde{s}_V(T\otimes e_j))=\sum_k b^{_T}_k\otimes T^V_k\otimes \delta^V(e_j)&=&\sum_{k,l} b^{_T}_k\otimes T^V_k\otimes e_l\otimes g^V_{lj}
      \\
     &=& 
     \sum_l \widetilde{s}_V(T\otimes e_l)\otimes g^V_{lj}
     \\
     &=&
     (\widetilde{s}_V\otimes \id_H)(\id_{E^V}\otimes \delta^V)(T\otimes e_j)
   \end{eqnarray*}
   and by linearity we conclude that $\widetilde{s}_V$ is a left $B$--linear $H$--colinear splitting of $\widetilde{m}$.

    By equation (\ref{ec.2.26}), $\widetilde{s}_V$ induces a left $B$--linear $H$--colinear splitting $$s_V: P^V\longrightarrow B\otimes P^V$$ of  $m'|_{B\otimes P^V}.$ In this way, by equation (\ref{f44}), the map $$s_M:=\bigoplus_{\delta^V\in \T}s_V:P \longrightarrow B\otimes P $$ is universal strong connection of $\zeta$.

    On the other hand, assume that $\zeta$ is equipped with a differential calculus. According to Proposition 3.3 of reference \cite{sald2}, 
    every element $\tau$ of $\E^V=\Mor(\delta^V,\Delta_\Hor)$ is of the form 
    \begin{equation}
        \label{3.ec.67}
        \tau=\sum^{d_V}_{k=1}\mu^{\tau}_k\,T^V_k\qquad \mbox{ with }\qquad \mu^{\tau}_k=\sum^{n_V}_{i=1}\tau(e_i)\,x^{V\ast}_{ki}\,\in\,\Omega^\bullet(B)
    \end{equation}
    and we can repeat the previous construction to obtain a universal strong connection $$S_M:=\bigoplus_{\delta^V\in \T}S_V:\Hor^\bullet\,P \longrightarrow \Omega^\bullet(B)\otimes \Hor^\bullet\,P $$ of $\zeta_\Hor$ which, by construction, satisfies $S_M|_{P}=s_M$.
\end{proof}  

Proposition (\ref{strongrema}) is how the assumption of Remark \ref{rema} (or, to be more precisely, the result of Theorem \ref{gen}) is related with Brzeziński--Majid formulation presented in \cite{libro}. Moreover, we have

\begin{Proposition}
    \label{horizontalstrong}
    Let $\zeta=(P,B,\Delta_P)$ be a qpb with a differential calculus. Then $$\Hor^\bullet\,P\cong \Omega^\bullet(B)\otimes_B\,P$$ as left $\Omega^\bullet(B)$--modules.
\end{Proposition}
\begin{proof}
    By equation (\ref{ec.3.61}) we get 
    \begin{equation}
        \label{ec.3.68}
        \Hor^\bullet\,P\;\cong \bigoplus_{\delta^{V}\,\in\, \T}\Hor^\bullet\,P^{V}\cong \bigoplus_{\delta^{V}\,\in\, \T} \E^V\otimes V 
    \end{equation}
    as $\Omega^\bullet(B)$--bimodules. In accordance with Proposition 3.3 of reference \cite{sald2}, the left $\Omega^\bullet(B)$--module morphism 
    \begin{equation}
        \label{ec.3.69}
        \Upsilon_V: \E^V\longrightarrow \Omega^\bullet(B)\otimes_B E^V, \qquad \tau \longmapsto \sum^{d_V}_{k=1}\mu^{\tau}_k\otimes_B T^V_k 
    \end{equation}
   is actually an isomorphism, with inverse 
   \begin{equation}
       \label{ec.3.70}
       \Upsilon^{-1}_V:\Omega^\bullet(B)\otimes_B E^V\longrightarrow \E^V,\qquad \sum_j \mu_j\otimes_B T_j\longmapsto \sum_j\mu_j\,T_j.
   \end{equation}
  Therefore, as left $\Omega^\bullet(B)$--modules we obtain 
  $$\Hor^\bullet\,P \,\cong\, \bigoplus_{\delta^{V}\,\in\, \T} \Omega^\bullet(B)\otimes_B E^V\otimes V \,\cong\,  \Omega^\bullet(B)\otimes_B \bigoplus_{\delta^{V}\,\in\, \T}  E^V\otimes V \cong \Omega^\bullet(B)\otimes_B P.$$
\end{proof}

Proposition~\ref{horizontalstrong} shows that, under the hypothesis of Remark \ref{rema}, every qpb $\zeta$ is \emph{strong} in the sense of Section 5.4.2 of reference \cite{libro}. In other words, for every $\delta^V\in \Obj(\Rep_H)$, every differential calculus on $\zeta$, and every qpc $\omega$, we can canonically induce a linear map from $E^V$ to $\Omega^1(B)\otimes_B E^V$ satisfying the left Leibniz rule (see equation (\ref{ec.2.23.1})), as we shall see shortly. In this situation, in Brzeziński--Majid formulation of qpb's (\cite{libro}), it is also common to say that every qpc  $\omega$   is strong.

\subsubsection{The Induced Quantum Linear Connection on Associated Quantum Vector Bundles}

First of all, since $$\Mor(\delta^V\oplus \delta^V,\Delta_P)\cong \Mor(\delta^V,\Delta_P)\oplus \Mor(\delta^W,\Delta_P)$$ and  $$\Mor(\delta^V\oplus \delta^V,\Delta_\Hor)\cong \Mor(\delta^V,\Delta_\Hor)\oplus \Mor(\delta^W,\Delta_\Hor)$$ and every element of $\Obj(\Rep_H)$ is a finite direct sum of elements of $\T$ (\cite{woro1}), the isomorphism $\Upsilon_V$ of equation (\ref{ec.3.69}) can be naturally extended to a $\Omega^\bullet(B)$--module isomorphism  
\begin{equation}
    \label{ec.3.71}
    \Upsilon_V:\E^V\longrightarrow \Omega^\bullet(B)\otimes_B E^V,
\end{equation}
for every $\delta^V$ $\in$ $\Obj(\Rep_H)$. As in Section 3.2, one may regard $\Omega^\bullet(B)\otimes_B E^V$ as the space of associated quantum vector bundle--valued differential forms on $B$. In addition, one may regard $\E^V$ as the space of basic differential forms on $P$ of type $\delta^V$.

On the other hand, for any qpc $\omega$, by the first part of equation (\ref{2.f31}), the covariant derivative $D^\omega$ induces the operator
\begin{equation}
    \label{ec.3.72}
    D^\omega: \E^V\longrightarrow \E^V, \qquad \tau\longmapsto D^\omega(\tau), 
\end{equation}
where $$D^\omega(\tau): V\longrightarrow \Hor^\bullet\,P,\qquad v\longmapsto D^\omega(\tau(v)).$$ 

 Let $\zeta=(P,B,\Delta_P)$ with a differential calculus and let $\delta^V$ $\in$ $\Obj(\Rep_H)$. For each qpc $\omega$, we define the linear map
\begin{equation}
    \label{ec.3.73}
    \nabla^\omega_V:E^V\longrightarrow \Omega^1(B)\otimes E^V, \qquad T\longmapsto \Upsilon_V(D^\omega(T)).
\end{equation}
 By the second part of equation (\ref{2.f31}) and equation (\ref{2.f32}), together with the fact that $\pi(\mathbbm{1})=0$, we have $$D^\omega(b\,T)=db\,T+b\,D^\omega(T)$$ for all $b$ $\in$ $B$, $T$ $\in$ $E^V$. In addition, since $$\Upsilon^{-1}_V(db\otimes_BT)=db\,T\;\;\Longrightarrow \;\; \Upsilon_V(db\,T)=db\otimes_B T,$$ it follows that 
 \begin{equation}
     \label{ec.3.73.1}
     \nabla^\omega_V(b\,T)=db\otimes_B T+ b\,\nabla^\omega_V(T);
 \end{equation}
 so $\nabla^\omega_V$ satisfies the left Leibniz rule. Note that equations (\ref{ec.3.73}), (\ref{ec.3.73.1}) hold for all qpc's $\omega$, without imposing any additional conditions. Of course, this is because $\zeta$ is strong (under the hypothesis of Remark \ref{rema}).

 According to Proposition 3.6 of reference \cite{sald2}, the following formula holds (see equation (\ref{def.3.2}))
 \begin{equation}
     \label{ec.3.73.2}
     d^{\nabla^\omega_V}=\Upsilon_V\circ D^\omega\circ \Upsilon^{-1}_V.
 \end{equation}

\begin{Remark}
\label{conection}
In the {\it classical} case, given a principal $G$--bundle $\pi:P\longrightarrow B$ and  a linear representation $\alpha:G\longrightarrow GL(V)$, there is a canonical isomorphism $\mathrm{GP}$ between basic differential forms on $P$ of type $\alpha$ and associated vector bundle--valued differential forms on $B$ \cite{saldgreg}. Moreover, this isomorphism allows to define the induced linear connection and its exterior derivative by $\mathrm{GP}(D^\omega)$ and $\mathrm{GP}\circ D^\omega\circ \mathrm{GP}^{-1}$ respectively, where $D^\omega$ is the covariant derivative of a principal connection $\omega$ \cite{diff,saldgreg}. The fact that $\Upsilon_{V}$ is an isomorphism and equations (\ref{ec.3.73}), (\ref{ec.3.73.2})  are all {\it non--commutative geometrical} counterparts of these results in differential geometry. 
\end{Remark}

Since the map $\nabla^\omega_V$ always satisfies the left Leibniz rule, it can, in the most general way, be regarded as a quantum linear connection. However, for the purposes of this paper, $\nabla^\omega_V$ is also required to satisfy the right Leibniz rule for some map $\sigma$ (see Definition~\ref{def3.1}). As noted in Section 5.4.2 of \cite{libro}, the map $\nabla^\omega_V$ satisfies the right Leibniz rule only for $\ast$--bimodule splittings of the Atiyah sequence. In Durdevich's formulation of qpb's, this is ensured by requiring that $\omega$ be regular (see Theorem \ref{seq1}).  In this way, in this paper we will work with the following definition.

\begin{Definition}
    Let $\zeta=(P,B,\Delta_P)$ be a qpb and let $\delta^V$ $\in$ $\Obj(\Rep_H)$. For a \emph{regular} quantum principal connection $\omega$, we define the induced quantum linear connection of the associated quantum vector bundle $E^V$ as the linear map $$\nabla^\omega_V:E^V\longrightarrow \Omega^1(B)\otimes E^V, \qquad T\longmapsto \Upsilon_V(D^\omega(T)).$$ 
\end{Definition}

It is worth mentioning that the required $B$--bimodule isomorphism 
\begin{equation}
    \label{ec.3.74}
    \sigma_V:\Omega^\bullet(B)\otimes_B E^V\longrightarrow E^V\otimes_B \Omega^\bullet(B)
\end{equation}
such that $$\nabla^\omega_V(T\,b)=\nabla^\omega_V(T)\,b+\sigma^{-1}_V(T\otimes_B db)$$ is given by
\begin{equation}
    \label{ec.3.75}
    \sigma_V:=\widehat{\Upsilon}_V\circ \Upsilon^{-1}_V,
\end{equation}
where $\widehat{\Upsilon}_V$ is the right $\Omega^\bullet(B)$--module isomorphism  given by 
\begin{equation}
       \label{ec.3.76}
       \widehat{\Upsilon}_V: \E^V\longrightarrow E^V\otimes_B\Omega^\bullet(B) , \qquad \tau \longmapsto \sum T^{\overline{V}\ast}_k \otimes_B  (\mu^{\tau^\ast}_k)^\ast,
   \end{equation}
where $\{T^{\overline{V}}_k \}$ is the union of the maps obtained in Theorem \ref{gen} for $\delta^{V_i}$ $\in$ $\T$ such that $\delta^{\overline{V}}=\oplus^r_{i=1} \delta^{V_i}$; and its inverse is given by   
\begin{equation}
       \label{ec.3.77}
       \widehat{\Upsilon}^{-1}_V:E^V\otimes_B \Omega^\bullet(B)\longrightarrow \E^V,\qquad \sum_j T_j\otimes_B \mu_j\longmapsto \sum_jT_j\,\mu_j.
   \end{equation}
In particular, we have
\begin{equation}
    \label{ec.3.77.1}
    \tau=\sum T^{\overline{V}\ast}_k \, (\mu^{\tau^\ast}_k)^\ast
\end{equation}
 for every $\tau$ $\in$ $\E^V$. For more details about the map $\widehat{\Upsilon}_V$, see Section 3.1 of reference \cite{sald2}.

According to Section 4.4 of reference \cite{michokrein}, the operator $\sigma_V$ satisfies
\begin{equation}
\label{f41.1}
\sigma^{-1}_{V}\circ (\id_{E^{V}}\otimes m_{\Omega^\bullet(B)}) \circ  (\sigma_{V}\otimes\id_{\Omega^\bullet(B)} )\\=(m_{\Omega^\bullet(B)}\otimes \id_{E^{V}})\;\,\circ\;\, (\id_{\Omega^\bullet(B)}\otimes \sigma^{-1}_{V}),
\end{equation}
where  $m_{\Omega^\bullet(B)}$ is the product map on $\Omega^\bullet(B)$. Thus, one can equip $\Omega^\bullet(B)\otimes_{B} E^{V}$ with a (graded) $\Omega^\bullet(B)$--bimodule structure, where the left multiplication is just 
\begin{equation}
    \label{f41.2}
    m_{\Omega^\bullet(B)}\otimes_B\id_{E^{V}}
\end{equation}
and the right multiplication is 
\begin{equation}
    \label{f41.3}
    \sigma^{-1}_{V}\circ (\id_{E^{V}}\otimes_{B}m_{\Omega^\bullet(B)})\circ (\sigma_{V}\otimes \id_{\Omega^\bullet(B)}).
\end{equation}
With this new structure, $\sigma_{V}$ becomes a (graded) $\Omega^\bullet(B)$--bimodule isomorphism \cite{michokrein}.

\subsection{Another Important Restriction}

Since, from a categorical point of view, we require regular quantum principal connections, we must impose conditions to guarantee the existence of such connections. This constitutes the second restriction that we impose in this paper. For this purpose, we will follow the theory presented in Section 6.5 of reference \cite{micho2} and reference \cite{michodif}.

Let $\zeta=(P,B,\Delta_P)$ be a qpb. Assume the existence of a graded $\ast$--algebra 
\begin{equation}
\label{ec.3.78}
    (\Omega^\bullet_\Hor,m,\mathbbm{1},\ast)\quad \mbox{ such that }\quad \Omega^{0}_\Hor=P
\end{equation}
with a graded $\ast$--subalgebra 
\begin{equation}
    \label{ec.3.79}
    (\Omega^\bullet_B,m,\mathbbm{1},\ast)
\end{equation}
with structure of graded differential $\ast$--algebra 
\begin{equation}
    \label{ec.3.80}
    d:\Omega^\bullet_B\longrightarrow \Omega^\bullet_B
\end{equation}
such that $\Omega^{0}_B=B$. Note that there is no need to assume that $(\Omega^\bullet_B,d,\ast)$ is generated by its degree~$0$ elements. Moreover, assume that the map $\Delta_P$ can be extended to a graded $\ast$--algebra morphism 
\begin{equation}
    \label{ec.3.81}
    \Delta_{\mathrm{H}}:\Omega^\bullet_\Hor\longrightarrow \Omega^\bullet_\Hor\otimes H
\end{equation}
such that $\Delta_{\mathrm{H}}$ is a $H$--corepresentation. In addition, we will assume that  $\Omega^\bullet_B$ is exactly the space of all $\Delta_{\mathrm{H}}$--invariant elements. In addition, assume that 
\begin{equation}
    \label{ec.3.82}
    \mathfrak{Der}\not=\emptyset,
\end{equation}
where 
\begin{equation}
    \label{ec.3.82.1}
    \mathfrak{Der}
\end{equation}
is the set of all first--order linear maps $$D:\Omega^\bullet_\Hor \longrightarrow \Omega^\bullet_\Hor$$ such that
\begin{enumerate}
    \item $\ast\circ D=D\circ \ast$.
    \item $D(\varphi\,\psi)=D(\varphi)\,\psi+(-1)^k\,\varphi\,D(\psi)$ for all $\varphi$ $\in$ $\Omega^k_\Hor$.
    \item $D|_{\Omega^\bullet(B)}=d$.
    \item $D\in \Mor(\Delta_\H,\Delta_H)$.
\end{enumerate}

In accordance with \cite{michodif}, there exists a bicovariant $\ast$--FODC on $H$
\begin{equation}
    \label{bico}
    (\Lambda,d)
\end{equation}
and a differential calculus on $\zeta$ for this $\ast$--FODC.

However, the construction of $(\Lambda,d)$ is somewhat restrictive (\cite{michodif}); let us \emph{weaken} it. Fix an element $D$ $\in$ $\mathfrak{Der}$. According to Lemma 2.2 of reference \cite{michodif}, there exists a linear map $$r:H\longrightarrow  \Omega^2_\Hor\quad \mbox{ such that }\quad D^2(\varphi)=-\varphi^{(0)}\,r(\varphi^{(1)}).$$  Thus, choose a bicovariant $\ast$--FODC on $H$
\begin{equation}
    \label{bico1}
    (\Gamma,d) \qquad \mbox{ with }\qquad \mathfrak{qg}^\#={\Ker(\epsilon)\over \mathcal{R}}
\end{equation}
 such that $\mathcal{R}$ is a subspace of the space of all elements annihilated by $r$. In this condition, the quotient map of $r$
\begin{equation}
    \label{curv1}
    R: \mathfrak{qg}^\#\longrightarrow \Omega^2_\Hor 
\end{equation}
is well--defined and it satisfies $$D^2(\varphi)=-\varphi^{(0)}\,R(\pi(\varphi^{(1)})).$$ Such a bicovariant $\ast$--FODC on $H$ always exists, since, for instance, $(\Lambda,d)$ provides a suitable example (see~\cite{michodif}). In accordance with Section 6.5 of reference \cite{micho2}, there is a graded differential $\ast$--algebra structure on (here, $\otimes$ is the tensor product of graded vector spaces)
\begin{equation}
    \label{ec.3.83}
    \widetilde{\Omega}^\bullet(P)=\Omega^\bullet_\H\otimes \mathfrak{qg}^{\#\wedge}
\end{equation}
given by
$$(\psi\otimes \theta)\cdot (\varphi\otimes \vartheta):=(-1)^{\partial\varphi\,\partial\theta} \psi\varphi^{(0)}\otimes (\theta \diamondsuit \varphi^{(1)})\vartheta,$$ $$(\psi\otimes \theta)^\ast=\psi^{(0)\ast}\otimes (\theta^\ast\diamondsuit \psi^{(1)\ast}),$$ $$d(\psi\otimes \theta)=D(\psi)\otimes \theta+(-1)^{\partial\psi}\psi\,(R(\theta)+d\theta)+(-1)^{\partial\psi}\psi^{(0)}\otimes \pi(\psi^{(1)})\theta,$$
where $\psi$, $\varphi$ $\in$ $\Omega^\bullet_\H$, $\theta$, $\vartheta$ $\in$ $\mathfrak{qg}^{\#\wedge}$, $\Delta_\H(\psi)=\psi^{(0)}\otimes \psi^{(1)}$, $\Delta_\H(\varphi)=\varphi^{(0)}\otimes \varphi^{(1)}$, $\theta \diamondsuit g$ for $g$ $\in$ $H$ is defined in equations (\ref{ec.3.8}), (\ref{ec.3.21.1}) and $\partial x$ denotes the grade of the element $x$. In addition, with the previous structure,  $\widetilde{\Omega}^\bullet(P)$ is generated by its degree--$0$ elements $\widetilde{\Omega}^0(P)=P$ and there exists an extension of $\Delta_P$ to a graded differential $\ast$--algebra morphism
\begin{equation}
    \label{ec.3.84}
\Delta_{\widetilde{\Omega}^\bullet(P)}:\widetilde{\Omega}^\bullet(P)\longrightarrow \widetilde{\Omega}^\bullet(P)\otimes \Gamma^{\wedge}.
\end{equation}
In other words, we get a differential calculus on $\zeta$. Furthermore, the space of horizontal forms $\Hor^\bullet\,P$ of $\widetilde{\Omega}^\bullet(P)$ is $$\Omega^\bullet_\Hor\otimes \mathbbm{1}\cong \Omega^\bullet_\Hor,$$ the space of base forms $\Omega^\bullet(B)$ of $\widetilde{\Omega}^\bullet(P)$ is $$\Omega^\bullet_B\otimes \mathbbm{1}\cong \Omega^\bullet_B$$ and under the previous canonical isomorphisms, we get $$\Delta_\Hor:=\Delta_{\widetilde{\Omega}^\bullet(P)}|_{\Hor^\bullet\,P}=\Delta_\H.$$ Of course, if the reader is interested in the details of the previous construction, we recommend consulting \cite{micho2}. For example, in Proposition 6.21 of \cite{micho2} we find that
\begin{Proposition}
    \label{regcon}
    The linear map $$\omega^c: \mathfrak{qg}^\#\longrightarrow \Omega^1(P),\qquad \theta\longmapsto \mathbbm{1}\otimes \theta $$ is a regular and multiplicative qpc such that $D^{\omega^c}=D$ and $R^{\omega^c}=R$.
\end{Proposition}

In summary, from the data
\begin{equation}
    \label{ec.3.85}
    \{\Delta_\H,\,\Omega^\bullet_\H,\,\{(\Omega^\bullet_B,d,\ast) \}, D,\,(\Gamma,d)\}
\end{equation}
we can induce a differential calculus on $\zeta$. At first sight, it may seem rather convoluted to determine the data in equation \eqref{ec.3.85}, instead of directly specifying a graded differential $\ast$--algebra generated by $P$ together with an extension $\Delta_{\Omega^\bullet(P)}$ of $\Delta_P$. However, as we will see in the following section, the data in equation \eqref{ec.3.85} can be obtained in a \emph{functorial manner}. In addition, Proposition \ref{regcon} guarantees the existence of regular qpc's, as we require.

\begin{Remark}
    \label{calculus}
   From now on, we assume that every differential calculus on a qpb is determined by some spaces $(\Gamma,d)$, $(\Omega^\bullet_\H,m,\mathbbm{1},\ast)$, $(\Omega^\bullet_B,d,\ast)$, together with a $H$--corepresentation $\Delta_\H$ and a map $D$ $\in$ $\mathfrak{Der}$ in the previous manner.
\end{Remark}

It is worth mentioning that, for a fixed quantum principal $H$--bundle $\zeta$ with differential calculus given by equations (\ref{bico1}), (\ref{ec.3.83}) and (\ref{ec.3.84}),  if $\omega^c\neq \omega'$ is another regular and multiplicative qpc, then its covariant derivative $D^{\omega'}$ is an element of $\mathfrak{Der}$ for which the map $R$ can be defined in $\mathfrak{qg}^\#$ (see Section 3.3). Hence, one can construct another differential calculus on $\zeta$ using $D^{\omega'}$ in such a way that, in this new differential calculus, $\omega'$ takes the form $$ \mathfrak{qg}^\#\longrightarrow \Omega^1(P),\qquad \theta\longmapsto \mathbbm{1}\otimes \theta$$ in the sense that the covariant derivative of this qpc is the operator $D^{\omega'}$. This may be interpreted as a kind of \emph{translation} from $\omega'$ to $\omega^c$. Accordingly, we define
\begin{equation}
\label{ec.3.86}
[(\Gamma,d),\zeta,\omega^c]
\end{equation}
to be the equivalence class of such \emph{translations}. Notice that, for every element of $[(\Gamma,d),\zeta,\omega^c]$, the space of quantum differential forms on $H$ is $\Gamma^\wedge$, the horizontal space is $\Omega^\bullet_\H$, and the space of base forms is $\Omega^\bullet_B$.

Now consider $\mathcal{O}$ the collection of all equivalence classes $[(\Gamma,d),\zeta,\omega^c]$. For each $[(\Gamma,d),\zeta,\omega^c]$ $\in$ $\mathcal{O}$, choose an element $((\Gamma,d),\zeta,\omega^c)$ $\in$ $[(\Gamma,d),\zeta,\omega^c]$ and let 
$$\Obj(\PB^{\omega^c}_{\Omega^\bullet(B)}) $$
be the collection of all such chosen $((\Gamma,d),\zeta,\omega^c)$, one for each equivalent class of $\mathcal{O}$.

\begin{Definition}
    \label{qpbcat2}
    Fix a graded differential $\ast$--algebra $(\Omega^\bullet(B)=\Omega^\bullet_B,d,\ast)$ with $\Omega^0(B)=B$. We define $\PB^{\omega^c}_{\Omega^\bullet(B)}$ as the category whose objects are $\Obj(\PB^{\omega^c}_{\Omega^\bullet(B)}).$ Furthermore, the morphisms of this category are pairs $$(h,F),$$ where $$h:H_1\oplus\Gamma_1\longrightarrow H_2\oplus\Gamma_2$$ is a graded linear map such that $h|_{H_1}:H_1\longrightarrow H_2$ is a $\ast$--Hopf algebra morphism and $$h(g_1dg_2)=h(g_1)d(h_2)$$ for every $g_1$, $g_2$ $\in$ $H_1$; and  $$F:\Hor^\bullet\,P_1\longrightarrow \Hor^\bullet\,P_2$$ is a graded $\Omega^\bullet(B)$--module morphism such that $$(\F\otimes h)\circ \Delta_{\Hor_1}=\Delta_{\Hor_2}\circ F.$$ 
\end{Definition}
The notions of monomorphism, epimorphism, and isomorphism should be clear.
 
Of course, the category $\PB^{\omega^c}_{\Omega^\bullet(B)}$ is not unique, since different choices in the selection of $\Obj(\PB^{\omega^c}_{\Omega^\bullet(B)})$ lead to different categories. However, by construction, all these categories are equivalent.

\section{The Quantum Association Functor with Quantum Connections}

With all the notions presented and discussed in the previous sections, finally we can define the desired functor $\A^{\omega^c}_{\zeta}$.

\subsection{Some Properties}

First of all, we have
\begin{Proposition}
\label{prop4.1}
Consider an element of $((\Gamma,d),\zeta,\omega^c)$ $\in$ $\Obj(\PB^{\omega^c}_{\Omega^\bullet(B)})$ and let $\delta^{V}$, $\delta^{W}$ $\in$ $\Obj(\Rep_H)$. If $f$ $\in$ $\Mor(\delta^{V},\delta^{W})$, the map
\begin{equation*}
A_f:E^{W} \longrightarrow E^{V},\qquad T \longmapsto T \circ f
\end{equation*}
is an element of $\Mor((E^{W},\nabla^{\omega^c}_W),(E^{V},\nabla^{\omega^c}_V))$.
\end{Proposition}
\begin{proof}
Since every element of $\Obj(\Rep_H)$ is the finite direct sum of elements of $\T$ (\cite{woro1}), it is sufficient to prove the proposition for $\delta^V$, $\delta^W$ $\in$ $\T$.

For every $T$ $\in$ $E^W$ (see equations (\ref{3.ec.67}), (\ref{ec.3.69})) $$\nabla^{\omega^c}_W(T)=\Upsilon_{W}(D^{\omega^c}(T))=\sum_k\mu^{D^{\omega^c}(T)}_k\otimes_{B}T^W_k;$$ so 
\begin{equation}
    \label{ec.4.1}
    (\id_{\Omega^{1}(B)}\otimes_{B}A_f)\nabla^{\omega^c}_W(T)=\sum_k\mu^{D^{\omega^c}(T)}_k\otimes_{B} (T^W_k\circ f).
\end{equation}
Moreover 
\begin{equation}
    \label{ec.4.2}
    \nabla^{\omega^c}_V(A_f(T))= \Upsilon_{V}(D^{\omega^c}(A_f(T)))= \sum_{k}\mu^{D^{\omega^c}(A_f(T))}_k\otimes_B T^V_k=\sum_{k}\mu^{D^{\omega^c}(T\circ f)}_k\otimes_B T^V_k.
\end{equation}
By applying $\Upsilon^{-1}_{V}$ in equation (\ref{ec.4.1}) one gets $$\sum_k\mu^{D^{\omega^c}(T)}_k\, (T^W_k\circ f)=\sum_k(\mu^{D^{\omega^c}(T)}_k\, T^W_k)\circ f=D^{\omega^c}(T)\circ f=D^{\omega^c}(T\circ f) $$ and by applying $\Upsilon^{-1}_{V}$ in equation (\ref{ec.4.1}) one obtains $$\sum_{k}\mu^{D^{\omega^c}(A_f(T))}_k\otimes_B T^V_k=\sum_{k}\mu^{D^{\omega^c}(T\circ f)}_k\, T^V_k=D^{\omega^c}(T\circ f).$$ Since $\Upsilon^{-1}_{V}$ is bijective, we conclude that $$(\id_{\Omega^{1}(B)}\otimes_{B}A_f)\circ \nabla^{\omega^c}_W=\nabla^{\omega^c}_V\circ A_f.$$

On the other hand (see equation (\ref{ec.3.76}))
\begin{equation*}
    (A_f\otimes_B\id_{\Omega^\bullet(B)})\sigma_W(\mu\otimes_B T)=(A_f\otimes_B\id_{\Omega^\bullet(B)})\widehat{\Upsilon}_W(\mu\, T)=\sum (T^{\overline{W}\ast}_k\circ f)\otimes_B (\mu^{(\mu\,T)^\ast}_k)^\ast
\end{equation*}
with $(T^{\overline{W}\ast}_k\circ f)(\overline{w})=(T^{\overline{W}}_k(f(w)))^\ast$ for all $\overline{w}$ $\in$ $\overline{W}$, so
\begin{equation}
 \label{ec.4.3}
  (A_f\otimes_B\id_{\Omega^\bullet(B)})\sigma_W(\mu\otimes_B T)=\sum (T^{\overline{W}}_k\circ f)^\ast\otimes_B (\mu^{(\mu\,T)^\ast}_k)^\ast;
\end{equation}
and
\begin{equation}
\label{ec.4.4}
    \sigma_V(\id_{\Omega^\bullet}\otimes_B A_f)(\mu\otimes_B T)=\sum (T^{\overline{V}\ast}_k \otimes_B (\mu^{(\mu\,T\circ f)^\ast}_k)^\ast.
\end{equation}
By applying $\widehat{\Upsilon}^{-1}_{V}$ in equation (\ref{ec.4.3}) one gets 
\begin{eqnarray*}
    \sum (T^{\overline{W}}_k\circ f)^\ast\, (\mu^{(\mu\,T)^\ast}_k)^\ast=\sum (\mu^{(\mu\,T)^\ast}_k\,T^{\overline{W}}_k\circ f)^\ast&=&  ((\mu\,T)^\ast\circ f)^\ast 
    \\
    &=& ((T^\ast \circ f)\,\mu^\ast)^\ast
    \\
    &=&
    ((T \circ f)^\ast\,\mu^\ast)^\ast
    =\mu\,(T\circ f)
\end{eqnarray*}
and by applying $\widehat{\Upsilon}^{-1}_{V}$ in equation (\ref{ec.4.4}) one obtains 
 $$\sum T^{\overline{V}\ast}_k \, (\mu^{(\mu\,T\circ f)^\ast}_k)^\ast= \mu\,(T\circ f),$$
 according to equation (\ref{ec.3.77.1}). By linearity and since $\widehat{\Upsilon}^{-1}_{V}$ is bijective, we obtain $$(A_f\otimes_B\id_{\Omega^\bullet(B)})\circ \sigma_W=\sigma_V\circ(\id_{\Omega^\bullet}\otimes_B A_f) $$
\end{proof}

In this way, we have

\begin{Definition}
\label{qfunctorq2}
    Consider an element of $((\Gamma,d),\zeta,\omega^c)$ $\in$ $\Obj(\PB^{\omega^c}_{\Omega^\bullet(B)})$. We define the quantum association functor $$ \A^{\omega^c}_{\zeta}:\Rep_H\longrightarrow \VB_{\Omega^\bullet(B)}$$ as the contravariant functor such that on objects is given by $$ \A^{\omega^c}_{\zeta}(\delta^V):=(E^V,\nabla^{\omega^c}_V)$$ and for a morphism $f$ $\in$ $\Mor(\delta^V,\delta^W)$, we  define 
    $$\A_\zeta(f):=A_f.$$
\end{Definition}

As in Section 2.5, we have

\begin{Theorem}
\label{productbar2}
    The quantum association functor $\A^{\omega^c}_\zeta$ is a contravariant bar functor.
\end{Theorem}

\begin{proof}
    The proof consists of a large but direct calculation verifying that all the properties of Definition \ref{b.6} in Appendix B  are satisfied for the natural isomorphisms (see equations (\ref{ec.2.21}), (\ref{ec.2.22}), (\ref{ec.2.26.1}))
    $$\phi_1:(B,d)\longrightarrow (E^\C_\triv, \nabla^{\omega^c}_\C), \qquad b\longmapsto T_b,$$ $$\phi_2(\delta^V,\delta^W):(E^V\otimes_B E^W, \nabla^{\omega^c}_V\otimes \nabla^{\omega^c}_W)\longrightarrow (E^{V\otimes W},\nabla^{\omega^c}_{V\otimes W}),\qquad T^V\otimes_B T^W\longrightarrow T^V\,T^W,$$ $$\mathrm{bf}_{\delta^V}:(\overline{E^V},\overline{\nabla^{\omega^c}_V})\longrightarrow (E^{\overline{V}},\nabla^{\omega^c}_{\overline{V}}),\qquad \overline{T}\longmapsto T^\ast$$ as in the proof of Theorem \ref{productbar1}, so will omit it.
\end{proof}

A similar statement can be found in Proposition 5.56 and 5.57 of reference \cite{libro}, but the domain of the association functor is the category of right $H$-crossed
modules (see Definition 2.22 of \cite{libro}).
 
Consider an element of $((\Gamma,d),\zeta,\omega^c)$ $\in$ $\Obj(\PB^{\omega^c}_{\Omega^\bullet(B)})$. We have already checked that it is possible to recreate the qpb and the horizontal space using the quantum association functor. Let  $\delta^V$ $\in$ $\T$. We define a first--order linear map 
\begin{equation}
\label{f46}
    D_V:\E^V\longrightarrow \E^V
\end{equation}
given by
\begin{equation*}
D_V:= \Upsilon^{-1}_{V}\circ d^{\nabla^{\omega}_{V}} \circ \Upsilon_{V}.
\end{equation*}
Notice that $D_V(\tau)=D^{\omega{^c}}(\tau)$. According to Section 4.4 of reference \cite{michokrein}, the map $D_V$ satisfies 
\begin{equation}
\label{f47}
D_{\C_\triv}=d|_{\Omega^\bullet(B)}, \qquad D_{V\otimes W}(\tau_1\cdot \tau_2)=D_{V}(\tau_1)\tau_2+(-1)^{k}\tau_1D_{W}(\tau_2), 
\end{equation}
\begin{equation}
    \label{f48}
    D_{\overline{V}}\circ \A_{\zeta_\Hor}(\overline{\id}^{\;-1}_{V})=\A_{\zeta_\Hor}(\overline{\id}^{\;-1}_{V})\circ D_{V} \quad\mbox{ and }\quad D_{V}\circ \A_{\zeta_\Hor}(f)=\A_{\zeta_\Hor}(f) \circ D_{W}
\end{equation}
for ever $f$ $\in$ $\Mor(\delta^{V},\delta^{W})$, where $$\overline{\id}_V:V\longrightarrow \overline{V}, \qquad v\longrightarrow \overline{v}.$$

By using these  properties one can induce a first--order linear map on the right--hand side of equation (\ref{ec.3.61}) which coincides with $D^{\omega^c}$. Now, using the method to create differential calculus on qpb presented in Section 3.5, we get a bicovariant $\ast$--FODC $(\Gamma,d)$ on $H$ and $\omega^c$ as the linear map given by $$\theta \longmapsto \mathbbm{1}\otimes \theta.$$ 

In summary, given $\A^{\omega^c}_\zeta$ (and hence, $\A_{\zeta_\Hor}$), we can recreate   $$((\Gamma,d),\zeta,\omega^c)\;\in\; \Obj(\PB^{\omega^c}_{\Omega^\bullet(B)}).$$ It is worth mentioning that this construction, up to isomorphisms, does not depend on $\T$.

To conclude this subsection, according to Proposition 3.6 of reference \cite{sald2}, we get $$d^{\nabla^{\omega^c}_{V}}= \Upsilon_{V}\circ D^{\omega^c}\circ \Upsilon^{-1}_{V}.$$ Since $\omega^c$ is multiplicative, by equation (\ref{cuadrado}) we have
$$R^{\nabla^{\omega^c}_{V}}(T)= \Upsilon_{V}\circ D^{\omega^c\,2}\circ T=-\Upsilon_{V} \circ T^{(0)}R^{\omega^c}(\pi(T^{(1)})),$$ where $\Delta(T(v))=T^{(0)}(v)\otimes T^{(1)}(v)$ for all $\delta^V$ $\in$ $\Obj(\Rep_H)$. 

\subsection{The Categorical Equivalence}

Now, we will proceed to prove the categorical equivalence. We shall begin with the following technical result.

\begin{Theorem}
\label{teo1}
Fix  a graded differential $\ast$--algebra $(\Omega^\bullet(B),d,\ast)$ with $\Omega^0(B)=B$ and consider a bicovariant $\ast$--FODC $(\Gamma,d)$ on a $\ast$--Hopf algebra $H^\infty=(H,m,\mathbbm{1},\Delta,\epsilon,S,\ast)$. Let $$\F: \Rep_{H} \longrightarrow \VB^{\nabla}_{\Omega^\bullet(B)}$$ be a contravariant bar functor.
Then, there exists a quantum principal $H$--bundle $\zeta$ over $B$ for which $\F$ is naturally isomorphic to $\A^{\omega^c}_{\zeta}$. 
\end{Theorem}
\begin{proof}
Let $\delta^V$ $\in$ $\T$ and consider $$\F(\delta^{V})=(\widehat{E}^V,\widehat{\nabla}_V).$$ Since $\F$ is a contravariant bar functor, there exists natural isomorphisms 
\begin{equation}
    \label{ec.4.7}
    \widehat{\phi}_1:(B,d)\longrightarrow (\widehat{E}^\C_\triv, \widehat{\nabla}_{\C}),
\end{equation}
\begin{equation}
\label{ec.4.8}
    \widehat{\phi}_2(\delta^V,\delta^W):(\widehat{E}^V\otimes_B \widehat{E}^W, \widehat{\nabla}_{V}\otimes \widehat{\nabla}_W)\longrightarrow (\widehat{E}^{V\otimes W},\widehat{\nabla}_{V\otimes W}),
\end{equation}
\begin{equation}
    \label{ec.4.9}
    \widehat{\mathrm{bf}}_{\delta^V}:(\overline{\widehat{E}^V},\overline{\widehat{\nabla}_V})\longrightarrow (\widehat{E}^{\overline{V}},\widehat{\nabla}_{\overline{V}})
\end{equation}
such that all the properties of a contravariant bar functor are satisfy. In addition, the composition of $\F$ with the forget functor $$\VB^\nabla_{\Omega^\bullet(B)}\longrightarrow \VB_B$$ is also a contravariant bar functor for the same natural isomorphisms $\widehat{\phi}_1,$ $\widehat{\phi}_2,$ $\widehat{\mathrm{bf}}$. In this way, we define the $B$--bimodule 
\begin{equation}
    \label{ec.4.10}
    P:=\bigoplus_{\delta^{V}\,\in\,\T}\widehat{E}^V\otimes V.  
\end{equation}
Notice that $H$ naturally coacts on $P$ via 
\begin{equation}
    \label{ec.4.11}
    \Delta_P:=\bigoplus_{\delta^{V}\,\in\,\T} \id_{\widehat{E}^V}\otimes \delta^{V}.
\end{equation}

For every $\delta^V$ $\in$ $\T$, we have $\widehat{E}^V\subseteq \Mor(\delta^V,\Delta_P)$. Indeed, for every $x$ $\in$ $E^V$, $x$ can be considered as the linear map $$x:V\longrightarrow P,\qquad v\longmapsto x\otimes v$$ and we obtain $x$ $\in$ $\Mor(\delta^V,\Delta_P)$. 

According to equation (\ref{ec.4.7}), there is a \emph{canonical} inclusion of $B$ on $P$ since $\delta^\C_\triv$ $\in$ $\T$. Furthermore, since every $\delta^V$ $\in$ $\Obj(\Rep_H)$ is the finite direct sum of elements of $\T$, every finite--dimensional $H$--corepresentation appears in $P$. In particular, this happens for $\delta^V\otimes \delta^W$ and $\delta^{\overline{V}}$ with $\delta^V$, $\delta^W$ $\in$ $\T$. Thus,  we can get a algebra structure on $P$ by means of 
\begin{equation}
    \label{ec.4.13}
    (x^V\otimes v)\cdot (x^W\otimes w):=\widehat{\phi}_2 (\delta^V,\delta^W)(x^V\otimes_B x^W)\otimes(v\otimes w)
\end{equation}
and a $\ast$ operation by means of 
\begin{equation}
    \label{ec.4.14}
    (x^V\otimes v)^\ast:=\mathrm{bf}_{\delta^V}(\overline{x^V})\otimes \overline{v}.
\end{equation}
These operations equip $P$ with the structure of a $\ast$--algebra and  in accordance with Theorem 3.4 of reference \cite{michokrein}, the triple
\begin{equation}
    \label{ec.4.15}
    \zeta=(P,B,\Delta_P)
\end{equation}
is a quantum principal $H$--bundle over $B$. Now, we can conclude that
\begin{equation}
    \label{ec.4.16}
    \widehat{E}^V=\Mor(\delta^V,\Delta_P).
\end{equation}

On the other hand, since $\F(\delta^V)=(\widehat{E}^V,\widehat{\nabla}_V)$  is object of $\VB^{\nabla}_{\Omega^\bullet(B)}$, there exists a graded $B$--bimodule isomorphism 
\begin{equation}
    \label{ec.4.17}
    \widehat{\sigma}^V:\Omega^\bullet(B)\otimes_B \widehat{E}^V\longrightarrow \widehat{E}^V\otimes_B \Omega^\bullet(B)
\end{equation}
such that $\widehat{\nabla}_V$ satisfies the left and right Leibniz rule (see Definition  \ref{def3.1}). 

Define a contravariant functor 
\begin{equation}
    \label{ec.4.18}
    \F^{\Hor}:\Rep_H\longrightarrow \VB_{\Omega^\bullet(B)}
\end{equation}
such that in objects is $$\F^{\Hor}(\delta^V)=\Omega^\bullet(B)\otimes_B\widehat{E^V}$$ with the (graded) $\Omega^\bullet(B)$--bimodule structure defined in equations (\ref{f41.2}), (\ref{f41.3}) for $\widehat{\sigma}^V$, and in morphism is $$\F^{\Hor}(f)=\widehat{A}^{\,\Hor}_f:=\id_{\Omega^\bullet(B)}\otimes_B \widehat{A}_f,$$ with $\F(f)=\widehat{A}_f$.  

 In this way, the natural isomorphisms
\begin{equation}
    \label{ec.4.19}
    \id_{\Omega^\bullet(B)}\otimes_B\widehat{\phi}_1,\qquad \id_{\Omega^\bullet(B)}\otimes_B\widehat{\phi}_2 \qquad \id_{\Omega^\bullet(B)}\otimes_B\widehat{\mathrm{bf}}_{\delta^V}
\end{equation}
turns $\F^\Hor$ into a contravariant bar functor and we can repeat the above process to get the qpb 
\begin{equation}
    \label{ec.4.20}
    \zeta_\Hor=(\Hor^\bullet\,P,\Omega^\bullet(B),\Delta_\Hor),
\end{equation}
where 
\begin{equation}
    \label{ec.4.21}
\Hor^\bullet\,P:=\bigoplus_{\delta^{V}\,\in\,\T}\Omega^\bullet(B)\otimes_B \widehat{E}^V\otimes V
\end{equation}
and
\begin{equation}
    \label{ec.4.22}
    \Delta_\Hor:=\bigoplus_{\delta^{V}\,\in\,\T} \id_{\Omega^\bullet(B)}\otimes_B \id_{\widehat{E}^V}\otimes \delta^{V}.
\end{equation}
In particular, and to be more precisely, we have
\begin{equation}
    \label{ec.4.23}
    \Omega^\bullet(B)\otimes_B E^V\cong \Mor(\delta^V,\Delta_\Hor)=:\widehat{\E}^V.
\end{equation}

The previous isomorphism agrees with the map (taking in consideration the corresponding structures) $$\mu\otimes_B T\longmapsto \mu\,T.$$ Denoting this map by $\Upsilon^{-1}_{V}$, the isomorphism between $\widehat{E^V}\otimes_{B}\Omega^\bullet(B)$ and $\widehat{\E}^V$ is given by $$\widetilde{\Upsilon}^{-1}_{V}:=\Upsilon^{-1}_{V} \circ \widehat{\sigma}^{-1}_{V},$$ so equation (\ref{ec.3.75}) is satisfied and $\widetilde{\Upsilon}^{-1}_{V}$ agrees with the map $T\otimes_B \mu\longmapsto T\mu$.

For every $\delta^{V}$ $\in$ $\T$, we define a first--order linear map 
\begin{equation}
    \label{ec.4.24}
    \widehat{D}_{V}:\widehat{\E}^V\longrightarrow \widehat{\E}^V 
\end{equation}
as the operator of equation (\ref{f46}). Since $\widehat{\E}^V=\Mor(\delta^V,\Delta_\Hor)$, we have  $$(D_V(\tau)\otimes \id_{H})\circ \delta^V=\Delta_\Hor\circ D_V(\tau) $$ for every $\tau$ $\in$ $\widehat{\E}^V$ and one can proves that  $\widehat{D}_{V}$ satisfies equations (\ref{f47}), (\ref{f48}) as well.

In this way, we define the first--order linear map
\begin{equation}
    \label{ec.4.25}
    \widehat{D}: \Hor^\bullet\,P\longrightarrow \Hor^\bullet\,P
\end{equation}
given by $$\widehat{D}(\tau^V\otimes v)=D_V(\tau^V)\otimes v$$ for every $\tau^V\otimes v$ $\in$ $\E^V\otimes V$ and is, actually, an element of $\mathfrak{Der}$ (see equation (\ref{ec.3.82.1}). 

Now, for a bicovariant $\ast$--FODC on $H$
$$(\Gamma,d)$$ such the map $R$ can be defined (see equation (\ref{curv1})), we can apply the method of Section 3.5 (the one of Section 6.5 of reference \cite{micho2}) to find a differential calculus on $\zeta$ for which the covariant derivative of the qpc $\omega^c$ is $\widehat{D}$. By construction, $\F$ and $\A^{\omega^c}_{\zeta}$ are natural isomorphic.
\end{proof}

\begin{Definition}
    \label{qgts}
    Let us denote by $\Obj(\GTS^{\nabla}_{\Omega^\bullet(B)})$ the set of all tuples $$((\Gamma,d),\F),$$ where $$\F:\Rep_H\longrightarrow \VB^\nabla_{\Omega^\bullet(B)}$$ is a contravariant bar functor and $(\Gamma,d)$ is the bicovariant $\ast$--FODC on $H$ for which $\F$ is natural isomorphic to $\A^{\omega^c}_\zeta$ by Theorem \ref{teo1}.

    On the other hand, we denote by $\Mor(((\Gamma_1,d_1),\F_1),((\Gamma_2,d_2),\F_2))$ the set of all tuples $$(h,\NT),$$ where $$h:H_1\oplus \Gamma_1\longrightarrow H_2\oplus\Gamma_2$$ is a graded linear map such that  $h|_{H_1}:H_1\longrightarrow H_2$ is a $\ast$--Hopf algebra morphism and $$h(g_1dg_2)=h(g_1)d(h_2)$$ for every $g_1$, $g_2$ $\in$ $H_1$; and 
    $$\NT: \F_1\longrightarrow \F_2\, \hat{h}$$ is a natural transformation such that $\NT$ commutes with $\bigoplus$,  $\bigotimes$ and $\bf{-}$, with $$\hat{h}:\Rep_{H_1}\longrightarrow \Rep_{H_2}$$ the covariant functor given by $$\hat{h}\, \delta^V=(\id_{V}\otimes h)\circ \delta^V$$ and $$\hat{h}f=f$$ for morphisms. 

    Finally, this define the category $$\GTS^{\nabla}_{\Omega^\bullet(B)},$$ which will be called \emph{the category of quantum gauge theory sectors over $\Omega^\bullet(B)$}.
\end{Definition}


 Now it is possible to generalize the quantum association functor  $\A^{\omega^c}_{\zeta}$ as a covariant functor  $$\A:\PB^{\omega^c}_{\Omega^\bullet(B)} \longrightarrow \GTS^{\nabla}_{\Omega^\bullet(B)}$$ such that on objects is defined by $$\A((\Gamma,d),\zeta,\omega^c)=((\Gamma,d),\A^{\omega^c}_{\zeta})$$ and on morphisms is defined by  $$\A(h,F):=(h,\hat{F})$$ where $$\hat{F}:\A^{\omega_1}_{\zeta_1}\longrightarrow \A^{\omega_2}_{\zeta_2}\, \hat{h} $$ is the natural transformation given by 
\begin{equation*}
\hat{F}\,\delta^V: \A^{\omega^c}_{\zeta_1}(\delta^V) \longrightarrow \A^{\omega^c}_{\zeta_2}(\hat{h}(\delta^V)),\qquad T  \longmapsto F\circ T.
\end{equation*}

In Theorem \ref{teo2}, we will prove that the category $\PB^{\omega^c}_{\Omega^\bullet(B)}$ is equivalent to $\GTS^{\nabla}_{\Omega^\bullet(B)}$. This is the purpose of this paper, and finally we have all the tools and all the context necessary to prove it.

\begin{Theorem}
\label{teo2}
Fix a graded differential $\ast$--algebra $(\Omega^\bullet(B)=\Omega^\bullet_B,d,\ast)$ with $\Omega^0(B)=B$. Then, the functor $\A$ provides an equivalence of categories from $\PB^{\omega^c}_{\Omega^\bullet(B)}$ to $\GTS^{\nabla}_{\Omega^\bullet(B)}$.
\end{Theorem}

\begin{proof}
By Theorem~\ref{teo1}, every  $((\Gamma,d),\F)$ $\in$ $\Obj(\GTS^{\nabla}_B)$ is isomorphic in $\GTS^{\nabla}_{\Omega^\bullet(B)}$ to $((\Gamma,d),\A^{\omega^c}_{\zeta})$, for an element $((\Gamma,d),\zeta,\omega^c)$ $\in$ $\Obj(\PB^{\omega^c}_{\Omega^\bullet(B)})$. Thus, we just need to show that $\A$ induces a bijection between $$\Mor((\Gamma_1,d_1),\zeta_1,\omega^c),((\Gamma_2,d_2),\zeta_2,\omega^c)) \quad \mbox{ and } \quad  \Mor(((\Gamma_1,d_1),\A^{\omega^c}_{\zeta_1}),((\Gamma_2,d_2),\A^{\omega^c}_{\zeta_2})).$$ According to our previous results, we know that $$P_i\cong\bigoplus_{\delta^{V_i}\,\in\,\T_i}E^{V_i}\otimes V_i,  $$ $$\Hor^\bullet\,P_i\cong\bigoplus_{\delta^{V_i}\,\in\,\T_i}\Omega^\bullet(B)\otimes_B E^{V_i}\otimes V_i,$$ and $$\Delta_{P_i}\cong \bigoplus_{\delta^{V_i}\,\in\,\T_i} \id_{E^{V_i}}\otimes \delta^{V_i}, \qquad  \Delta_{\Hor_i} \cong \bigoplus_{\delta^{V_i}\,\in\,\T_i} \id_{\Omega^\bullet(B)}\otimes_B\id_{E^{V_i}}\otimes \delta^{V_i},$$ where  $\T_i$ is a complete set of mutually non--equivalent irreducible $H_i$--corepresentations for $i=1$, $2$. 

Taking $$(h,\NT)\;\in\;\Mor(((\Gamma_1,d_1),\A^{\omega^c}_{\zeta_1}),((\Gamma_2,d_2),\A^{\omega^c}_{\zeta_2})),$$ let us define
\begin{equation*}
F:\bigoplus_{\delta^{V_1}\,\in\,\T_1}\Omega^\bullet(B)\otimes_B E^{V_1}\otimes V_1 \longrightarrow \bigoplus_{\delta^{V_2}\,\in\,\T_2}\Omega^\bullet(B)\otimes_B E^{V_2}\otimes V_2
\end{equation*}
such that $$F(\mu\otimes_B T\otimes v):=\bigoplus_{\delta^{V_2}\in \T_2}\mu\otimes_B T_{V_2}\otimes v_{V_2}$$ if $T$ $\in$ $E^{V_1}$, and $$\NT\,\delta^{V_1} (T)=\bigoplus_{\delta^{V_2}\in \T_2}T_{V_2} \;\in\;  \Mor(\hat{h}\,\delta^{V_1},\Delta_P)\subseteq \bigoplus_{\delta^{V_2}\in \T_2}E^{V_2},$$ $$v= \bigoplus_{\delta^{V_2}\in \T_2}v_{V_2}\, \in \, V_1=\bigoplus_{\delta^{V_2}\in \T_2} V_2,$$ where in the last three expressions we have used the same finite number of corepresentations $\{ \delta_2\}$ $\in$ $\T_2$.

Clearly, $F$ is a graded $\Omega^\bullet(B)$--module morphism and hence $$(h,F)\;\in\;\Mor((\Gamma_1,d_1),\zeta_1,\omega^c),((\Gamma_2,d_2),\zeta_2,\omega^c))$$ and by construction, $\hat{F}=\NT$. This implies that $\A(h,F)=(h,\NT)$. 

On the other hand, if $$\A(h_1,F_1)=(h_1,\hat{F}_1)=(h_2,\hat{F}_2)=\A(h_2,F_2),$$ we automatically get $h_1=h_2$, and since $\hat{F}_1=\hat{F}_2$ we obtain  $F_1\circ T=F_2\circ T$ for all $T$ $\in$ $E^{V_1}$ and all $\delta^{V_1}$ $\in$ $\T_1$. By considering the decomposition of $\Hor^\bullet\,P_1$ into the direct sum, it follows that $F_1=F_2$; thus $$(h_1,F_1)=(h_2,F_2).$$ 
\end{proof}

In Differential Geometry, fix a (compact) manifold $M$ and consider $\mathbf{PB}^{\omega}_M$, the category of all principal $G$--bundles $$\pi: GM\longrightarrow M$$ equipped with principal connections $\omega$. Then, there is a categorical equivalence between $\mathbf{PB}^{\omega}_M$ and the category $\mathbf{GTS}^{\nabla}_M$ of gauge theory sectors (see~\cite{saldgreg} or the brief summary given in Appendix~A). The category $\mathbf{GTS}^{\nabla}_M$ consists of all association functors of principal bundles with principal connections; this is the reason for its name.

More precisely, each association functor (or, equivalently, each principal $G$--bundle $\pi: GM \longrightarrow M$ with a principal connection $\omega$) $$\A^\omega_{GM}:\mathbf{Rep}_G\longrightarrow \mathbf{VB}^\nabla_M  $$ defines a gauge theory on $M$, where $\mathbf{Rep}_G$ is the category of finite--dimensional linear representations of $G$, and $\mathbf{VB}^\nabla_M$ is the category of vector bundles over $M$ with linear connections. In concrete, for each such association functor $\A^\omega_{GM}$, we obtain all associated vector bundles $$\pi_V:V^\alpha M\longrightarrow M$$ of $\pi: GM\longrightarrow M,$ together with the linear connections $\nabla^\omega_V$ induced by $\omega$.

For instance, consider the principal bundle 
$$\pi: GM:=\R^4\times U(1)\longrightarrow \R^4,\qquad (x,\mathrm{e}^{it})\longmapsto x,$$ equipped with a principal connection $\omega$. Then $$\A^\omega_{\R^4\times U(1)}:\mathbf{Rep}_{U(1)}\longrightarrow \mathbf{VB}^\nabla_{\R^4} $$
defines the gauge theory corresponding to the electromagnetic theory on $M$ for the principal connection $\omega$.

For example, if $\omega$ is a solution of the Yang--Mills equation with current\footnote{An equation given in $\mathfrak{u}(1)M:=\A^\omega_{\mathbb{R}^4\times U(1)}(\ad)$, where $\ad$ denotes the adjoint action of $U(1)$ on its Lie algebra $\mathfrak{u}(1)$; see reference \cite{diff}.}, then any global smooth section of the associated bundle $\pi_V:V^\alpha M\longrightarrow M$ is interpreted as a matter field interacting with the electromagnetic field produced by $\omega$ via $\nabla^\omega_V$ (see \cite{diff}).

As in the \emph{classical} case, in Non--Commutative Geometry, by Theorem~\ref{teo2} we know that $\PB^{\omega}_{\Omega^\bullet(B)}$ is equivalent to $\GTS^{\nabla}_{\Omega^\bullet(B)}$. Hence, every object of $\GTS^{\nabla}_{\Omega^\bullet(B)}$ arises as the quantum association functor of some qpb $\zeta=(P,B,\Delta_P)$, equipped with a bicovariant $\ast$--FODC $(\Gamma,d)$ on $H$ and the corresponding qpc $\omega^c$. This is the reason for the name of\footnote{Note that, in Differential Geometry, $\Omega^\bullet(B)$ and $(\Gamma,d)$ are canonically fixed by considering differential forms on the corresponding spaces. Furthermore, in Differential Geometry there is no need to restrict to a particular class of principal connections. In contrast, the generality of Non--Commutative Geometry forces us to choose $\Omega^\bullet(B)$, $(\Gamma,d)$, and $\omega^c$ in order to reconstruct the categorical equivalence.} $\GTS^{\nabla}_{\Omega^\bullet(B)}$.

More precisely, each quantum association functor (or, equivalently, each quantum principal $H$--bundle $\zeta=(P,B,\Delta_P)$ equipped with a bicovariant $\ast$--FODC $(\Gamma,d)$ on $H$ and the qpc $\omega^c$) $$\A^{\omega^c}_{\zeta}:\Rep_H\longrightarrow \VB^{\nabla}_{\Omega^\bullet(B)} $$ defines a gauge theory on $B$ equipped with the graded differential $\ast$--algebra $\Omega^\bullet(B)$.

In particular, for each such quantum association functor $\A^{\omega^c}_{\zeta}$, we obtain all associated quantum vector bundles $E^V$ of $\zeta$, together with the quantum linear connections $\nabla^{\omega^c}_V$ induced by $\omega^c$.

To conclude this section, it is worth emphasizing that we have imposed several restrictions in order to derive the categorical equivalence of Theorem \ref{teo2}. However, in general, these restrictions are not necessary for the definition of the quantum association functor.

For example, for every qpc $\omega$ (under the assumption of Remark \ref{rema}, or, if one does not wish to use Remark \ref{rema}, for every strong qpc; see \cite{libro}), the quantum association functor $$\A^{\omega}_{\zeta}:\Rep_H\longrightarrow \VB^{\nabla}_{\Omega^\bullet(B)} $$ can be defined, provided that the category $\VB^{\nabla}_{\Omega^\bullet(B)}$ is suitably modified so that the qlc’s $\nabla$ satisfy only the left Leibniz rule (that is, ignoring the existence of the map $\sigma$ and the right Leibniz rule). Although, in this situation, $\A^{\omega}_{\zeta}$ is no longer a contravariant bar functor. 

Even in the situations described above, we can still regard $\A^{\omega}_{\zeta}$ as defining a gauge theory on $B$ equipped with the graded differential $\ast$--algebra $\Omega^\bullet(B)$. We encourage the interested reader to consult~\cite{saldym,saldcon} for further details. In these references, we develop the theory of Yang--Mills fields coupled to scalar matter fields and to fermionic matter fields using $\A^{\omega}_{\zeta}$ (for arbitrary qpc’s, not only regular ones, since these paper are developed under the assumption of Remark~\ref{rema}). However, in those works, the categorical language is omitted and we only use the association $$\Rep_H\;\;\longmapsto\;\; \VB^{\nabla}_{\Omega^\bullet(B)}.$$ An interesting direction for future research is to study gauge theories on quantum flag manifolds (\cite{porque2,porque3}) for every qpc, as in references \cite{saldym,saldcon}. This will be addressed in forthcoming publications.

On the other hand, in Differential Geometry, connections on vector bundles and connections on principal bundles can be considered as sections of jets bundles. In reference \cite{porque}, the authors present an study of jets bundles in non--commutative geometry for a functorial point of view. An interesting line of research is to reproduce the previous \emph{classical} fact in non--commutative geometry using reference \cite{porque} and the line of research opens in this paper.  

\section{Examples}

In this section, we present two classes of qpb’s for which the quantum association functor $\A^{\omega^c}_{\zeta}$ can be defined. Of course, these classes of qpb’s are not the only ones for which the theory applies. For example, trivial quantum principal bundles (in the sense of Section~6.3 of reference \cite{micho2}) provide another class of examples for our theory.

\subsection{Quantum Principal $U(1)$--Bundles over the Non--Commutative $n$--Torus}

Fix a real antisymmetric $n\times n$ matrix $\Xi=(\Xi_{kj})$. Now, let us consider the quantum $n$--torus (\cite{con}) 
\begin{equation}
    \label{ec.5.1}
    \torus^n_\Xi.
\end{equation}
for $n$ $\in$ $\N$, $n\geq 2$. A generic element of $\torus^n_\Xi$ is  a formal sum $$\sum a_{m_1\cdots m_n} \u^{m_1}_1\,\u^{m_2}_2\cdots \u^{m_n}_n,$$ where $\{ a_{m_1\cdots m_n}\}_{m_i\in \Z}$ form a rapid decay sequence in $\C$ (decay faster than the inverse of any polynomial in $(m_1,...,m_n)$), and the operators  $\u_1$,..., $\u_n$ satisfy
\begin{equation}
    \label{ec.5.2}
    \u^\ast_k=\u^{-1}_k,\qquad \u_k\,\u_j=\mathrm{e}^{2\pi i\,\Xi_{kj}}\,\u_j\,\u_k
\end{equation}
for all $k$, $j$ $\in$ $\{1,...,n\}$. By taking the space $C(\torus^n_\Xi)$ defined as the universal $C^\ast$--algebra of $\torus^n_\Xi$, it is well--known that $\torus^n_\Xi$  is stable under holomorphic calculus \cite{con}. 

Now, let us consider the canonical $\ast$--Hopf algebra 
\begin{equation}
    \label{ec.5.3}
    (\mathcal{U}(1):=\C[z,z^\ast=z^{-1}],m,\mathbbm{1},\Delta,\epsilon,S,\ast)
\end{equation}
associated with the Lie group $U(1)$ (see Example 1.2 of reference \cite{appendix}). The coproduct, the counit and the antipode are given by
\begin{equation}
\label{3.f4.1.1}
\Delta(z)=z\otimes z,\quad \epsilon(z)=1, \quad S(z)=z^\ast, \quad S(z^\ast)=z.
\end{equation}

In order to develop a concrete example, we will focus in the \emph{usual} qvb over $\torus^n_\Xi$  (see Section 6.2 of reference \cite{tomas}). This is a quantum principal $\mathcal{U}(1)$--bundle 
\begin{equation}
    \label{ec.5.4}
    \zeta_n=(P:=\torus^{n+1}_\Xi,\torus^n_\Xi,\Delta_P),
\end{equation}
where  $$\Delta_P:P\longrightarrow P\otimes \mathcal{U}(1)$$ is the $\ast$--algebra morphism defined by
\begin{equation}
    \label{ec.5.5}
    \Delta_P(\u_i)=\u_i\otimes \mathbbm{1},\qquad \Delta_P(\u_{n+1})=\u_{n+1}\otimes z
\end{equation}
for $i=1,...,n$.

The next step is to define all the \emph{ingredients} of equation (\ref{ec.3.85}) in order to define a differential calculus on $\zeta_n$ as in Section 3.5.

Consider
\begin{equation}
    \label{ec.5.6}
      \delta_j: \torus^n_\Xi\longrightarrow \torus^n_\Xi
\end{equation}
the $\ast$--preserving derivation in $\torus^n_\Xi$ defined by $$\delta_j(\mathbbm{u}_k)=\left\{
\begin{array}{ll}
 2\pi i\, \mathbbm{u}_j & \text{ if } j=k \\
  \,\,0 & \text{ if } j \not= 0.
\end{array}
\right.$$ The derivation  $\delta_j$ is the non--commutative counterpart of the operator $\displaystyle {\partial \over \partial x_j }$ in the {\it classical} case (\cite{con}). Now, we can define a graded differential $\ast$--algebra as follows: consider the space $\C^n$ and its exterior algebra $$\bigwedge^\bullet \C^{n}=\bigoplus^n_{k=0}\wedge^k \C^{n}.$$ Then, we define 
\begin{equation}
    \label{ec.5.7}
    \Omega^\bullet(\torus^n_\Xi):= \torus^n_\Xi\otimes \bigwedge^\bullet \C^{n}
\end{equation}
with the $\ast$--algebra structure given by 
\begin{equation}
    \label{ec.5.8}
    (x_1\otimes v_1)\cdot (x_2\otimes v_2)=x_1\cdot x_2\otimes v_1\wedge v_2,\qquad (x\otimes v)^\ast=x^\ast\otimes v^\ast
\end{equation}
for $x$ $\in$ $\torus^n_\Xi$, $v$ $\in$ $\C$, and the differential 
\begin{equation}
    \label{ec.5.9}
   d: \Omega^\bullet(\torus^n_\Xi)\longrightarrow \Omega^\bullet(\torus^n_\Xi),\qquad x\longmapsto \delta_1(x)\,d\mathbbm{U}_1+\cdots \delta_n(x)\,d\mathbbm{U}_{n}
\end{equation}
for $x$ $\in$ $\Omega^0(\torus^n_\Xi)=\torus^n_\Xi\cong \torus^n_\Xi\otimes \C,$ where
\begin{equation}
    \label{ec.5.10}
   d\mathbbm{U}_j:=-{i\over 2\pi}\u^\ast_jd\u_j=(0,...,0,\mathbbm{1},0...,0).
\end{equation}
Notice that the differential can also be viewed as
\begin{equation}
    \label{ec.5.11}
    d(x)= (\delta_1(x),\delta_2(x),\cdots \delta_n(x)).
\end{equation}
For $x\,d\mathbbm{U}_1\cdots d\mathbbm{U}_k$ $\in$ $\Omega^k(\torus^n_\Xi)$ with $x$ $\in$ $\torus^n_\Xi$, we have 
\begin{equation}
    \label{ec.5.12}
    d(x\,d\mathbbm{U}_1\cdots d\mathbbm{U}_k)=dx\cdot d\mathbbm{U}_1\cdots d\mathbbm{U}_k.
\end{equation}
Of course, $(\Omega^\bullet(\torus^n_\Xi),d,\ast)$ is a graded differential $\ast$--algebra generated by its degree 0 elements $\Omega^0(\torus^n_\Xi)=\torus^n_\Xi$ \cite{con}.

In this way, consider the graded $\ast$--algebra 
\begin{equation}
    \label{ec.5.13}
(\Omega^\bullet_\Hor:=\Omega^\bullet(\torus^n_\Xi)\,P,m,\mathbbm{1},\ast).
\end{equation}
Of course $$(\Omega^\bullet_B:=\Omega^\bullet(\torus^n_\Xi),m,\mathbbm{1},\ast)$$ is a graded $\ast$--subalgebra with structure of graded differential $\ast$--algebra. Moreover, we define
\begin{equation}
    \label{ec.5.14}
    \Delta_\H:\Omega^\bullet_\Hor\longrightarrow \Omega^\bullet_\Hor\otimes \U(1) \quad \mbox{ given by}\quad \Delta_\H(\mu\,p)=\mu\,p^{(0)}\otimes p^{(1)},
\end{equation}
where $\Delta_P(p)=p^{(0)}\otimes p^{(1)}$ with $p$ $\in$ $P$. Clearly $\Omega^\bullet_B$ is exactly the space of all $\Delta_\H$--invariant elements and extends $\Delta_P$.

Let
\begin{equation}
    \label{ec.5.15}
    \mu=\sum^n_{j=1}t_j\,d\mathbbm{U}_j\qquad \mbox{ with }\qquad t_j \,\in\,\R.
\end{equation}
Notice that $\mu$ is a close element (i.e., $d\mu=0$) of the graded centre of $\Omega^\bullet_\Hor$.  Then, the linear map 
\begin{equation}
\label{ec.5.16}
    D_\mu: \Omega^\bullet_\Hor\longrightarrow \Omega^\bullet_\Hor
\end{equation}
such that $$D_\mu|_{\Omega^\bullet_B}=d,\qquad D_\mu(\u_{n+1})=-\u_{n+1}\,\mu=-\mu\,\u_{n+1}$$ and $$D_\mu(\varphi\,\psi)=D_\mu(\varphi)\,\psi+(-1)^k\,\varphi\,D_\mu(\psi)$$ for $\varphi$ $\in$ $\Omega^k_\Hor$, $\psi$ $\in$ $\Omega^\bullet_\Hor$, is well--defined and is an element of $\mathfrak{Der}$. In other words,  we have $\{ D_\mu\} \subseteq \mathfrak{Der}$ and clearly, $\{ D_\mu\}$ is in bijection with $\R^n$. 

Each operator $D_\mu$ satisfies
\begin{equation}
    \label{ec.5.17}
    D^2_\mu=0.
\end{equation}
Fix an operator $D_\mu$. In this way, for every bicovariant $\ast$--FODC on $\U(1)$
\begin{equation}
    \label{ec.5.18}
    (\Gamma,d) \qquad \mbox{ with }\qquad \mathfrak{qu}(1)^\#={\Ker(\epsilon)\over \mathcal{R}}, 
\end{equation}
the linear map $$ R_\mu:\mathfrak{qu}(1)^\#\longrightarrow \Omega^2_\Hor$$ of equation (\ref{curv1}) can be defined by setting $R_\mu=0$. Therefore, we can obtain a differential calculus on $\zeta_n$ by means of the method of Section 3.5 in such a way that
\begin{equation}
    \label{ec.5.19}
    \omega^c: \mathfrak{qu}(1)^\#\longrightarrow \widetilde{\Omega}^1(P),\qquad \theta\longrightarrow \mathbbm{1}\otimes \theta.
\end{equation}
is a regular and multiplicative qpc, and $D_\mu$ is its covariant derivative.

In summary, for each operator $D_\mu$ and each bicovariant $\ast$--FODC on $H$, there exists a differential calculus on $\zeta_n$ available for our theory. 

In order to develop a concrete example, let us take the \emph{canonical} choices, i.e., $$D:=D_0$$ and $$(\Gamma_\class,d)$$ be the \emph{classical} bicovariant $\ast$--FODC of differential forms on $\U(1)$  (see Example 1.2 of reference \cite{appendix}). In this case, the quantum dual Lie algebra is the complexification of the dual space of the Lie algebra  $\mathfrak{u}(1)$ of $U(1)$, i.e.,
\begin{equation}
    \label{ec.5.20}
    \mathfrak{qu}(1)^\#=\mathfrak{u}(1)^\#_\C=\mathrm{span}_\C\{ \pi(z)\},
\end{equation}
where $\pi: \U(1)\longrightarrow \mathfrak{u}(1)^\#_\C$ is the quantum germs map. In addition, $\pi(g)$ is the differential at the neutral element $e$ $\in$ $U(1)$ of the $\C$--valued function $g$ $\in$ $\U(1)$ (see \cite{appendix}).

Under these \emph{canonical} choices, the differential calculus on $\zeta_n$ obtained is isomorphic to the differential calculus on $\zeta_n$ for which the space of quantum differential forms of $P$ is
\begin{equation}
    \label{ec.5.21}
    \Omega^\bullet(P):=\Omega^\bullet(\torus^{n+1}_\Xi),
\end{equation}
where
\begin{equation}
    \label{ec.5.22}
\Delta_{\Omega^\bullet(P)}:\Omega^\bullet(P)\longrightarrow  \Omega^\bullet(P)\otimes \Gamma^\wedge
\end{equation}
is given by 
\begin{equation*}
    \Delta_{\Omega^\bullet(P)}|_{P}=\Delta_P,  \qquad\qquad \Delta_{\Omega^\bullet(P)}(d\u_i)=d\u_i\otimes \mathbbm{1},
\end{equation*}
$$ \Delta_{\Omega^\bullet(P)}(d\u_{n+1})=d\u_{n+1}\otimes z+\u_{n+1}\otimes dz. $$
Furthermore, in this differential calculus, the linear map $$\omega_\mu:\mathfrak{u}(1)^\#_\C\longrightarrow \Omega^1(P) \qquad \mbox{ such that }\qquad \omega_\mu(\pi(z))=\u^\ast_{n+1}d\u_{n+1}+\mu $$ is a regular and multiplicative qpc, and its covariant derivative is the operator $D_\mu$. Hence, constructing another differential calculus on $\zeta_n$ by the method of Section~3.5 for $D_\mu$ and $(\Gamma_\class,d)$ amounts to a \emph{translation} from $\omega_\mu$ to $\omega^c=\omega_0$. This illustrates the concept of \emph{translation} proposed at the end of Section 3.5.

Let us return to the differential calculus on $\zeta_n$ given by the canonical choices $D$, $(\Gamma_\class,d)$. In this way, in light of Theorem \ref{teo2},
\begin{equation}
    \label{ec.5.23}
    ((\Gamma_\class,d),\zeta_n,\omega^c=\omega_0)
\end{equation}
is equivalent to the quantum association functor 
\begin{equation}
    \label{ec.5.24}
    \A^{\omega^c}_{\zeta_n}: \Rep_{\U(1)}\longrightarrow \VB^\nabla_{\Omega^\bullet(\torus^{n}_\Xi)}.
\end{equation}

It is well--known that a complete set of mutually non--equivalent irreducible (finite--dimensional) $\U(1)$--corepresentations is given by
\begin{equation}
    \label{ec.5.25}
    \T:=\{\delta^m: \C\longrightarrow \C\otimes \U(1)\}_{m\in \Z},
\end{equation}
where $$\delta^m: \C\longrightarrow \C\otimes \U(1),\qquad w\longrightarrow w\otimes z^m.$$ In light of Proposition 6.3 of reference \cite{saldcon},  for each $m$ $\in$ $\Z$, the set $\{T_m\}$ with 
\begin{equation}
    \label{ec.5.26}
    T_m: \C\longrightarrow P, \qquad w\longmapsto w\,\u^m_{n+1}
\end{equation}
is exactly the set of operators of Theorem \ref{gen} for the qpb $\zeta_n$. Thus, we conclude that the associated qvb $$E^m:=\Mor(\delta^m,\Delta_P)$$ is  a $1$--dimensional left/right free $\torus^n_\Xi$--module.

Let $\delta^m$ $\in$ $\T$ and consider its associated qvb $E^m$. Then, for every $T$ $\in$ $E^m$ we have
\begin{equation}
    \label{ec.5.27}
    T=b^{_T}\,T_m\qquad \mbox{ with }\qquad b^{_T}=T(1)\,\u^{m\,\ast}_{n+1} \,\in\,\torus^n_\Xi.
\end{equation}
Thus $$D^{\omega^c}(T_m(1))=D(T_m(1))=D(\u^m_{n+1})=0\;\;\;\Longrightarrow \;\;\; D^{\omega^c}(T)=db^{_T}\,T_m$$ and hence 
\begin{equation}
    \label{ec.5.28}
    \nabla^{\omega^c}_m(T)=db^{_T}\otimes_{\torus^n_\Xi} T_m.
\end{equation}
Notice that 
\begin{equation}
    \label{ec.5.29}
    \nabla^{\omega^c}_m(T_m)=0.
\end{equation}
The functor $\A^{\omega^c}_{\zeta_n}$ is completely characterized by equations (\ref{ec.5.27}), (\ref{ec.5.28}) since every finite--dimensional $\U(1)$--corepresentation is the finite direct sum of elements of $\T$ and $\A^{\omega^c}_{\zeta_n}$ preserves the direct sum.

With the previous characterization of the functor $\A^{\omega^c}_{\zeta_n}$, it should be clear how to recreate $((\Gamma_\class,d),\zeta_n,\omega^c=\omega_0)$ using the methods presented in the text. For example, notice that equation (\ref{f44}) holds: $$P=\torus^{n+1}_\Xi\cong \bigoplus_{m\in \Z} E^m\otimes \C\cong \torus^n_\Xi\oplus \bigoplus_{m\not=0} E^m\otimes \C\cong \torus^n_\Xi\oplus \bigoplus_{m\not=0}  \torus^n_\Xi\{\u^m_{n+1} \}.$$ As another example, according to equation (\ref{ec.5.28}), the operator $$D_m: \E_m\longrightarrow \E_m $$ of equation (\ref{f46}) is given by $$D_m(\tau)=d\mu^\tau\,T_m\qquad \mbox{ where }\qquad \mu^\tau=\tau(1) \,\u^{m\,\ast}_{n+1}\,\in\,\Omega^\bullet(\torus^n_\Xi)$$ for all $\tau$ $\in$ $\E_m=\Mor(\delta_m,\Delta_\Hor)$, and \emph{gluing} of all them in $$\Omega^\bullet_\Hor=\Omega^\bullet(\torus^n_\Xi)\,P\,\cong \bigoplus_{m\in \Z} \E^m\otimes \C\cong \Omega^\bullet(\torus^n_\Xi)\oplus \bigoplus_{m\not=0} \E^m\otimes \C\cong \Omega^\bullet(\torus^n_\Xi)\oplus \bigoplus_{m\not=0}   \Omega^\bullet(\torus^n_\Xi)\{\u^m_{n+1} \}$$ we recover an operator on  $\Omega^\bullet_\Hor$ such that in $\Omega^\bullet(\torus^n_\Xi)$ is $d$, satisfies the graded Leibniz rule and in $\u_{n+1}$ is zero (see equation (\ref{ec.5.29}) for $m=1$). In other words, we recover the operator $D=D_0$.

Thus,
\begin{equation}
    \label{ec.5.30}
    \{ \zeta_n\mid n \,\in\,\N,\;n\geq 2 \}
\end{equation}
is a collection of qpb's suitable for the theory of this paper. 

\subsection{Homogeneous Quantum Principal Bundles}

Homogeneous quantum principal bundles are one of the most well--studied examples of qpb's and the reader can check the basics in, for example, references  \cite{micho2,stheve,libro}. Two relevant references providing a modern perspective on these spaces are given in \cite{porque1,raemon0}.

It is worth mentioning that Remark \ref{rema} holds for these qpb's provided that Remark \ref{hopf1} is satisfied. In particular, in Proposition $5.3$ of reference \cite{sald2} it is shown the specific form of the maps $\{T^V_k \}$ in this case.

Let  $$(P,m,\mathbbm{1},\Delta,\epsilon,S,\ast)$$ be a $\ast$--Hopf algebra and consider   $$(H,m,\mathbbm{1},\Delta',\epsilon',S',\ast) $$  be a  $\ast$--Hopf subalgebra. This implies the existence of a $\ast$--Hopf algebra epimorphism $$j:P\longrightarrow H$$ and consider the linear map $$\Delta_P:=(\id_P\otimes j)\circ \Delta: P\longrightarrow P\otimes H.$$ Defining $$B:=\{b \in P \mid \Delta_P(b)=b\otimes \mathbbm{1} \},$$ the triple $$\zeta=(P,B,\Delta_P)$$ is a qpb called {\it homogeneous quantum principal bundle}. 

Let 
\begin{equation}
    \label{ec.5.31}
    (\Omega^1_\Hor,d)
\end{equation}
 be a right covariant $\ast$--FODC on $P$ and let us take the universal differential envelope $\ast$--calculus
\begin{equation}
    \label{ec.5.32}
    (\Omega^\bullet_\Hor,d,\ast)
\end{equation}
of $(\Omega^1_\Hor,d)$, and consider the linear map $${_{\Omega^{1\wedge}_\Hor}}\Phi: \Omega^\bullet_\Hor\longrightarrow \Omega^\bullet_\Hor\otimes P $$ of equation (\ref{ec.3.13}) (which exists for right covariant $\ast$--FODC's \cite{stheve}); which is a graded differential $\ast$--algebra, extends $\Delta$ and satisfies 
\begin{equation}
    \label{ec.5.33}
    {_{\Omega^{1\wedge}_\Hor}}\Phi\circ d=(d\otimes \id_P)\circ \Delta.
\end{equation}
Hence, we define
\begin{equation}
    \label{ec.5.34}
    \Delta_\H:=(\id_{\Omega^\bullet_\Hor}\otimes j)\circ {_{\Omega^{1\wedge}_\Hor}}\Phi:\Omega^\bullet_\Hor\longrightarrow \Omega^\bullet_\Hor\otimes H.
\end{equation}
This map is a graded differential $\ast$--algebra which extends $\Delta_P$ and we can define
\begin{equation}
    \label{ec.5.35}
    \Omega^\bullet_B:=\{\mu\,\in\,\Omega^\bullet_\Hor\mid \Delta_\Hor(\mu)=\mu\otimes \mathbbm{1}\}.
\end{equation}
This is a graded differential $\ast$--subalgebra of $(\Omega^\bullet_\Hor,d,\ast)$. In addition, by equation (\ref{ec.5.33}), it follows that $$\Delta_\H\circ d=(d\otimes \id_H)\circ \Delta_\H.$$ Therefore, $d$ $\in$ $\mathfrak{Der}$.

In this way, for every bicovariant $\ast$--FODC on $H$ 
\begin{equation}
    \label{ec.5.36}
    (\Gamma,d) \qquad \mbox{ with }\qquad \mathfrak{qg}^\#={\Ker(\epsilon')\over \mathcal{R}'},
\end{equation}
the linear map  $$R:\mathfrak{qg}^\#\longrightarrow \Omega^2_\Hor$$ of equation (\ref{curv1}) can be defined by setting $R=0$, since $d^2=0$. Hence, we can obtain a differential calculus on $\zeta$ by means of the method of Section 3.5 in such a way that
\begin{equation}
    \label{ec.5.37}
    \omega^c:\mathfrak{qg}^\#\longrightarrow \widetilde{\Omega}(P),\qquad \theta \longrightarrow \mathbbm{1} \otimes \theta
\end{equation}
is a regular and multiplicative qpc, and $d$ is its covariant derivative.

With this differential structure on $\zeta$, it follows from Theorem \ref{teo2} that
\begin{equation}
    \label{ec.5.38}
    ((\Gamma,d),\zeta,\omega^c)
\end{equation}
is equivalent to the quantum association functor
\begin{equation}
    \label{ec.5.39}
    \A^{\omega}_\zeta: \Rep_H\longrightarrow \VB^\nabla_{\Omega^\bullet_B}.
\end{equation}
Therefore, homogeneous quantum principal bundles forms a collection of qpb's suitable for the theory of this paper.

However, in certain cases, it is possible to choose another space $\Omega^\bullet_\Hor$. In fact, let $\zeta$ be a homogeneous qpb and consider a left covariant $\ast$--FODC on $P$
\begin{equation}
    \label{ec.5.40}
    (\Omega^1(P),d) \qquad \mbox{ with }\qquad \mathfrak{qp}^\#={\Ker(\epsilon)\over \mathcal{R}},
\end{equation}
and consider the space
\begin{equation}
\label{ec.5.41}
\mathcal{R}':=j(\mathcal{R}).
\end{equation}
Then, $\mathcal{R}'$ defines a $\ast$--FODC on $H$
\begin{equation}
    \label{ec.5.42}
    (\Gamma',d)\qquad \mbox{ with }\qquad \mathfrak{qg}^\#={\Ker(\epsilon')\over \mathcal{R}'}.
\end{equation}
Therefore, in accordance with Proposition 2.2 of reference \cite{raemon0}, this defines a quantum principal bundle in the sense of Section 5 of reference \cite{libro} if and only if
\begin{equation}
    \label{ec.5.43}
    (\id\otimes j)\Ad(\mathcal{R})\subseteq \mathcal{R}\otimes H.
\end{equation}
In particular, in this situation, $(\Gamma',d)$ is bicovariant. According to Section 5.4 of reference \cite{libro}, the map $\Delta_P$ can be extended to a $\ast$--preserving morphism
$$\Delta_{\Omega^1(P)}: \Omega^1(P)\longrightarrow (\Omega^1(P)\otimes H)\oplus (P\otimes \Gamma')$$  such that $$\Delta_{\Omega^1(P)}(p_1dp_2)=\Delta_P(p_1)(dp^{(0)}_2\otimes p^{(1)}_2+p^{(0)}\otimes dp^{(1)})$$ with $\Delta_P(p_2)=p^{(0)}_2\otimes p^{(1)}_2$.

Consider now the universal differential envelope $\ast$--calculus 
\begin{equation}
    \label{ec.5.44}
    (\Omega^\bullet(P),d,\ast)
\end{equation}
of $(\Omega^1(P),d)$, and the universal differential envelope $\ast$--calculus 
\begin{equation}
    \label{ec.5.45}
    (\Gamma'^\wedge,d,\ast)
\end{equation}
of $(\Gamma',d)$. By Proposition \ref{prop3.1}, the maps $\Delta_P$, $\Delta_{\Omega^1(P)}$ can be extended to a graded differential $\ast$--algebra morphism
\begin{equation}
    \label{ec.5.46}
    \Delta_{\Omega^\bullet(P)}:\Omega^\bullet(P)\longrightarrow \Omega^\bullet(P)\otimes \Gamma'^\wedge
\end{equation}
and we get a differential calculus on $\zeta$. Let us define
\begin{equation}
    \label{ec.5.47}
    \Omega^\bullet_\Hor:=\Hor^\bullet\,P:=\{\varphi\in\Omega^\bullet(P)\mid \Delta_{\Omega^\bullet(P)}(\varphi)\in \Omega^\bullet(P)\otimes H \},
\end{equation}
\begin{equation}
    \label{ec.5.48}
\Delta_\H:=\Delta_\Hor:=\Delta_{\Omega^\bullet(P)}|_{\Hor^\bullet\,P},
\end{equation}
\begin{equation}
    \label{ec.5.49}
    \Omega^\bullet_B:=\Omega^\bullet(B):=\{\mu \in \Omega^\bullet(P)\mid  \Delta_{\Omega^\bullet(P)}(\mu)=\mu\otimes \mathbbm{1}\}.
\end{equation}

On the other hand, equations (\ref{ec.5.42}), (\ref{ec.5.43}) guarantees the existence of a projection map
\begin{equation}
    \label{ec.5.50}
    \rho: \mathfrak{qp}^\#\longrightarrow \mathfrak{qg}^\#
\end{equation}
defined by the formula
\begin{equation}
    \label{ec.5.51}
    \rho\circ \pi=\pi'\circ j,
\end{equation}
where $$\pi:P\longrightarrow \mathfrak{qp}^\#$$ is the quantum germs map associated with $(\Omega^1(P),d)$ (see equation (\ref{ec.3.4})), and $$\pi':H\longrightarrow \mathfrak{qg}^\#$$ is the quantum germs map associated with $(\Gamma',d)$ (see equation (\ref{ec.3.4})). In particular,   
$\Ker(\rho)$ is a $\ast$--$\ad$--invariant right $P$--submodule of $\mathfrak{qp}^\#$. Fix now a $\ast$--$\ad$--invariant complement $L$ of $\Ker(\rho)$ in $\mathfrak{qp}^\#$.  By Lemma 6.14 of reference \cite{micho2}, the linear map
\begin{equation}
    \label{ec.5.52}
    \omega':\mathfrak{qg}^\#\longrightarrow \Omega^1(P),\qquad \theta\longmapsto \rho|_L^{-1}(\theta)
\end{equation}
is a quantum principal connection.

Assume that there exists at least, one regular qpc $\widehat{\omega}$ (in some cases, $\omega'$ is regular). According to Section 4 of reference \cite{micho2}, for a given regular qpc $\omega$, there always exists a non--trivial quotient differential calculus for which $\omega$ becomes multiplicative. Hence, we may assume that $\widehat{\omega}$ is also multiplicative. In light of Section 3.3, the covariant derivative $$D:=D^{\widehat{\omega}}$$ of $\widehat{\omega}$ is an element of $\mathfrak{Der}$. In addition, since $\widehat{\omega}$ is multiplicative, the map $R$ of equation~(\ref{curv1}) can be defined by setting $R:=R^{\widehat{\omega}}$. Hence, we can construct a differential calculus on $\zeta$ using the method of Section~3.5 in such a way that $\widehat{\omega}$ is \emph{sent} to $\omega^c$ in the sense that $D$ is the covariant derivative of $\omega^c$.

With this differential calculus on $\zeta$, it follows from Theorem \ref{teo2} that
\begin{equation}
    \label{ec.5.53.1}
    ((\Gamma',d),\zeta,\omega^c)
\end{equation}
is equivalent to the quantum association functor
\begin{equation}
    \label{ec.5.54}
    \A^{\omega^c}_\zeta: \Rep_H\longrightarrow \VB^\nabla_{\Omega^\bullet_B}.
\end{equation}

To finalize this paper, let us present a concrete example for which $\omega'$ is regular and multiplicative.

Let us take the $\ast$--Hopf algebra $\SU_q(2)$ of the quantum $SU(2)$ group for $q$ $\in$  $(-1,1)-\{0\}$. We will use the original Woronowicz's notation for $\SU_q(2)$ presented in \cite{woro1}, that is, the $\ast$--algebra $$(\SU_q(2),\cdot,\mathbbm{1},\ast)$$ is generated by two symbols $\{\alpha,\gamma\}$ satisfying
\begin{equation}
\label{2.f1.6}
\begin{aligned}
\alpha^{\ast}\alpha+\gamma^{\ast}\gamma=\mathbbm{1},\qquad \alpha\alpha^{\ast}+q^{2}\gamma\gamma^{\ast}=\mathbbm{1},\qquad \gamma\gamma^{\ast}=\gamma^{\ast}\gamma\\
q\gamma\alpha=\alpha\gamma,\quad q\alpha^\ast\gamma^\ast=\gamma^\ast\alpha^\ast, \quad q\gamma^{\ast}\alpha=\alpha\gamma^{\ast},\quad q\alpha^\ast\gamma=\gamma\alpha^\ast
\end{aligned}
\end{equation}
 and the following relations for the coproduct, the counit and the antipode:   
\begin{equation}
\label{2.f1.7}
\begin{split}
\Delta(\alpha)&=\alpha\otimes\alpha-q\gamma^\ast\otimes\gamma, \quad \Delta(\gamma)=\gamma\otimes\alpha+\alpha^\ast\otimes\gamma,\quad \epsilon(\alpha)=1, \quad \epsilon(\gamma)=0\\
&S(\alpha)=\alpha^\ast, \qquad S(\alpha^\ast)=\alpha,\qquad S(\gamma)=-q\gamma,\qquad S(\gamma^\ast)=-q^{-1}\gamma^\ast.
\end{split}
\end{equation}

Now, let us take the canonical $\ast$--Hopf algebra associated with the Lie group $U(1)$ (see equations (\ref{ec.5.3}), (\ref{3.f4.1.1})), and the $\ast$--Hopf algebra epimorphism
\begin{equation}
\label{3.f4.2.1}
j:\SU_q(2) \longrightarrow \U(1)
\end{equation}
such that $$j(\alpha)=z,\qquad j(\gamma)=0.$$  By considering the linear map
\begin{equation}
\label{3.f4.3}
\Delta{_{\SU_q(2)}}:=(\id_{\SU_q(2)}\otimes j)\circ \Delta:\SU_q(2)\longrightarrow \SU_q(2)\otimes \U(1).
\end{equation}
and the space $$\S^2_q:=\{ b\in \SU_q(2) \mid \Delta{_{\SU_q(2)}}(b)=b\otimes \mathbbm{1} \},$$  we get a homogeneous quantum principal $\mathcal{U}(1)$--bundle 
\begin{equation}
    \label{ec.hopf}
    \zeta_{\mathrm{Hopf}}=(\SU_q(2),\S^2_q,\Delta{_{\SU_q(2)}})
\end{equation}
commonly called  the {\it quantum Hopf fibration} or the $q$–{\it Dirac monopole bundle} \cite{micho2,stheve,libro}.

Consider now the {\it $3D$ Woronowicz differential calculus of} $\SU_q(2)$ (\cite{woro2,micho1}). This is a left covariant $\ast$--FODC
\begin{equation}
\label{3.f4.5}
(\Omega^1(\SU_q(2)),d)
\end{equation}
defined by
\begin{equation}
\label{3.f4.6.1}
\mathcal{R}_3:=\langle \{ \gamma^2,\,\gamma^{\ast\,2},\,\gamma\gamma^\ast,\,\alpha\gamma-\gamma,\, \alpha\gamma^\ast-\gamma^\ast,\,q^2\alpha+\alpha^\ast-(1+q^2)\mathbbm{1} \}\rangle\; \subseteq \; \Ker(\epsilon).
\end{equation}
Here, $\langle x\rangle$ denotes the right $\SU_q(2)$--ideal generated by $x$. The $\ast$--FODC $(\Omega^1(\SU_q(2)),d)$ receives this name because, as Woronowicz showed in \cite{woro2}, the $\C$--vector space $$\mathfrak{qsu^\#}(2):=\dfrac{\Ker(\epsilon)}{\mathcal{R}_3}$$ has dimension $3$ and the set
\begin{equation}
\label{3.f4.7}
\beta:=\{\eta_3:=\pi(\alpha-\alpha^\ast),\quad \eta_+:=\pi(\gamma),\quad \eta_-:=\pi(\gamma^\ast) \}
\end{equation}
is a linear basis of $\mathfrak{qsu^\#}(2)$, where $\pi$ is the corresponding quantum germs map. Furthermore, $\beta$ is a left $\SU_q(2)$--basis of $\Omega^1(\SU_q(2))$ (\cite{woro2,stheve}), i.e., $$\Omega^1(\SU_q(2))=\SU_q(2)\,\eta_-+\SU_q(2)\,\eta_++\SU_q(2)\,\eta_3.$$

Now, let us consider the right $\U(1)$--ideal  
\begin{equation}
\label{3.f4.15}
\mathcal{R}':=j(\mathcal{R}_3) \subseteq \Ker(\epsilon').
\end{equation}
The $\ast$--FODC  $$(\Gamma',d)$$ induced by $\mathcal{R}'$ is bicovariant and if $\pi'$ is the corresponding quantum germs map, then 
\begin{equation}
\label{3.f4.16}
\beta':=\{\varsigma:=\pi'(z-z^\ast)\}
\end{equation}
is a linear basis of the $\C$--vector space $$\mathfrak{qu^\#}(1):=\dfrac{\Ker(\epsilon')}{\mathcal{R}'}.$$ According to Section 6.4 of reference \cite{micho2}, $(\Gamma',d)$ is not the \emph{classical} $\ast$--FODC on $\U(1)$ of the previous subsection. Finally, since equation (\ref{ec.5.44}) is satisfied, the map $\Delta_{\SU_q(2)}$ can be extended to $\Delta_{\Omega^1(\SU_q(2))}$. In particular, we have
$$\Delta_{\Omega^1(\SU_q(2))}(\eta_-)=\eta_-\otimes z^{2\,\ast},\qquad \Delta_{\Omega^1(\SU_q(2))}(\eta_+)=\eta_+\otimes z^{2}$$ and $$\Delta_{\Omega^1(\SU_q(2))}(\eta_3)=\eta_3\otimes \mathbbm{1}+\mathbbm{1}\otimes \varsigma.$$

Taking the universal differential envelope $\ast$--calculus $$(\Omega^\bullet(\SU_q(2)),d,\ast)$$ of $(\Omega^1(\SU(2)),d)$, and the universal differential envelope $\ast$--calculus
$$(\Gamma'^\wedge,d,\ast)$$ of $(\Gamma',d)$, by Proposition  \ref{prop3.1} we get the map $$ \Delta_{\Omega^\bullet(\SU_q(2))}:\Omega^\bullet(\SU_q(2))\longrightarrow \Omega^\bullet(\SU_q(2))\otimes \Gamma'^\wedge$$ and this defines a differential calculus on $\zeta_{\mathrm{Hopf}}$.

According to \cite{saldhopf,tomas}, we have that
$$\Hor^\bullet\,\SU_q(2)=\bigoplus^2_{k=0}\Hor^k\,\SU_q(2),$$ with $$\Hor^0\,\SU_q(2)=\SU_q(2),\qquad \Hor^1\,\SU_q(2)=\SU_q(2)\,\eta_-+ \SU_q(2)\,\eta,$$ $$\Hor^2\,\SU_q(2)=\SU_q(2)\,\eta_-\,\eta_+ $$ and $$\Omega^\bullet(\S^2_q)=\bigoplus^2_{k=0}\Omega^k(\S^2_q)$$ with $$\Omega^0(\S^2_q)=\S^2_q,\qquad \Omega^2(\S^2_q)=\S^2_q\,\eta_-\,\eta_+,$$ $$\Omega^1(\S^2_q)=\{x\eta_-+y\eta_+\mid \Delta_{\SU_q(2)}(x)=x\otimes z^2,\;\;\Delta_{\SU_q(2)}(y)=y\otimes z^{2\,\ast} \}.$$

Finally, the projection map of equation (\ref{ec.5.51})
\begin{equation}
    \label{ec.5.65}
    \rho: \mathfrak{qsu^\#}(2)\longrightarrow\mathfrak{qu^\#}(1) 
\end{equation}
is given by $$\rho(\eta_-)=\rho(\eta_+)=0,\qquad \rho(\eta_3)=\varsigma$$ and $\Ker(\rho)=\mathrm{span}_\C\{\eta_-,\eta_+ \}$. There is a \emph{canonical} $\ast$--$\ad$--invariant complement $L$ of $\Ker(\rho)$, which is given by $$L=\mathrm{span}_\C\{\eta_3\} $$ and hence, the qpc $\omega'$ of equation (\ref{ec.5.53}) is given by 
\begin{equation}
    \label{ec.5.66}
    \omega'(\varsigma)=\eta_3,
\end{equation}
which is regular and multiplicative \cite{micho2}. It is worth mentioning that, according to Proposition 2.2 of reference \cite{saldhopf}, $\omega'$ is the only regular qpc. 

Finally, the covariant derivative of $ \omega'$ is the operator (\cite{saldhopf})
\begin{equation}
    \label{ec.5.67}
    D:=D^{\omega'}:\Hor^\bullet\,\SU_q(2)\longrightarrow \Hor^\bullet\,\SU_q(2)
\end{equation}
that satisfies the graded Leibniz rule, for $x$ $\in$ $\SU_q(2)$ is given by $$D(x)=x^{(1)}[\pi_-(x^{(2)})+\pi_+(x^{(2)})]$$  with $\Delta(x)=x^{(1)}\otimes x^{(2)}$, and $$D(\eta_-)=D(\eta_+)=0.$$ Here, $\pi_\pm:=\rho_\pm\circ \pi$, where $$\rho_\pm:\mathfrak{qsu^\#}(2)\longrightarrow \mathrm{span}_\C\{\eta_\pm\}$$ is the canonical projection. Explicitly $$D(\alpha)=-q\gamma^\ast\eta_+,\qquad D(\alpha^\ast)=-q\gamma\eta_-=D(\alpha)^\ast,$$  $$D(\gamma)=\alpha^\ast\eta_+,\qquad D(\gamma^\ast)=\alpha\eta_-=D(\gamma)^\ast,\qquad D(\mathbbm{1})=0.$$ In addition 
\begin{equation}
    \label{ec.5.68}
    D^2(\varphi)=-\varphi^{(0)}\,R(\pi(\varphi^{(1)})),
\end{equation}
where
\begin{equation}
    \label{ec.5.69}
    R:\mathfrak{qu^\#}(1)\longrightarrow \Omega^2(\SU_q(2))
\end{equation}
is given by (\cite{micho2}) $$R(\varsigma)=d\eta_3=q(1+q^2)\,\eta_-\,\eta_+.$$

For $$\Omega^\bullet_\Hor:=\Hor^\bullet\,\SU_q(2), \qquad \Omega^\bullet_{\S^2_q}:=\Omega^\bullet(\S^2_q),\qquad \Delta_\H:=\Delta_\Hor,$$ we obtain $D$ $\in$ $\mathfrak{Der}$ and the operator $R$ of equation (\ref{curv1}) can be defined. Hence,  we can apply the method presented in Section 3.5 to get another differential calculus on $\zeta_{\mathrm{Hopf}}$. It is worth mentioning that, according to \cite{michodif}, in this specific situation,  the resulting differential calculus is isomorphic to the original one and under this isomorphism, $\omega'$ turns into $\omega^c$. However, now it is clear that we can apply our theory to $$((\Gamma',d),\zeta_{\mathrm{Hopf}},\omega^c)$$ and the quantum association functor $$\A^{\omega^c}_{\zeta_{\mathrm{Hopf}}}: \Rep_{\U(1)}\longrightarrow \VB^\nabla_{\Omega^\bullet(\S^2_q)}.$$ For example, taking into account equation (\ref{ec.5.25}), is well--known that (\cite{30})
$$\SU_q(2)=\bigoplus_{n\in \Z} \{x\in \SU_q(2)\mid \Delta_{\SU_q(2)}(x)=x\otimes z^n \} \cong \bigoplus_{n\in \Z} E^n\otimes \C,$$ where $$E^n:=\Mor(\delta^n,\Delta_{\SU_q(2)})  $$ is the associated qvb of $\delta^n$.

It is worth mentioning that the same process can be carried out for any homogeneous qpb equipped with a differential calculus given by the condition of equation~(\ref{ec.5.43}), provided that there exists at least one regular qpc.



\appendix

\section{On the Classical Case}

In this appendix we present a brief summary of the theory presented in reference \cite{saldgreg}. Fix a smooth (not necessary compact) manifold $M$. Let $G$ be a (not necessary compact) Lie group and consider $\mathbf{MF}_G$ the category of (not necessary compact) manifolds $F$ endowed with a smooth left action $\star_F:G\times F\longrightarrow F$. Morphisms in this category are smooth $G$--equivariant maps. Furthermore, consider the category $\mathbf{FB}_M$ of fiber bundles $$\pi: FM\longrightarrow M $$ over $M$. Arrows in this category are fiber bundle morphisms. In this way, every principal $G$--bundle $$\pi:GM\longrightarrow M $$ over $M$ defines a covariant functor $$\A_{GM}:\mathbf{MF}_G\longrightarrow \mathbf{FB}_M$$ given by the associated fiber bundle. In other words, for every $\star_F$ $\in$ $\Obj(\mathbf{MF}_G)$,  its image under $\A_{GM}$ is the fiber bundle $$\pi_F: F^\star M\longrightarrow M, \qquad [x,f]\longmapsto \pi(x),$$ where $$F^\star M:=GM\times_G F:=(GM\times F)/G,$$ where the action of $G$ on $GM\times F$ is given by $$(x,f,A)\longmapsto (x\cdot g,\star_F(A^{-1},f)).$$ In addition, given a $\varphi$ $\in$ $\Mor(\star_F,\star_{F'})$, its image under $\A_{GM}$ is the fiber bundle morphism $$\A_{GM}(\varphi):F^\star M\longrightarrow F'^{\star}M, \qquad [x,f]\longmapsto [x,\varphi(f)].$$ 

The functor $\A_{GM}$ is called the \emph{association functor}. In addition, if $$\omega: TGM\longrightarrow \mathfrak{g} $$ is a principal connection of $\pi:GM\longrightarrow M$, then the association functor can be promoted to
\begin{equation}
    \label{a.1}
    \A^\omega_{GM}:\mathbf{MF}_G\longrightarrow \mathbf{FB}^\nabla_M,
\end{equation}
where $\mathbf{FB}^\nabla_M$ is the category of fiber bundles over $M$ with general connections and parallel fiber bundle morphisms.

The association functor $\A^\omega_{GM}$ can be restricted to
\begin{equation}
    \label{a.2}
    \A^\omega_{GM}:\mathbf{Rep}_G\longrightarrow \mathbf{VB}^\nabla_M,
\end{equation}
where $\mathbf{Rep}_G$ is the category of finite--dimensional linear representations of $G$ and $\mathbf{VB}^\nabla_M$ is the category of vector bundles over $M$ with linear connections.

Reference~\cite{saldgreg} presents a study of the functors defined in equations~(\ref{a.1}) and~(\ref{a.2}). A particularly important result is given in Proposition~4.6, where the image of the functor in equation~(\ref{a.2}) is characterized under the assumption that $G$ and $M$ are simply connected. 

However, the main results of~\cite{saldgreg} are presented in Section~5. First, Theorem~5.1 shows that, given a cartesian covariant functor $$\F:\mathbf{MF}_G\longrightarrow \mathbf{FB}^\nabla_M$$ which agrees with the product functor $$M\times: \mathbf{MF}\longrightarrow \mathbf{FB}^\nabla_M$$ in the subcategory $\mathbf{MF}$ of $\mathbf{MF}_G$ consisting of manifolds with trivial $G$--actions, there exists a principal $G$--bundle over $M$ such that the functor of equation~(\ref{a.1}) is naturally isomorphic to $\F$.

Now, let us define the following two
 rather special categories associated with a given smooth manifold $M$. Objects in
 the category $\mathbf{PB}^\omega_M$ of principal bundles with principal connections on
 $M$ are triples $$(G,GM,\omega)$$ formed by a Lie group $G$ and a
 principal $G$--bundle $\pi:GM\longrightarrow M$ over $M$ endowed with a principal connection
 $\omega$. Every morphism between two such objects
 $$
  (\varphi_{\mathrm{grp}},\varphi):(G,GM,\omega)
  \longrightarrow (\hat G,\hat GM,\hat\omega\;)
 $$
 consists of a parallel homomorphism $\varphi:GM\longrightarrow\hat GM$
 of fiber bundles which is $G$--equivariant over the Lie group homomorphism
 $\varphi_{\mathrm{grp}}:G\longrightarrow\hat G$. 

On the other hand, objects in the category
 $\mathbf{GTS}^\nabla_M$ of gauge theory sectors on $M$ are tuples $$(G,\F)$$ formed by a Lie group $G$ and a cartesian covariant functor $$\F:\mathbf{MF}_G\longrightarrow\mathbf{FB}^\nabla_M$$ which   agrees with the product functor on the subcategory
 $\mathbf{MF}\,\subset\,\mathbf{MF}_G$. In
 $\mathbf{GTS}^\nabla_M$, morphisms are again tuples
 $$
  (\varphi_{\mathrm{grp}},\Phi):
  (G,\F)\longrightarrow (\hat G,\hat\F)
 $$
 consisting of a Lie group homomorphism $\varphi_{\mathrm{grp}}:G\longrightarrow\hat G$
 and a natural transformation $$\Phi:\F\,\circ
 \,\varphi^*_{\mathrm{grp}}\longrightarrow \hat\F$$ between the two functors $\mathbf{MF}_{\hat G}
 \longrightarrow \mathbf{FB}^\nabla_M$ involved, where the action pull back functor
 $$
  \varphi^*_{\mathrm{grp}}:\mathbf{MF}_{\hat G}\longrightarrow \mathbf{MF}_G,\qquad
  (\hat F,\star_{\hat G})\longmapsto (\hat F,\star_G)
 $$
 induced by $\varphi_{\mathrm{grp}}$ lets $G$ act via $$g\,\star_Gf\,:=\,\varphi_{\mathrm{grp}}(g)
 \,\star_{\hat G}f$$ on a $\hat G$--manifold $\hat F$.

In this way,  we can consider the association functor as a functor
 $$
\A:\;\;\mathbf{PB}^\omega_M\;\longrightarrow\;\mathbf{GTS}^\nabla_M
 $$
 with $$(G,GM,\omega)\longmapsto(G,\A^\omega_{GM})$$ on objects, and for every morphism $(\varphi_{\mathrm{grp}},\varphi)$ in the
 source category $\mathbf{PB}^\omega_M$, the induced morphism of $\A$ in the target category $\mathbf{GTS}^\nabla_M$ is the natural transformation $\Phi_\varphi$ defined for
 $\hat F\,\in\,\Obj(\mathbf{MF}_{\hat G})$ by
 $$
  \Phi_\varphi(\,\hat F\,):\;\;
  GM\,\times_G\hat F\;\longrightarrow\;\hat GM\,\times_{\hat G}\hat F,
  \qquad[\,g,\,\hat f\,]\;\longmapsto\;[\,\varphi(g),\,\hat f\,]\ .
 $$ In this way, in Corollary 5.2 of \cite{saldgreg} it is proven that $\A$ is an equivalence of categories.

In particular, for compact spaces, it is possible to work only with the restricted association functor of equation~(\ref{a.2}). In this case, for every object $(G,\F)$ of $\mathbf{GTS}^\nabla_M$, the functor $\F$ is of the form $$\F=A^\omega_{GM}:\mathbf{Rep}_G\longrightarrow \mathbf{VB}^\nabla_M,$$ for some principal $G$--bundle $\pi:GM\longrightarrow M$. This is the way we present this category at the end of Section~4.2. As explained there, the name of the category $\mathbf{GTS}^\nabla_M$ reflects the fact that each of its objects determines a gauge theory on $M$. In this sense, Theorem~\ref{teo1} can be regarded as the non--commutative geometrical counterpart \emph{at purely algebraic level} of Theorem~5.1 in~\cite{saldgreg}, while Theorem~\ref{teo2} corresponds to the non--commutative geometrical counterpart \emph{at purely algebraic level} level of Corollary~5.2 in~\cite{saldgreg}.

\section{Some Definitions in Category Theory}
In this appendix, we present some definitions in category theory used throughout the paper. This appendix is based on reference \cite{barcategories}.

\begin{Definition}
\label{b.1}
A monoidal category is a tuple 
$(\mathbf{C},\otimes,\mathbbm{1}_{\mathbf{C}},a,l,r)$ consisting of

\begin{enumerate}
\item A category $\mathbf{C}$.

\item A bifunctor
\[
\otimes : \mathbf{C}\times\mathbf{C} \longrightarrow \mathbf{C}
\]
called the \emph{tensor product}.

\item An object $\mathbbm{1}_{\mathbf{C}}$ $\in$ $\Obj(\mathbf{C})$ called the \emph{unit object}.

\item A natural isomorphism (the \emph{associator})
\[
a_{X,Y,Z}:(X\otimes Y)\otimes Z \longrightarrow X\otimes(Y\otimes Z)
\]
for all $X,Y,Z$ $\in$ $\Obj(\mathbf{C})$.
\item Natural isomorphisms (the \emph{unit constraints})
\[
l_X:\mathbbm{1}_{\mathbf{C}}\otimes X \longrightarrow X,
\qquad
r_X:X\otimes \mathbbm{1}_{\mathbf{C}} \longrightarrow X
\] for all objects $X,Y,Z$ $\in$ $\Obj(\mathbf{C})$, satisfying the following coherence conditions:
\begin{enumerate}

\item The pentagon axiom diagram is satisfies (see Section 2 of reference \cite{librogroups}).

\item  The triangle axiom  diagram is satisfies (see Section 2 of reference \cite{librogroups}).

\end{enumerate}
\end{enumerate}
\end{Definition}

A monoidal category $\mathbf{C}$  is called \emph{strict} if $a$, $l$ and $r$ are the identities maps for all objects of $\mathbf{C}$. It is worth mentioning that every monoidal category can be \emph{strictified}.

\begin{Definition}
\label{b.2}
A monoidal additive category is a monoidal category 
$(\mathbf{C},\otimes,\mathbbm{1}_{\mathbf{C}},a,l,r)$ such that:
\begin{enumerate}
\item The category $\mathbf{C}$ is an additive category; that is:
\begin{enumerate}
\item For any objects $X,Y \in \Obj(\mathbf{C})$, the set 
$\Mor(X,Y)$ is an abelian group.
\item Composition of morphisms is bilinear with respect to this group structure.
\item The category admits a zero object and finite direct sums $$\oplus:\mathbf{C}\times \mathbf{C}\longrightarrow \mathbf{C}.$$
\end{enumerate}

\item The category $\mathbf{C}$ is monoidal and the tensor product
\[
\otimes : \mathbf{C} \times \mathbf{C} \longrightarrow \mathbf{C}
\]
is additive in each variable. More precisely, for all  
$X,Y,Z$ $\in$ $\Obj(\mathbf{C})$ and morphisms $f,g : X \longrightarrow Y$, $h,k : Y \longrightarrow Z$,
we have
\[
(f+g)\otimes h = f\otimes h + g\otimes h,
\qquad
f \otimes (h+k) = f\otimes h + f\otimes k .
\]
\end{enumerate}
\end{Definition}

Given a monoidal additive category $\mathbf{C}$, we will denote by
$${\bf{flip}}: {\mathbf{C}}\times{\mathbf{C}}\longrightarrow {\mathbf{C}}\times{\mathbf{C}}$$ the functor that sends $(X,Y)$ $\in$ $\Obj({\mathbf{C}}\times{\mathbf{C}})$ to $(Y,X)$ $\in$ $\Obj({\mathbf{C}}\times{\mathbf{C}})$. In this way, we have \cite{barcategories}.

\begin{Definition}
\label{b.3}
A bar category is a monoidal additive category $\mathbf{C}$
together with the following additional data:

\begin{enumerate}

\item A functor
\[
{\bf{-}}:\mathbf{C}\longrightarrow\mathbf{C},
\]
called the bar functor $$X\longmapsto \overline{X},\qquad (f:X\longrightarrow Y)\longmapsto (\overline{f}: \overline{X}\longrightarrow  \overline{Y}).  $$

\item A natural isomorphism $\mathrm{bb}$ between the identity functor and the double bar functor.

\item A natural isomorphism $\Xi^{\oplus}$ between ${\bf{-}}\circ \oplus$ and $\oplus \circ ({\bf{-}}\times {\bf{-}})$.

\item A natural isomorphism $\Xi^{\otimes}$ between ${\bf{-}}\circ \otimes$ and $\otimes \circ ({\bf{-}}\times {\bf{-}})\circ {\bf{flip}}$.

\item An isomorphism
\[
\star:\mathbbm{1}_{\mathbf{C}}\longrightarrow \overline{\mathbbm{1}_{\mathbf{C}}}.
\]

\end{enumerate}

These data satisfy the following coherence conditions:

\begin{enumerate}

\item The double bar functor is compatible with the bar functor,
\[
\overline{\mathrm{bb}_X} = \mathrm{bb}_{\overline{X}}.
\]

\item If $a$ is the associator, then $$a_{\overline{Z},\overline{Y},\overline{X}}\circ (\Xi^\otimes_{Y,Z}\otimes \id)\circ \Xi^\otimes_{X,Y\otimes Z}\circ \overline{a_{X,Y,Z}}=(\id\otimes \Xi^\otimes_{X,Y})\circ \Xi^\otimes_{X\otimes Y,Z} $$

\item The composition
\[
\overline{X}
\xrightarrow{\overline{r^{-1}_X}}
\overline{X\otimes \mathbbm{1}_{\mathbf{C}}}
\xrightarrow{\Xi^\otimes_{X,\mathbbm{1}_{\mathbf{C}}}}
\overline{\mathbbm{1}_{\mathbf{C}}}\otimes \overline{X} \xrightarrow{\star^{-1}\otimes \id}\mathbbm{1}_{\mathbf{C}}\otimes \overline{X}
\]
coincides with $l^{-1}_{\overline{X}}$; and
 \[
\overline{X}
\xrightarrow{\overline{l^{-1}_X}}
\overline{\mathbbm{1}_{\mathbf{C}}\otimes X}
\xrightarrow{\Xi^\otimes_{\mathbbm{1}_{\mathbf{C}},X}}
\overline{X}\otimes \overline{\mathbbm{1}_{\mathbf{C}}}  \xrightarrow{\id\otimes \star^{-1}} \overline{X}\otimes \mathbbm{1}_{\mathbf{C}}
\]
coincides with $r^{-1}_{\overline{X}}$.
\end{enumerate}
\end{Definition}

\begin{Definition}
\label{additivefunctor}
Let $\mathbf{C}$ and $\mathbf{D}$ be additive categories. A contravariant functor
\[
\F:\mathbf{C}\longrightarrow\mathbf{D}
\]
is called an additive functor if for every pair of objects 
$X,Y$ $\in$ $\Obj(\mathbf{C})$ the induced map
\[
\F_{X,Y}:\Mor(X,Y)
\longrightarrow
\Mor(F(Y),F(X)), \qquad f\longrightarrow \F(f)
\]
is a morphism of abelian groups and there is a natural isomorphism between $$\F(X\oplus Y)\qquad \mbox{ and } \qquad \F(X)\oplus'\F(Y).$$ 
\end{Definition}

\begin{Definition}
    \label{b.5}
Let $(\mathbf{C},\otimes,\mathbbm{1}_{\mathbf{C}},a,l,r)$ and 
$(\mathbf{D},\otimes',\mathbbm{1}_{\mathbf{D}},a',l',r')$ be monoidal categories.

A monoidal contravariant functor 
\[
(\F,\phi_2,\phi_1):\mathbf{C}\longrightarrow\mathbf{D}
\]
consists of
\begin{enumerate}

\item A contravariant functor
\[
\F:\mathbf{C}\longrightarrow\mathbf{D}.
\]

\item A natural transformation
\[
\phi_2(X,Y):\F(X)\otimes' \F(Y)\longrightarrow \F(X\otimes Y)
\]
for all  $X,Y$ $\in$ $\Obj(\mathbf{C})$.

\item A morphism
\[
\phi_1:\mathbbm{1}_{\mathbf{D}}\longrightarrow \F(\mathbbm{1}_{\mathbf{C}})
\]
such that the following conditions hold:
\begin{enumerate}
\item  For all $X,Y,Z$ $\in$ $\Obj(\mathbf{C})$ $$\F(a_{X,Y,Z})\circ \phi_2(X,Y\otimes Z)\circ (\id\otimes'\phi_2(Y,Z))\circ a'_{F(X),F(Y),F(Z)}=\phi_2(X\otimes Y,Z)\circ (\phi_2(X,Y)\otimes \id). $$
\item For all $X$ $\in$ $\Obj(\mathbf{C})$ $$\F(r_X)\circ r'_{\F(X)}=\phi_2(X,\mathbbm{1}_{\mathbf{C}})\circ (\id\otimes \phi_1) $$ and $$\F(l_X)\circ l'_{\F(X)}=\phi_2(\mathbbm{1}_{\mathbf{C}},X)\circ ( \phi_1\otimes \id).$$
\end{enumerate}
\end{enumerate}
Furthermore, we say that the functor $\F$ is strict if 
\[
\phi_2(X,Y) \qquad \mbox{ and }\qquad \phi_1
\]
are isomorphisms, for all objects $X, Y$ $\in$ $\mathbf{C}$.
\end{Definition}

Finally, we have

\begin{Definition}
\label{b.6}
Let $\mathbf{C}$ and $\mathbf{D}$ be two bar categories. A contravariant bar functor consists of a strict additive monoidal contravariant functor 
\[
\F:\mathbf{C}\longrightarrow\mathbf{D},
\]
and a natural isomorphism
\[
\mathrm{fb}_X:\overline{\F(X)}\longrightarrow \F(\overline{X})
\]
for every $X$ $\in$ $\Obj(\mathbf{C})$, such that the following compatibility conditions hold:

\begin{enumerate}
\item $\F(\star^{-1})\circ \phi_1 = \mathrm{bf}_{\mathbbm{1}_{\mathbf{C}}}\circ \overline{\phi_1}\circ \star'$.

\item For all $X$ $\in$ $\Obj(\mathbf{C})$
\[
\F(\mathrm{bb}^{-1}_X)
=
\mathrm{fb}_{\overline{X}}\circ \overline{\mathrm{fb}_X}\circ \mathrm{bb}'_{\F(X)} .
\]

\item For all objects $X,Y$ $\in$ $\Obj(\mathbf{C})$ $$ \phi_2(\overline{Y},\overline{X})\circ (\mathrm{fb}_Y\otimes \mathrm{fb}_X)\circ \Xi^{\otimes' }_{\F(X),\F(Y)}\circ\overline{\phi^{-1}_2(X,Y)}= \F(\Xi^{\otimes\,-1}_{X,Y})\circ \mathrm{fb}_{X\otimes Y}.$$
\end{enumerate}
\end{Definition}

\section{A $C^\ast$--algrebaic Generalization of $\A_\zeta$}

As discussed in the first section, by Theorem~\ref{theoremrep} and Remark~\ref{rema}, the underlying context of this paper is that of $C^\ast$--algebras, although all the developments and results are presented in purely algebraic terms, based on geometric considerations. Furthermore, in Section 3.2 of references \cite{sald2} it is proven that
\begin{equation}
    \label{c.1}
    \langle-,-\rangle: E^V\times E^V\longrightarrow B, \qquad (T_1,T_2)\longmapsto \sum_{i} T_1(e_i)\,T_2(e_2)^\ast
\end{equation}
is a non--degenerate Hermitian structure on $E^V$ as left $B$--module, where $\{ e_i\}$ is an orthonormal basis of $V$ (with respect to the inner product that makes $\delta^V$ unitary); and in Theorem 3.12 of \cite{sald2} it is proven that $E^V$ is a left pre--Hilbert $C^\ast$--module over the pre--$C^\ast$--algebra $B$  with respect to this Hermitian structure. Of course, if $B$ is a $C^\ast$--algebra, $E^V$ is a left Hilbert $C^\ast$--module (\cite{sald2}). In particular, in this situation, $\langle-,-\rangle$ is a $B$--valued inner product. 

In this way, the entire theory remains closely connected to the $C^\ast$--algebra framework, and in this short appendix we pass to that level, establishing a link between the theory of this paper and Durdevich's formulation of qpb's with previous studies on the subject \cite{11,br2,ss,ss1}.

Let $$\Delta_{\overline{P}}: \overline{P}\longrightarrow \overline{P}\otimes_{\mathrm{min}} \G$$ be a free $C^\ast$--dynamical system (\cite{ss,ss1}), where $(\overline{P},||\;\;||)$ is a $C^\ast$--algebra, $\G$ is a compact quantum group with coproduct $$\Delta:\G\longrightarrow \G\otimes_{\mathrm{min}}\G,$$ with $\otimes_{\mathrm{min}}$ denoting the minimal tensor product of $C^\ast$--algebras. 

Consider $H^\infty=(H,m,\mathbbm{1},\Delta,\epsilon,S,\ast)$ the canonical dense $\ast$--Hopf algebra of $\G$ (the linear span of the matrix coefficients of all finite--dimensional quantum $\G$--representations), and let us define by $$P=\{p\in \overline{P}\mid \Delta_{\overline{P}}(p)\,\in\,\overline{P}\otimes H\}.$$ Here, $\otimes$ is the algebraic tensor product.  According to \cite{11}, since $\Delta_{\overline{P}}$ is free, we have that $$\zeta=(P,B,\Delta_P:=\Delta_{\overline{P}}|_{P})$$ is a quantum principal $H$--bundle over $$B:=\{b\in P\mid \Delta_P(b)=b\otimes \mathbbm{1} \}$$ and according to Lemma 1.2 of reference \cite{ss1}, $B$ is actually a non--degenerate $C^\ast$--subalgebra of $\overline{P}$. Therefore, Remark \ref{rema} holds and Theorem \ref{gen} is satisfied.

In reference \cite{ss1} is presented a $C^\ast$--algebraic generalization of the theory of Section 2.5. In fact, consider the category $$\mathbf{Rep}_\G$$ of finite--dimensional (unitary) \emph{quantum} representations of $\G$, and consider the category $$\mathbf{Corr}\,B$$ of $C^\ast$--correspondence on $B$. By a  $C^\ast$--correspondence on $B$, we mean a left Hilbert $C^\ast$--module over $B$ together with a non--degenerate right action of $B$ on it. Then, the spectral association functor $$\A^{C^\ast}_{\zeta}:\mathbf{Rep}_\G\longrightarrow \mathbf{Corr}\,B $$ is defined as the weak unitary monoidal contravariant functor such that in object is given by $$\A^{C^\ast}_\zeta(\delta^V)=\overline{P}\;\square_\G\,V,$$ where $\square_\G$ denotes the cotensor product of $\overline{P}$ and $V$, and $$\A^{C^\ast}_{\zeta}(T)=\id_{\overline{P}}\otimes_{\mathrm{min}}\,T $$ on morphisms. Furthermore, the $B$--valued inner product of $\overline{P}\;\square_\G\,V$ is given by $$\langle X,Y\rangle:=\sum_i x_i\,y^\ast_i$$ if $X=\displaystyle\sum_i x_i\otimes_{\mathrm{min}} e_i$, $X=\displaystyle\sum_i y_i\otimes_{\mathrm{min}} e_i$,  where $\{ e_i\}$ is an orthonormal basis of $V$ (with respect to the inner product that makes $\delta^V$ unitary) \cite{ss1}.

On the hand, it is easy to see that the isomorphism of Proposition 6.1 of reference \cite{br2} is also an isomorphisms of $C^\ast$--correspondence over $B$ between $$\overline{P}\;\square_\G\,V^\#\qquad  \mbox{ and }\qquad  \Mor(\delta^V,\Delta_{\overline{P}})$$ with the $B$--valued inner product of equation (\ref{c.1}) and hence, the spectral association functor can be defined in terms of interwiner maps $\Mor(\delta^V,\Delta_{\overline{P}})$ as in Section 2.5.

It is worth noting that the category $\mathbf{Rep}_\mathcal{G}$ is the $C^\ast$--algebraic generalization of $\Rep_H$, while $\mathbf{Corr}\,B$ is the $C^\ast$--algebraic generalization of $\VB_B$. Moreover, the concept of  weak unitary monoidal contravariant functor is the $C^\ast$--algebraic generalization of the notion of  bar contravariant functor. 

In this way, the main result of~\cite{ss1}, stated in Theorem~2.2, shows that if $\mathcal{G}$ is a reduced compact quantum group, then there is a bijection between isomorphism classes $(\overline{P},B,\Delta_{\overline{P}})$ and isomorphism classes of weak unitary monoidal contravariant functors  $$\F:\mathbf{Rep}_\G \longrightarrow \mathbf{Corr}\,B,$$ which is  the $C^\ast$--algebraic generalization (at degree~$0$ and of course, without quantum connections,  as in Section 2.5), of the main result of this paper, i.e.,  Theorem~\ref{teo2}.

Finally, it is worth mentioning that, since Remark~\ref{rema} holds, $\zeta$ is strong (in the sense of Section~5.2 of~\cite{libro}), in accordance with Proposition~\ref{horizontalstrong}. Therefore, for every differential calculus on $\zeta$ and every qpc $\omega$, the linear map $$\nabla^\omega_V:=\Upsilon_V\circ D^\omega: E^V\longrightarrow \Omega^\bullet(B)\otimes_B E^V $$ satisfies the left Leibniz rule for every $\delta^V\in \Obj(\Rep_H)$, as verified in Section~3.4.1. Moreover, in light of Theorem 3.16 of reference \cite{sald2}, $\nabla^\omega_V$ is Hermitian. In other words, 
\begin{equation}
    \label{c.2}
    \langle \nabla^\omega_V(T_1),T_2\rangle+\langle T_1,\nabla^\omega_V(T)\rangle =d \langle T_1,T_2\rangle
\end{equation}
for all $T_1$, $T_2$ $\in$ $E^V$, where we have extended the $B$--valued inner product of equation (\ref{c.1}) to 
\begin{equation*}
\begin{aligned}
\langle-,-\rangle:(\Omega^\bullet(B)\otimes_B E^V)\times (\Omega^\bullet(B)\otimes_B E^V)\longrightarrow  \Omega^\bullet(B)
\end{aligned}
\end{equation*}
by means of  $$\langle \mu_1\otimes_B T_1,\mu_2\otimes_B T_2\rangle=\mu_1 \,\langle T_1,T_2\rangle\,\mu^\ast_2.$$ It is worth noticing that equation (\ref{c.2}) holds for every differential calculus on $\zeta$, every qpc $\omega$ (not necessarily regular nor multiplicative) and every $\delta^V\in \Obj(\Rep_H)$.

\end{document}